\documentclass[francais]{ensmath}

\usepackage{epsfig}
\usepackage[mathscr]{eucal}
\usepackage{graphpap}

\usepackage{amsmath, amssymb}
\usepackage{epstopdf}
\newtheorem{question}{Question}

\makeatletter 
\renewcommand{\thefigure}{\ifnum \c@section>\z@ \thesection.\fi\@arabic\c@figure}
\@addtoreset{figure}{section}

\@addtoreset{equation}{section}
\renewcommand\thetable{\thesection.\arabic{table}}
\@addtoreset{table}{section}
\makeatother

\def\:{{\nobreak\hspace{0.5ex}:}}
\def\;{{\nobreak\hspace{0.5ex};}}
\def\?{{\nobreak\hspace{0.5ex}?}}
\font\cyr=wncyr10
\def\og{{\cyr\char'074}}
\def\fg{{\cyr\char'076}}

\newcommand{\abs}[1]{\left\vert#1\right\vert}
\newcommand{\arbre}{\mathcal{A}}

\newcommand{\braid}{{B}}

\newcommand{\Compl}{\mathbb C}
\newcommand{\Cl}{\Phi}

\newcommand{\code}{\mathrm{c}}
\newcommand{\Conj}{\Psi}
\newcommand\cred{{\mathrm{cr}}}

\newcommand{\disque}{\mathbb D}
\newcommand{\domfond}{\mathcal D}
\newcommand\dec{\text{d\'ec}}

\newcommand{\feuilles}{\mathcal{F}}

\renewcommand\ge{\geqslant}
\newcommand{\GLZ}{\SLpm}
\newcommand{\geod}{\lambda}

\newcommand{\Hy}{\mathcal H}
\newcommand{\HOM}{P}

\newcommand{\Ideal}{\mathcal I}

\renewcommand\le{\leqslant}
\newcommand{\longueur}{\ell}
\newcommand{\lid}{\left\langle}
\newcommand{\lbid}{\left\{}

\newcommand{\mama}[4]{\left( \begin{matrix} #1 & #2 \\ #3 & #4 \end{matrix} \right)}

\newcommand{\MMZ}{\mathrm{M}_{2}(\relatifs)}
\newcommand{\module}[1]{\left\vert#1\right\vert}
\newcommand{\monj}{j}
\newcommand{\monome}{p}
\newcommand{\mot}{w}
\newcommand{\Murasomme}{\mathbin{\mathtt{\#}}}

\newcommand{\ordre}{o}
\newcommand{\grandordre}{\mathcal O}

\newcommand{\pgcd}{\text{pgcd}}
\newcommand{\PSL}{\mathrm{PSL}_2(\relatifs)}
\newcommand{\PSLR}{\mathrm{PSL}_2(\reels)}
\newcommand{\poids}{\pi}

\newcommand\pr{{\mathrm{pr}}}

\newcommand{\rationnels}{\mathbb Q}

\newcommand{\rid}{\right\rangle}
\newcommand{\rbid}{\right\}}
\newcommand{\reels}{\mathbb R}
\newcommand{\resslr}{h}
\newcommand{\resmod}{f}
\newcommand{\rond}{\bigcirc}

\newcommand{\SLR}{\mathrm{SL}_2(\reels)}
\newcommand{\SLZ}{\mathrm{SL}_2(\relatifs)}
\newcommand{\SLpm}{\mathrm{GL}_2(\relatifs)}
\newcommand{\Sph}{\mathbb S}

\newcommand{\surf}{\Sigma}
\newcommand{\surfmod}{\Sigma_{mod}}

\newcommand{\tangente}{tangente}
\newcommand{\tore}{T}
\newcommand{\trip}{pas}
\newcommand{\trefoil}{
	\begin{picture}(3.5,0)
	\put(-1,-0.5){\includegraphics[height=3mm]{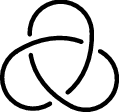}}
	\end{picture}}

\newcommand{\universel}{U}
\newcommand{\un}{T^{1}}
\newcommand{\USmod}{\un\Sigma_{mod}}
\newcommand{\USmodx}{\un_x\Sigma_{mod}}

\newcommand{\xx}{L}

\newcommand{\yy}{R}

\newcommand{\relatifs}{\mathbb Z}

\begin{document}

\title{\sc Les n{\oe}uds de Lorenz}

\author[Pierre Dehornoy]{Pierre Dehornoy\footnote{Unit\'e de math\'ematiques pures et appliqu\'ees, \'Ecole normale sup\'erieure de Lyon, 46 all\'ee d'Italie, 69364 Lyon Cedex 07, France, {\tt pierre.dehornoy@umpa.ens-lyon.fr}}}

\date{9 novembre 2009}

\maketitle 

\begin{abstract}
Cet article est un expos\'e de synth\`ese sur les n\oe uds de Lorenz. 
Nous d\'ecrivons la construction initiale, d\'emontrons de fa\c{c}on pr\'ecise un certain nombre de propri\'et\'es classiques, en particulier le fait que la cl\^oture d'une tresse positive est un n\oe ud fibr\'e, et explicitons la correspondance de Ghys entre n\oe uds modulaires et n\oe uds de Lorenz. 
Par ailleurs, nous montrons deux propri\'et\'es nouvelles, \`a savoir que, suivant la correspondance de Ghys, les images des orbites triviales forment un sous-groupe du groupe des classes, et, d'autre part, que l'image r\'eciproque d'un \'el\'ement du groupe des classes et de son inverse sont des orbites isotopes.
\end{abstract}

\begin{likeabstract}[Abstract]
This article is a survey on Lorenz knots. We describe the original construction, prove several classical properties, in particular the fact that the closure of a positive braid is a fibered knot, and describe Ghys'correspondance between modular knots and Lorenz knots. We also prove two new properties, namely that following Ghys' correspondance, the images of trivial orbits of the Lorenz flow form a subgroup of the class group, and that the reverse images of an element in the class group and of its inverse are isotopic orbits.
\end{likeabstract}


Les n{\oe}uds de Lorenz sont les n{\oe}uds qui peuvent \^etre r\'ealis\'es comme orbites p\'eriodiques du flot g\'eom\'etrique de Lorenz. 
Tous les n{\oe}uds ne sont pas des n{\oe}uds de Lorenz\: par exemple, le n{\oe}ud trivial et le n{\oe}ud de tr\`efle en sont, mais pas le n{\oe}ud de huit\; plus g\'en\'eralement, on ne compte que vingt n{\oe}uds de Lorenz parmi les n{\oe}uds ayant au plus seize croisements. 
Pourtant, la famille des n{\oe}uds de Lorenz est vaste, puisque tous les n{\oe}uds toriques et tous les n{\oe}uds alg\'ebriques sont des n{\oe}uds de Lorenz. 
Introduits par J.\,Birman et R.\,Williams en 1983, ces n{\oe}uds ont fait l'objet de travaux assez nombreux, car ils poss\`edent de riches propri\'et\'es et de multiples facettes. 
On sait depuis les travaux de R.\,Ghrist qu'il existe un flot non singulier dans~$\reels^3$ tel que tout n{\oe}ud est r\'ealis\'e comme orbite p\'eriodique dudit flot; mais le r\'esultat ne fournit en g\'en\'eral aucun contr\^ole des n{\oe}uds ainsi engendr\'es. 
A l'oppos\'e, la combinatoire sp\'ecifique du flot de Lorenz, qui s'exprime par des liens avec les mots de Lyndon et les diagrammes de Young, donne acc\`es pour les n{\oe}uds de Lorenz \`a des param\`etres topologiques comme le genre ou l'indice de tresse qui sont difficiles \`a calculer dans le cas g\'en\'eral. 
De la sorte, la famille des n{\oe}uds de Lorenz appara\^it comme une \'etape naturelle dans une \'etude des n{\oe}uds partant des plus simples et allant progressivement vers davantage de complexit\'e\: plus riche que celle des n{\oe}uds toriques ou celle des n{\oe}uds alg\'ebriques, la famille de n{\oe}uds de Lorenz est n\'eanmoins assez petite pour qu'une description effective des propri\'et\'es de ses \'el\'ements soit possible. 
R\'ecemment, l'int\'er\^et port\'e \`a ces n{\oe}uds a \'et\'e renforc\'e par la d\'ecouverte par \'E.\,Ghys d'une correspondance naturelle entre les n{\oe}uds de Lorenz et les n{\oe}uds dits modulaires, qui apparaissent comme orbites p\'eriodiques du flot g\'eod\'esique sur la surface modulaire, un certain orbifold li\'e \`a l'arithm\'etique des corps quadratiques et \`a la g\'eom\'etrie hyperbolique.

Ce que nous faisons dans ce texte, c'est principalement de passer en revue les propri\'et\'es des n\oe uds de Lorenz connues \`a ce jour. 
La plupart des d\'emonstrations sont seulement esquiss\'ees, \`a l'exception de quelques-unes qui, dans la litt\'erature, apparaissent de fa\c{c}on sch\'ematique ou dans un contexte diff\'erent, et que nous pr\'esentons ici de fa\c{c}on d\'etaill\'ee. 
Par ailleurs, nous \'etablissons deux r\'esultats nouveaux dans la lign\'ee de la correspondance de Ghys, \`a savoir, d'une part, que, pour tout ordre d'un anneau d'entiers quadratiques, les classes d'id\'eaux dont l'image par la correspondance de Ghys est un n{\oe}ud trivial forment un sous-groupe du groupe des classes de la forme $(\relatifs{/}2\relatifs)^d$ dont la multiplication peut \^etre d\'ecrite de fa\c{c}on explicite, et, d'autre part, que deux \'el\'ements inverses l'un de l'autre dans le groupe des classes ont pour image le m\^eme n{\oe}ud par la correspondance de Ghys.
 
Le plan de l'article est le suivant. 
Dans une premi\`ere partie, on effectue un bref rappel historique sur l'origine des n\oe uds de Lorenz.
Puis on introduit les n{\oe}uds de Lorenz \`a partir du patron de Lorenz, qui est une surface branch\'ee de~$\reels^3$ li\'ee \`a la dynamique d'un certain syst\`eme diff\'erentiel chaotique. 
Un interm\`ede est alors d\'edi\'e \`a le th\'eorie des surfaces branch\'ees. 
On d\'eduit ensuite un codage des n{\oe}uds de Lorenz \`a l'aide des mots de Lyndon, qui sont des mots particuliers sur un alphabet \`a deux lettres. 
Ce codage permet d'\'enum\'erer syst\'ematiquement tous les n{\oe}uds de Lorenz, et il est \`a la base de toute l'\'etude ult\'erieure. 
Notre exposition ici est bas\'ee sur l'article original de J.\,Birman et R.\,Williams~\cite{A:B-W}, \`a l'exception de la section~\ref{S:Patrons}, et des sections~\ref{S:Stabilisation} et~\ref{S:Young} o\`u nous d\'eveloppons une nouvelle variante moins redondante du codage initial \`a l'aide de diagrammes de Young.

La deuxi\`eme partie est consacr\'ee \`a l'\'etendue de la famille des n{\oe}uds de Lorenz et \`a leurs principales propri\'et\'es topologiques. 
On montre que tout n{\oe}ud de Lorenz est premier. 
On y explique que tous les n{\oe}uds toriques et tous les n{\oe}uds alg\'ebriques sont des n{\oe}uds de Lorenz, le cas des n{\oe}uds alg\'ebriques se d\'eduisant d'une \'etude plus g\'en\'erale des n{\oe}uds de Lorenz qui sont des n{\oe}uds satellites. 
On montre ensuite que tout n\oe ud de Lorenz est un n{\oe}ud fibr\'e, et on en d\'eduit une formule pour le genre des n\oe uds de Lorenz, ainsi qu'un crit\`ere permettant d'exclure certain n\oe uds de la liste des n\oe uds de Lorenz. 
Les r\'esultats de cette partie apparaissent principalement dans les articles~\cite{A:B-W},~\cite{A:ElRifai2}, et~\cite{A:Williams}. 
Pour ce qui est du caract\`ere fibr\'e, il est d\'eduit dans~\cite{A:B-W} de r\'esultats g\'en\'eraux de Stallings~\cite{A:Stallings}. 
Il peut \'egalement \^etre d\'eduit de travaux de K.\,Murasugi~\cite{A:Murasugi} et D.\,Gabai~\cite{A:Gabai1, A:Gabai2} sur les fibrations. 
Dans le cas de n{\oe}uds qui sont cl\^otures de tresses positives, l'argument  est sp\'ecialement simple et visuel. 
Comme il semble n'avoir jamais \'et\'e r\'edig\'e, tout au moins de fa\c{c}on concise, nous le d\'etaillons~ici.

Dans la troisi\`eme partie, nous poursuivons l'\'etude topologique des n{\oe}uds de Lorenz. 
On introduit la famille des tresses de Birman-Williams, qui fournit une nouvelle fa\c{c}on d'exprimer un n{\oe}ud de Lorenz comme cl\^oture d'une tresse positive, et on en d\'eduit une d\'etermination explicite de l'indice de tresse d'un n{\oe}ud de Lorenz comme nombre de brins d'une tresse de Birman--Williams associ\'ee. 
Le r\'esultat original est d\^u \`a J.\,Franks, H.\,Morton et R.\,Williams~\cite{A:F-W}. 
Nous proposons ici une version compl\`ete de la d\'emonstration bas\'ee sur l'introduction d'une variante du polyn\^ome HOMFLY.

La quatri\`eme partie est consacr\'ee \`a la correspondance de Ghys entre les n\oe uds de Lorenz et les orbites p\'eriodiques du flot g\'eod\'esique sur la surface modulaire, lesquelles sont intimement li\'ees \`a des objets tr\`es usuels de la th\'eorie des nombres depuis Gauss, \`a savoir les classes de conjugaison de matrices dans~$\SLZ$, les classes de formes quadratiques \`a coefficients entiers et les classes d'id\'eaux dans des corps quadratiques. 
Pour expliquer ce lien entre topologie et arithm\'etique, nous partons de la correspondance naturelle entre les classes d'id\'eaux dans certains corps quadratiques et les classes de conjugaison dans~$\SLZ$, qu'on peut caract\'eriser en termes de mots de Lyndon. 
On en d\'eduit le lien entre les n\oe uds modulaires et les n\oe uds de Lorenz qui constitue la correspondance de Ghys. 
Notre exposition dans cette partie est bas\'ee sur~\cite{A:Cohn} et \cite{A:Ghys3}. 
Par ailleurs, c'est l\`a que nous \'etablissons les deux r\'esultats nouveaux mentionn\'es plus haut concernant les orbites repr\'esentant l'\'el\'ement neutre du groupe des classes (th\'eor\`eme~\ref{T:groupetrivial}) et les orbites repr\'esentant des \'el\'ements inverses l'un de l'autre (proposition~\ref{T:inverse}).

Enfin, dans une br\`eve cinqui\`eme partie, nous mentionnons quelques questions ouvertes mettant en jeu les n{\oe}uds de Lorenz. 
L'une des plus fascinantes concerne la possibilit\'e de d\'efinir directement sur les n{\oe}uds de Lorenz une multiplication qui serait la contrepartie de celle des classes d'id\'eaux d'un ordre d'un corps quadratique.

{\bf Remerciements}. Je remercie vivement \'Etienne Ghys pour m'avoir init\'e aux n\oe uds de Lorenz, m'avoir guid\'e pendant ce travail, et pour de passionnantes conversations. Je remercie \'egalement Hugh Morton pour avoir partag\'e son exp\'erience sur les n\oe uds de Lorenz et r\'epondu \`a plusieurs de mes questions. Je remercie Maxime Bourrigan et les deux referees d'une version pr\'ec\'edente pour de nombreuses remarques, ayant en particulier permis de lever certaines impr\'ecisions.


\section{Description combinatoire des n{\oe}uds de Lorenz}
\label{S:Combinatoire}

Cette premi\`ere partie introduit les n{\oe}uds de Lorenz et leur codage par des mots de Lyndon et des diagrammes de Young. 
Dans la section~\ref{S:Systeme}, nous rappelons le contexte du syst\`eme diff\'erentiel chaotique de Lorenz, puis, dans la section~\ref{S:Patron}, nous introduisons le patron de Lorenz qui permet de mod\'eliser g\'eom\'etriquement le comportement des solutions du syst\`eme de Lorenz. 
La section~\ref{S:Patrons} est un interm\`ede introduisant la th\'eorie des patrons qui a \'et\'e d\'eriv\'ee de l'\'etude des n\oe uds de Lorenz.
La section~\ref{S:Codage} est consacr\'ee au codage des orbites p\'eriodiques du flot de Lorenz, donc des n{\oe}uds de Lorenz, \`a l'aide de mots de Lyndon. 
Dans la section~\ref{S:TresseLorenz}, nous montrons que tout n{\oe}ud de Lorenz est cl\^oture d'une tresse de permutation particuli\`ere. 
Dans la section~\ref{S:Stabilisation}, nous montrons comment r\'eduire la redondance du codage en introduisant les mots de Lyndon dits minimaux, puis, dans la section~\ref{S:Young}, nous passons au langage des diagrammes de Young.


\subsection{Le syst\`eme de Lorenz, Guckenheimer et Williams}
\label{S:Systeme}

En 1963, le m\'et\'eorologiste E.\,N.\,Lorenz a exhib\'e un syst\`eme dynamique aux propri\'et\'es remarquables~\cite{A:Lorenz}. 
Issu d'un mod\`ele simplifi\'e de convection atmosph\'erique, il est d\'ecrit par les \'equations diff\'erentielles
	\begin{equation}
	\label{E:Lorenz}
	\dot x=-10x+10y, \quad \dot y=rx-y-xz, \quad \dot z=-8/3z+xy, 
	\end{equation}
o\`u $r$ est un param\`etre r\'{e}el proche de~24. 
Ce syst\`eme induit un flot d\'eterministe sur~$\reels^3$\: le pass\'e et le futur d'un point de~$\reels^3$ sont d\'etermin\'es par sa seule position. 
Par contre, \eqref{E:Lorenz} est un syst\`eme chaotique, au sens o\`u une petite perturbation du point de d\'epart change l'allure globale de l'orbite qui en est issue. 
Bien avant Lorenz, l'existence de syst\`emes chaotiques avait \'et\'e observ\'ee par H.\,Poincar\'e dans le cadre du probl\`eme des trois corps~\cite{A:PoincarŽ} et par J.\,Hadamard pour le flot g\'eod\'esique sur des surfaces \`a courbure n\'egative~\cite{A:Hadamard}. 
Pourtant, l'exemple de Lorenz, peut-\^etre \`a cause de son origine m\'et\'eorologique, ou peut-\^etre \`a cause de la simplicit\'e des \'equations qui le d\'ecrivent, a suscit\'e un grand engouement, se trouvant m\^eme \`a l'origine du populaire \og\,effet papillon\,\fg~\cite{A:Lorenz2}.

\begin{figure}[ht]
	\begin{center}
	\begin{picture}(120,50)(0,0)
	\put(30,-27){\includegraphics*[width=0.45\textwidth]{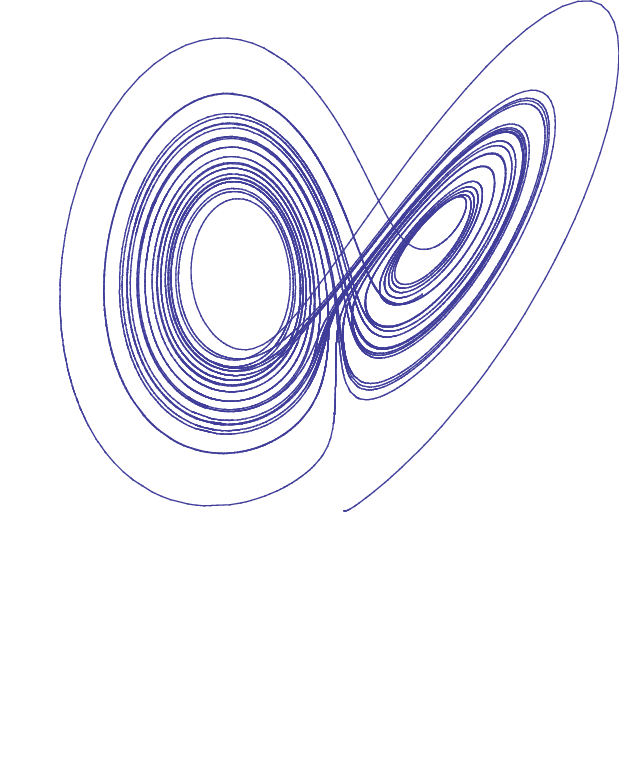}}
	\end{picture}
	\end{center}
	\caption{\small \sf Une simulation num\'erique des \'equations de Lorenz \eqref{E:Lorenz}. 
	Partant d'un point g\'en\'erique, le flot s'accumule sur une surface et tourne autour de deux points critiques.} 
	\label{F:LorenzNum}
\end{figure}

Quand on effectue des simulations num\'eriques des \'equations de Lorenz~\eqref{E:Lorenz}, on observe le ph\'enom\`ene suivant\: l'orbite issue de presque tout point semble plonger tr\`es vite vers une surface branch\'ee, puis, une fois au voisinage de celle-ci, l'orbite d\'ecrit une suite de tours autour de deux points critiques, dans un ordre apparemment al\'eatoire (figure~\ref{F:LorenzNum}). 
Ayant observ\'e ce ph\'enom\`eme, J.\,Guckenheimer a sugg\'er\'e un mod\`ele g\'eom\'etrique simple pour d\'ecrire les orbites de ce syst\`eme~\cite{A:Guckenheimer}, \`a savoir une surface de~$\reels^3$ munie d'un semi-flot imitant la dynamique observ\'ee num\'eriquement et  le long de laquelle s'enroulent les orbites. 
Ce mod\`ele de Guckenheimer a ensuite \'et\'e \'etudi\'e et enrichi dans des travaux communs avec R.\,Williams~\cite{A:G-W, A:Williams2}, o\`u il est notamment \'etabli que le mod\`ele est {\it persistant}, c'est-\`a-dire qu'une perturbation de classe $C^0$ du champ de vecteurs m\`ene \`a un syst\`eme du m\^eme type que celui de Lorenz. 
La corr\'elation entre les mod\`eles de Lorenz et de Guckenheimer, et en particulier le fait que les orbites de Lorenz s'accumulent effectivement le long du patron \'etudi\'e par Guckenheimer et Williams, n'a en fait jamais \'et\'e compl\`etement prouv\'ee. En 1999, W. Tucker~\cite{A:Tucker1, A:Tucker2} a montré que les orbites du flot de Lorenz s'accumulent le long d'un patron, qui est du type de celui prédit par Lorenz. Cependant, rien ne garantit que ce patron contienne autant d'orbites que celui de Guckenheimer-Williams. Ce dernier est donc une sorte de sur-patron du flot de Lorenz, c'est-à-dire qu'il décrit toutes les orbites du flot de Lorenz, et peut-être un peu plus\footnote{Avec le vocabulaire de la définition~\ref{D:code}, il n'est pas garanti que toutes les suites en~$L$ et~$R$ apparaissent comme codes d'orbites du flot de Lorenz, voir aussi la remarque~\ref{R:PremierRetour}}.

Le flot g\'eom\'etrique de Guckenheimer et Williams contient une infinit\'e d\'enombrable d'orbites p\'eriodiques. 
Au d\'ebut des ann\'ees 1980, J.\,Birman et R.\,Williams ont initi\'e l'\'etude de ces derni\`eres du point de vue de la th\'eorie des n\oe uds~\cite{A:B-W}.
Cette approche suit l'id\'ee, attribu\'ee \`a Poincar\'e, d'\'etudier un syst\`eme dynamique {\it via} ses orbites p\'eriodiques. 
De ce point de vue, il est trop t\^ot pour d\'ecider si cette t\^ache a \'et\'e r\'ealis\'ee. 
Par contre, nous allons voir que les n\oe uds de Lorenz sont une famille riche du point de vue de la th\'eorie des n\oe uds.

\subsection{Le patron de Lorenz}
\label{S:Patron}

Dans la suite de cet article, nous d\'elaissons les syst\`emes diff\'erentiels et repartons du mod\`ele g\'eom\'etrique sugg\'er\'e par Lorenz, Guckenheimer et Williams. 
On rappelle qu'\'etant donn\'e un champ de vecteurs~$X$ sur une surface~$S$, on appelle \emph{flot} de~$X$ l'application~$\Phi$ de~$\reels \times S$ dans~$S$ dont la d\'eriv\'ee en $t = 0$ co\"\i  ncide avec~$X$ et qui, pour tout point~$P$ de~$S$, v\'erifie 
$$\Phi(t_1 + t_2, P) = \Phi(t_1,\Phi(t_2,P)),$$
c'est-\`a-dire que $\Phi(t, P)$ d\'ecrit l'orbite partant de~$P$ le long du champ~$X$, lorsqu'une telle
application existe.  Un \emph{semi-flot} est une application semblable, mais d\'efinie seulement sur~$\reels_+ \times S$.

\begin{figure}[t!]
	\begin{center}
	\begin{picture}(120,98)(0,3)
	\put(0,0){\includegraphics*[scale=.7]{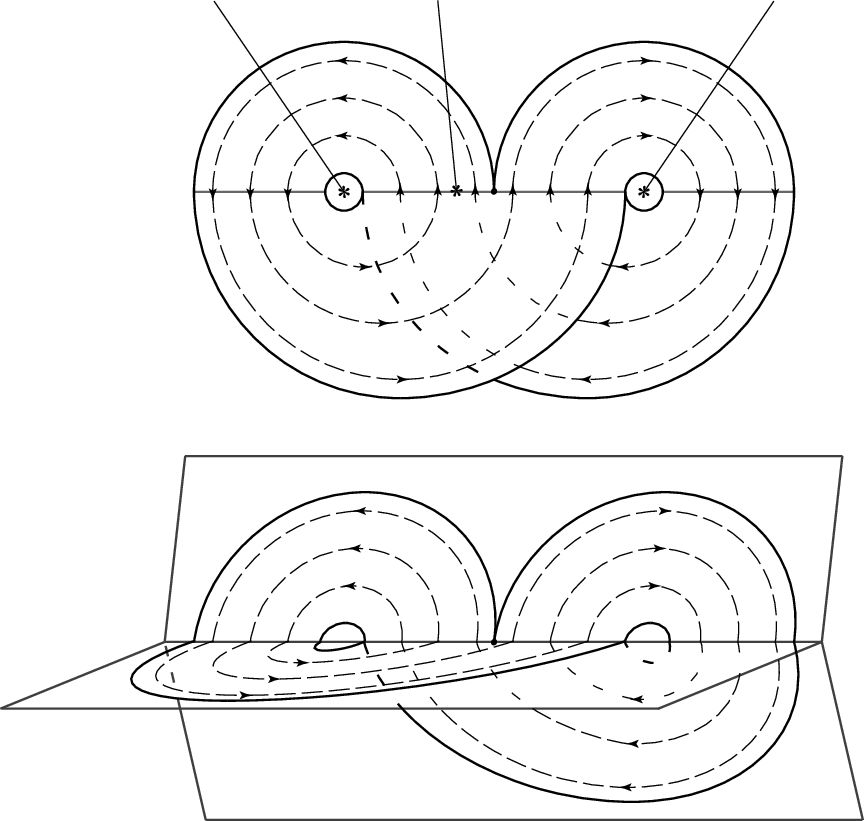}}
	\put(42.7,75.4){0} 
	\put(73,75.4){1} 
	\put(56.4,71.5){$M$}
	\put(48,98){axe $A$} 
	\put(11,98){point critique gauche}
	\put(78,98){point critique droit}
	\put(12,73.8){axe $A_1$}
	\put(95,73.8){axe $A_2$} 
	\put(21.5,85){$B_1$} 
	\put(91,86){$B_2$}
	\put(26,55){$B_3$} 
	\put(91,60){$B_4$}
	\end{picture}
	\end{center}
	\caption{\small \sf Le patron de Lorenz. Le dessin du haut repr\'{e}sente le patron projet\'e sur un plan horizontal. Les segments $A, A_1$ et $A_2$ sont les axes de recollement des bandes $B_1, \dots, B_4$, le segment $A$ \'{e}tant l'axe de branchement du patron. Le flot de Lorenz est figur\'{e} en pointill\'{e}s. L'abscisse du point $M$ est $1/2$, la surface s'y s\'{e}pare et le flot n'y est pas d\'{e}fini. Le dessin du bas montre le patron sous la forme d'un livre \`a trois pages.} 
	\label{F:patron}
\end{figure}

\begin{definition}
(Voir figure~\ref{F:patron}.) 
$(i)$ Le~{\it patron de Lorenz} est la surface branch\'ee de~$\reels^3$ obtenue comme suit. 
On part de quatre bandes rectangulaires ouvertes $B_1, \dots, B_4$, munies chacune d'un champ de vecteurs parall\`ele \`a deux c\^ot\'es oppos\'es. 
Ces bandes sont d\'eform\'ees et recoll\'ees le long de trois axes $A, A_1$ et~$A_2$ comme montr\'e sur la figure. 
La surface branch\'ee obtenue est munie du champ de vecteurs obtenu en recollant les champs des quatre bandes. 

$(ii)$ On appelle \emph{semi-flot de Lorenz} le semi-flot engendr\'e par le champ de vecteurs du patron de Lorenz. 
L'axe~$A$ est identifi\'e au segment $]0,1[$ de sorte  que l'abscisse du point~$M$ o\`u la surface se s\'epare soit~$1/2$, et que l'application de premier retour sur l'axe $A$ en suivant le flot, qu'on note~$\pr$, soit identifi\'ee \`a la fonction $t\mapsto 2t\mod 1$. 
\end{definition}

Par construction, le patron de Lorenz est une surface ouverte, branch\'{e}e le long de~$A$. 
Le champ de vecteurs est d\'efini en tout point, sauf au point~$M$ o\`{u} la surface se s\'{e}pare. 
Par cons\'equent, le semi-flot est d\'efini sur toute la surface priv\'ee des pr\'eorbites du point~$M$. 
En particulier, il n'est pas d\'efini pour les points de~$A$ dont l'abscisse est de la forme $m/2^n$ pour $m$ et $n$ entiers. 
En suivant le semi-flot \`{a} partir d'un point~$P$ de~$A$ diff\'{e}rent de $M$, on recoupe~$A$ apr\`es un temps fini. 
Par contre, le pass\'e n'est pas d\'efini\: le long de~$A$, chaque point a deux pr\'eorbites, on a donc bien affaire \`a un semi-flot. 
N\'eanmoins, par habitude autant que par commodit\'e, nous utilisons ici l'appellation {\it flot}.

Dans la suite, on s'int\'eresse aux orbites du flot de Lorenz, et particuli\`erement aux orbites \emph{p\'eriodiques}. 
Consid\'erons par exemple l'orbite issue du point~$A$ d'abscisse~$1/7$. 
En suivant le semi-flot, elle fait le tour de la boucle gauche du patron, passe par le point $2/7$, suit la boucle gauche, arrive au point $4/7$, suit la boucle droite, et, enfin, revient au point~$1/7$. 
On a donc une orbite p\'eriodique qui coupe trois fois l'axe du patron. 
Sur la figure~\ref{F:patron} sont repr\'{e}sent\'{e}es les orbites passant par les points d'abscisses~$1/7$ et~$3/7$. 
On verra que le flot de Lorenz contient une infinit\'{e} d'orbites p\'{e}riodiques. 
On appelle {\it p\'{e}riode} d'une telle orbite le nombre d'intersections de celle-ci avec l'axe du patron. 
Ce qui nous int\'eresse ici est le fait que de telles orbites p\'eriodiques peuvent former des n\oe uds non triviaux.

\begin{definition}
\label{D:noeudLorenz}
Un \emph{n\oe ud} est une classe d'isotopie de plongement du cercle $\Sph^1$ orient\'e dans l'espace~$\reels^3$ orient\'e. Un \emph{entrelacs} est une classe d'isotopie de plongements de plusieurs cercles dans l'espace
\footnote{On troque parfois l'espace ambiant $\reels^3$ contre son compactifi\'e d'Alexandroff, la sph\`ere~$\Sph^3$. 
Ce sera notamment le cas dans les sections~\ref{S:Topologie} et~\ref{S:Arithmetique} de ce survol. 
Ceci ne change pas les classes d'isotopies consid\'er\'ees.}.

Un n{\oe}ud $K$ est appel\'{e} {\it n{\oe}ud de Lorenz} s'il existe une orbite p\'{e}riodique $\gamma$ du flot de Lorenz  sur le patron de Lorenz repr\'esentant~$K$, c'est-\`a-dire dont~$K$ est la classe d'isotopie.
\end{definition}

\`A partir d'un n\oe ud~$K$, on obtient naturellement le n\oe ud \emph{oppos\'e} $\tilde K$ en renversant l'orientation de~$\Sph^1$, et deux n\oe uds \emph{miroirs} par sym\'etrie de $K$ et $\tilde K$ par rapport \`a un plan.
En g\'en\'eral, ces quatre n\oe uds ne sont pas isotopes. 
Cependant, nous verrons que tout n\oe ud de Lorenz est isotope \`a son oppos\'e (proposition~\ref{T:inverse}). 
Par contre, leurs images-miroir ne sont pas des n\oe uds de Lorenz (corollaire~\ref{T:nonamphichiraux}), en particulier, elles ne leurs sont pas isotopes.


\subsection{Interm\`ede: la th\'eorie des patrons et le patron universel de Ghrist}
\label{S:Patrons}

Apr\`es l'introduction du patron g\'eom\'etrique de Lorenz par Guckenheimer, et l'\'etude de ses orbites p\'eriodiques initi\'ee par Birman et Williams, une {\it th\'eorie des patrons} ({\it templates} en anglais) a \'et\'e lanc\'ee, en tant que branche de la th\'eorie des n\oe uds, en particulier par R.\,Ghrist, Ph.\,Holmes et M.\,Sullivan~\cite{A:G-H-S}.
Le but est de relier des propri\'et\'es des orbites p\'eriodiques d'un patron \`a la seule forme du patron et \`a son plongement dans~$\reels^3$. 
Par exemple, un analogue du th\'eor\`eme d'Alexander selon lequel tout entrelacs peut \^etre r\'ealis\'e comme cl\^oture d'une tresse existe pour les patrons, ainsi que des analogues des mouvements de Reidemeister.

\begin{figure}[htb]
	\begin{center}
	\includegraphics*[width=0.6\textwidth]{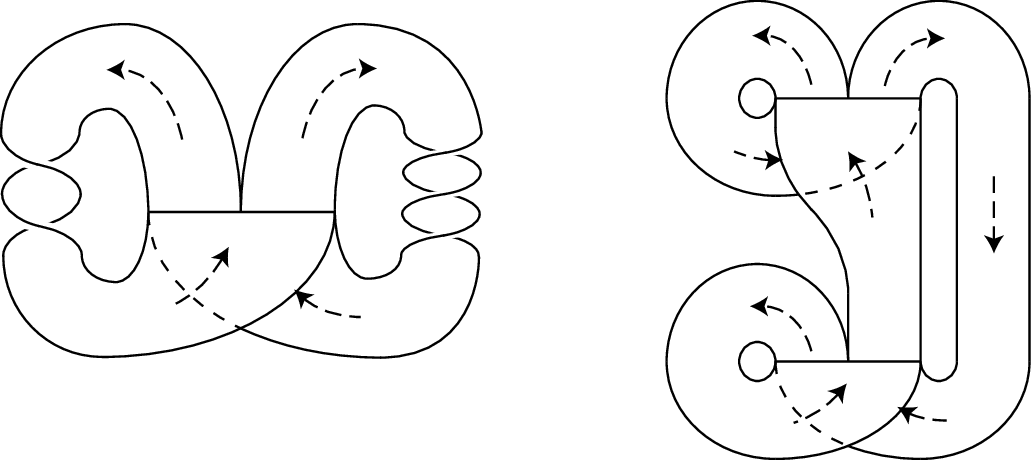}
	\end{center}
	\caption{\small \sf \`A gauche, un patron de Lorenz g\'en\'eralis\'e $L(-2, 3)$. \`A droite, le patron universel de Ghrist.}
	\label{F:patrongeneral}
\end{figure}

Notons que la classe d'isotopie d'une orbite p\'{e}riodique du flot g\'eom\'etrique de Lorenz d\'{e}pend du plongement du patron dans~$\reels^3$. 
On obtient une infinit\'e de plongements diff\'erents en ajoutant un nombre arbitraire de tours sur chacune des deux branches du patron. 
Plus pr\'ecis\'ement, on associe \`a chaque couple d'entiers relatifs $(p,q)$ un patron $L(p,q)$ obtenu \`a partir du patron de Lorenz en rajoutant $p$ demi-tours sur la branche gauche, et $q$ demi-tours sur la branche droite, comme sur la figure~\ref{F:patrongeneral}. 
\footnote{Si $p$ ou $q$ est impair, on obtient ainsi un patron non orientable.}
Les n\oe uds associ\'es \`a des orbites p\'eriodiques du flot d\'efini sur ces attracteurs de Lorenz g\'en\'eralis\'es d\'ependent fortement des entiers~$p$ et $q$. 
Cependant, \`a l'aide d'arguments visuels de d\'eformations et de d\'ecoupages simillaires \`a ceux de la figure~\ref{F:deformation}, M.\,Sullivan montre que pour tout $n$ entier, tout entrelacs r\'ealis\'e par une collection finie d'orbites p\'eriodiques de $L(0,n)$ l'est \'egalement par une collection d'orbites de $L(0,n-2)$~\cite{A:Sullivan1}. 
Ainsi, rajouter des tours n\'egatifs sur une des deux branches du patron de Lorenz ne peut qu'enrichir la famille des n\oe uds r\'ealis\'es.

En 1996, Ghrist~\cite{A:Ghrist} a exhib\'e un patron \emph{universel}~$\universel$, au sens o\`u tout entrelacs est r\'ealis\'e par un ensemble d'orbites p\'eriodiques de~$\universel$.
En fait, on retrouve $U$ comme \emph{sous-patron} de $L(0,-1)$ et de $L(0,-2)$, impliquant que les patrons $L(0,-n)$ sont universels pour tout entier $n$ strictement~positif.
Signalons que si tout n\oe ud est r\'ealis\'e par une orbite p\'eriodique de $U$, il n'est pas facile, pour un n\oe ud donn\'e, de pr\'edire la longueur minimale d'une telle orbite. Les bornes données par la  démonstration de Ghrist sont exponentielles, et ne permettent par exemple pas de savoir si le n\oe ud de huit est r\'ealis\'e par une orbite p\'eriodique de $U$ coupant moins de 11.358.338~fois un des axes de~$U$! \footnote{Il semblerait qu'une recherche informatique exhaustive exhibe un n\oe ud de huit coupant seulement 25 fois l'un des axes du patron, voir la th\`ese de Vadim Meleshuk~\cite{Meleshuk}.} 

Dans cet article, et plus g\'{e}n\'{e}ralement quand on parle de n{\oe}uds de Lorenz, on ne s'int\'{e}resse qu'au plongement $L(0,0)$ du patron de Lorenz montr\'{e} sur la figure~\ref{F:patron} qui, lui, n'est pas universel.


\subsection{Codage des orbites}
\label{S:Codage}

On cherche maintenant \`a d\'ecrire l'ensemble des orbites p\'eriodiques du flot  de Lorenz. 
Pour cela, on associe \`a toute orbite du flot de Lorenz un mot sur un alphabet \`a deux lettres. Restreinte aux orbites p\'eriodiques et aux mots dits de Lyndon
\footnote{Les mots de Lyndon apparaissent dans un tout autre contexte comme codes d'une base de l'alg\`ebre de Lie libre \cite{A:B-P,A:C-F-L,A:Lyndon}. 
Cependant, \cite{A:Reut} rappelle qu'ils \'enum\`erent naturellement les classes de formes quadratiques et les classes de conjugaison du groupe~$\PSL$. 
Comme nous allons montrer dans la section~\ref{S:Arithmetique} que les n\oe uds de Lorenz, \`a travers leur correspondance avec les n\oe uds dits modulaires, admettent une \'enum\'eration simillaire, nous utilisons d\`es maintenant l'appellation \emph{mot de~Lyndon}.} 
, la correspondance ainsi d\'efinie est bijective.

\begin{definition}
\label{D:code}
Soit $\gamma$ une orbite du flot de Lorenz. 
En suivant $\gamma$ \`a partir d'un point~$P$ de l'axe $A$ et en notant $\xx$ chaque fois que $\gamma$ coupe l'axe~$A$ apr\`es avoir fait un tour de l'attracteur gauche du patron, et $\yy$ chaque fois que $\gamma$ coupe l'axe~$A$ apr\`es un tour de l'attracteur droit, on associe \`a $\gamma$  un mot infini en les lettres $\xx$ et~$\yy$, appel\'e {\it code} de $\gamma$ et not\'e $\code(\gamma, P)$.
\end{definition}

Une orbite $\gamma$ est p\'eriodique de p\'eriode~$n$ si et seulement si le code de $\gamma$ est lui-m\^eme \emph{p\'eriodique de p\'eriode~$n$}, au sens o\`u il existe un (unique) mot $w$ de longueur~$n$ tel que $\code(\gamma,P)$ est le mot infini $wwww...$ 
On dit alors que $w$ est un {\it code r\'eduit} de~$\gamma$ \emph{relativement
\`a}~$P$, et on le note~$\cred(\gamma, P)$. 

Le code d'une orbite $\gamma$ d\'epend du choix du point de d\'epart $P$\: par exemple, si on part du point de premier retour $\pr(P)$, le code associ\'e $\code(\gamma, \pr(P))$ est le mot obtenu \`a partir de $\code(\gamma, P)$ en effa\c cant la premi\`ere lettre. 
Pour associer \`a chaque orbite p\'eriodique du flot de Lorenz un code r\'eduit ind\'ependant du point de d\'epart, on fixe un ordre lexicographique sur les mots en~$\xx$ et~$\yy$ et on choisit comme repr\'esentant distingu\'e celui des divers codes r\'eduits de l'orbite, qui est le plus petit vis-\`a-vis de cet ordre. 
Ceci revient \`a attribuer un rang aux lettres successives d'un code, correspondant \`a l'ordre selon lequel on rencontre les points d'intersection de l'orbite avec l'axe du patron de Lorenz quand on le parcourt de gauche \`a droite. 

\begin{definition}
\label{D:ordreLyndon}
Soit $w$ un mot de longueur $n$ sur l'alphabet~$\{\xx, \yy\}$. 
On appelle \emph{d\'ecal\'e} de $w$ le mot obtenu en d\'epla\c cant la premi\`ere lettre de $w$ en queue de mot. On le note $\dec(w)$.
Pour~$1 \le i, j \le n$, on d\'eclare que $i$ est \emph{avant}~$j$ dans l'\emph{ordre de Lyndon} associ\'e \`a~$w$ si $\dec^{i-1}(w)$ pr\'ec\`ede~$\dec^{j-1}(w)$ dans l'ordre lexicographique \'etendant~$\xx < \yy$. 
On note alors $i <_w j$.
\end{definition}

Par exemple, si $w$ est le mot $\xx\xx\yy\xx\yy$, on a $1 <_w 4$ car le mot~$\dec^3(\xx\xx\yy\xx\yy)$ est $\xx\yy\xx\xx\yy$, qui est plus grand lexicographement que~$\xx\xx\yy\xx\yy$. 
On trouverait de m\^eme $1 <_w 4 <_w 2 <_w 5 <_w 3$. 

\begin{lemma}
\label{L:ordreLyndon}
Soit $\gamma$ une orbite p\'eriodique de p\'eriode~$n$ du flot de Lorenz issue d'un point $P$ de l'axe~$A$.
Pour $1 \le i,j \le n$, le point $\pr^{i-1}(P)$ est \`a gauche du point $\pr^{j-1}(P)$ sur l'axe~$A$ si et seulement si on a $i <_{\cred(\gamma,P)} j$.
\end{lemma}

\begin{proof}
Notons $w$ le mot $\cred(\gamma,P)$ et $x$ l'abscisse du point~$P$. 
Par construction, $\pr$ est l'application $t \mapsto 2t \mod 1$\; comme l'orbite~$\gamma$ est p\'eriodique, elle ne passe jamais par le point $M$ d'abcisse~$1/2$, donc $x$ admet un d\'eveloppement en base~$2$ unique. 
De plus, la position de $\pr^{i-1}(P)$ dans la moiti\'e gauche ou droite de l'axe~$A$ d\'etermine le $i$-i\`eme chiffre du d\'eveloppement dyadique de~$x$, le cas ambig\"u \'etant justement exclu.
Par cons\'equent, si on remplace les $\xx$ par des 0 et les $\yy$ par des 1 dans~$w$,  l'abscisse de $\pr^{i-1}(P)$ est $\overline{0,w_iw_iw_i\ldots}$, o\`u $w_i$ d\'esigne le mot~$\dec^{i-1}(w)$. 
Or l'ordre canonique sur les r\'eels correspond \`a l'ordre lexicographique \'etendant $0<1$ sur le d\'eveloppement dyadique, donc $\pr^{i-1}(P)$ est \`a gauche de $\pr^{j-1}(P)$ si et seulement si $w_i$ pr\'ec\`ede $w_j$ dans l'ordre lexicographique.
\end{proof}

\begin{remark}
Soit $P$ un point d'intersection de $\gamma$ avec $A$ et $x$ son abscisse. 
Alors la condition $\pr^n(P)=\nobreak P$ se reformule en $2^nx=x\mod 1$, soit $x=m/(2^n-1)$ pour un certain~$m$ entier. 
Les abscisses des points d'intersection de~$\gamma$ avec $A$ sont alors les rationnels~$2^im/(2^n-1)\!\mod 1$. 
Inversement, on v\'erifie que l'orbite issue tout point d'abscisse de cette forme est p\'eriodique
\footnote{Si $k$ est impair, d'apr\`es le th\'eo\`eme d'Euler, $k$ divise $2^{\varphi(k)}-1$, donc tout nombre rationnel \`a d\'enominateur impair est de la forme $\frac{m}{2^n-1}$.}.
\end{remark}

\begin{definition}
\label{D:motLyndon}
Un mot~$w$ de longueur~$n$ sur l'alphabet~$\{\xx, \yy\}$ est dit {\it de Lyndon}
si $w$ est {\it primitif}, c'est-\`{a}-dire s'il n'est pas de la forme $u^k$ avec $k\ge 2$, et si $w$ pr\'ec\`ede lexicographiquement $\dec^i(w)$ pour tout $1 \le i < n$.
\end{definition}

Par exemple, les mots $\xx, \yy, \xx\yy, \xx\xx\yy, \xx\yy\yy, \xx\xx\xx\yy, \xx\xx\yy\yy, \xx\yy\yy\yy, \xx\xx\yy\xx\yy$ sont des mots de Lyndon. 
Par contre $\xx\yy\xx$ n'est pas de Lyndon puisque sont d\'ecal\'e $\xx\xx\yy$ le pr\'ec\`ede lexicographiquement.
Comme tout mot primitif admet un unique d\'ecal\'e qui soit de Lyndon, on d\'eduit que le nombre $\ell_n$ de mots de Lyndon de longeur~$n$ v\'erifie $\ell_n = \frac 1 n \left(2^n - \sum_{k \vert n} k\ell_k \right)$, d'o\`u $\ell_n \sim 2^n/n$.

\begin{proposition}[\cite{A:B-W}]
\label{T:Lyndon}
Pour chaque orbite p\'eriodique~$\gamma$ de p\'eriode~$n$ du flot de
Lorenz, il existe un et un seul mot de Lyndon de longueur~$n$
codant~$\gamma$, et, inversement, tout mot de Lyndon de longueur~$n$ est
code r\'eduit d'une orbite p\'eriodique de p\'eriode~$n$ du flot de
Lorenz.
\end{proposition}

\begin{proof}
Soit~$\gamma$ une orbite p\'eriodique de p\'eriode~$n$ et $P$ l'intersection de  $\gamma$ avec l'axe~$A$ la plus proche du point critique de gauche. 
Alors le code r\'eduit $\cred(\gamma,P)$ est un mot de Lyndon. 
En effet, si $\cred(\gamma,P)$ est de la forme $u^k$, alors le point $\pr^{\vert u \vert}(P)$ est confondu avec $P$, donc on a~$k=1$, et par cons\'equent le mot $\cred(\gamma,P)$ est primitif.
D'autre part, d'apr\`es le lemme~\ref{L:ordreLyndon}, les codes r\'eduits associ\'es \`a d'autres points d'intersection de~$\gamma$ avec $A$ sont apr\`es $\cred(\gamma,P)$ dans l'ordre lexicographique.

L'application ainsi construite est injective, car l'application~$\pr$ de premier retour sur~$A$ est dilatante, donc les orbites issues de deux points distincts de~$A$ ne peuvent co\"\i  ncider sur un nombre infini de tours cons\'ecutifs autour de l'attracteur, et donc leurs codes non plus.
Elle est surjective, puisque pour tout mot infini~$w$ sur l'alphabet~$\{\xx, \yy\}$ non ultimement constant, et en particulier pour tout mot p\'eriodique de p\'eriode au moins~2, $w$ est le code de l'orbite issue du point d'abscisse~$\overline{0,w}$, o\`u on a remplac\'e les lettres $\xx$ par des $0$ et les $\yy$ par des 1.
\end{proof}

D\'esormais, on appelle \emph{code de Lyndon} d'une orbite p\'eriodique l'unique mot de Lyndon qui la code. 
Par exemple, le code de Lyndon de l'orbite issue du point d'abscisse $1/7$ est le mot de Lyndon~$\xx\xx\yy$. 
De m\^eme, celui de l'orbite issue de~$5/31$ est le mot de Lyndon~$\xx\xx\yy\xx\yy$.

\begin{remark}
\label{R:PremierRetour}
L'hypoth\`ese que l'application~$\pr$ est identifi\'ee \`a $t\mapsto 2t \mod 1$ est cruciale. 
Si, par exemple, les deux bandes $B_1$ et $B_2$ du patron de Lorenz de la figure~\ref{F:patron} sont coll\'ees, non avec les parties $]0, 1/2[$ et $]1/2,1[$ de l'axe~$A$, mais avec les parties
$]0,k_0[$ et $]k_1,1[$, avec $k_0 < 1/2 < k_1$, alors tous les mots possibles sur l'alphabet~$\{\xx, \yy\}$ ne sont pas r\'ealisables comme suite des tours d'une orbite du syst\`eme de Lorenz modifi\'e. 
On a alors un syst\`eme $L_{k_0, k_1}$ dont l'ensemble des orbites est un sous-ensemble de celui qu'on vient de d\'ecrire. 
Les n\oe uds de Lorenz correspondent au cas le plus g\'en\'eral, o\`u tous les mots sont r\'ealis\'es par une orbite\footnote{Guckenheimer et Williams montrent dans~\cite{A:G-W} que, si le flot g\'eom\'etrique de Lorenz n'est pas {\it structurellement stable}, au sens o\`u une petite perturbation du champ de vecteur le d\'efinissant ne le transforme pas en un champ conjugu\'e, ce flot est n\'eanmoins {\it persistant}, au sens o\`u une petite perturbation l'envoie sur un champ d\'efinissant un flot de type $L_{k_0, k_1}$.}.
\end{remark}


\subsection{Tresses de Lorenz}
\label{S:TresseLorenz}

D'apr\`es un th\'{e}or\`eme c\'el\`ebre d'Alexander~\cite{A:Alexander}, tout n{\oe}ud est cl\^oture d'une tresse
\footnote{voir \cite{A:B-Z} pour une introduction \`a la th\'eorie des tresses}. 
Pour toute orbite p\'eriodique~$\gamma$ du flot de Lorenz, on d\'ecrit maintenant une tresse particuli\`ere dont la cl\^oture est~$\gamma$. 

\begin{figure}[htb]
	\begin{center}
	\begin{picture}(100,40)(0,0)
	\put(0,5){\includegraphics[scale=1.1]{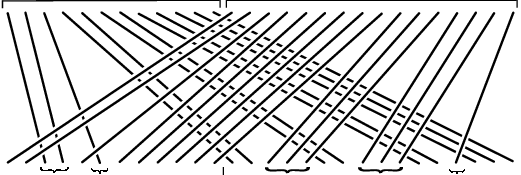}}
	\put(17,39){$p=12$} 
	\put(61,39){$n-p=16$} 
	\put(38, 2.5){$t=9$}
	\put(9,2.5){$n_1$} 
	\put(17,2.5){$n_2$} 
	\put(51.5,2.5){$m_6$}
	\put(69,2.5){$m_4$} 
	\put(83,2.5){$m_2$}
	\end{picture}
	\end{center}
	\caption{\small \sf Un exemple de tresse de Lorenz. 
	Les $p$ premiers brins vont vers la droite et les $n-p$ derniers vont vers la gauche en passant par-dessus les premiers.
	Le \trip~$t$ (d\'efinition~\ref{D:trip}) est le nombre de brins du premier groupe qui finissent dans le second groupe. 
	L'indice $n_i$ (d\'efinition~\ref{D:Young}) est le nombre de brins du premier groupe qui passent par-dessous $i+1$ brins du second groupe, l'indice $m_j$ est le nombre  de brins du second groupe qui passent par-dessus $j+1$ brins du premier groupe.}
	\label{F:permutation}
\end{figure}

On rappelle qu'\`a toute tresse $b$ \`a $n$ brins on associe une permutation $\pi$ de $\{1, \ldots, n\}$ en appelant $\pi(k)$ la position finale du brin initialement en position $k$. 
Chaque permutation~$\pi$ est associ\'ee \`a une infinit\'e de tresses.
Par contre, il en existe une et une seule dont tous les croisements sont orient\'es positivement\footnote{Il existe deux conventions d'orientation des croisements pour les n\oe uds et pour les tresses, l'une utilisée majoritairement par les théoriciens des n\oe uds (voir par exemple~\cite[p. 4]{A:G-P}), l'autre par les dynamiciens (voir~\cite[p. 6]{A:G-H-S}). Nous adoptons la convention des noueurs et décrétons que les croisements de la figure~\ref{F:permutation} où tous les brins sont orientés de haut en bas sont positifs. Cette convention implique par exemple que si un entrelacs se projette de sorte que les croisements entre ses composantes sont positifs, alors chaque composante coupe positivement une surface bordée par une autre composante. La convention des dynamiciens, qui est opposée, implique que si on part de courbes planaires de la forme~$(t,x_i(t))$ et qu'on les rel\`eve dans l'espace en~$(t,x_i(t), x'_i(t))$, alors les projections horizontales se coupent positivement.}, c'est-\`a-dire que le brin venant d'en haut \`a droite passe au-dessus de celui venant de gauche, et telle que deux brins quelconques se croisent au plus une fois. 
Cette tresse est appel\'ee {\it tresse de permutation associ\'ee \`a~$\pi$}.

\begin{theorem}[\cite{A:B-W}, voir figure~\ref{F:permutation}]
\label{T:tresseLorenz}
Soit $\gamma$ une orbite p\'eriodique du flot de Lorenz coupant $n$~fois l'axe~$A$ et coupant $p$~fois la partie $]0, 1/2[$ de cet axe. 
Alors $\gamma$ est la cl\^oture d'une (unique) tresse de permutation \`{a} $n$ brins. 
La permutation~$\pi$ associ\'{e}e est un cycle, et elle v\'erifie
\begin{equation}
\label{Eq:PermutationLorenz}
1<\pi(1)<\pi(2)<\ldots<\pi(p)=n>\pi(n)>\pi(n-1)>\ldots>\pi(p+1)=1.
\end{equation}

R\'eciproquement, toute tresse de permutation telle qu'il existe une permutation $\pi$ qui est un cycle et qui v\'erifie~\eqref{Eq:PermutationLorenz} se referme en un n\oe ud de Lorenz.
\end{theorem}

Dans le contexte ci-dessus, la tresse de permutation associ\'ee \`a~$\gamma$ est appel\'ee la \emph{tresse de Lorenz} de~$\gamma$, et la permutation~$\pi$ est appel\'ee la \emph{permutation de Lorenz} de~$\gamma$.

\begin{likeproof}[Esquisse de la d\'emonstration]
Il s'agit de d\'{e}crire en termes de tresses la partie centrale du patron de Lorenz tel qu'il est repr\'{e}sent\'{e} sur la figure~\ref{F:patron}. 
On obtient une tresse de la forme de celle de la figure~\ref{F:permutation}.
Comme l'application $\pr$ est croissante sur $]0, 1/2[$, deux brins quelconques parmi les $p$ premiers ne se coupent pas et on a $1<\pi(1)<\pi(2)<\ldots<\pi(p)$. 
Comme $\pr$ est aussi croissante sur $]1/2, 1[$, deux brins parmi les $n-p$ derniers ne se coupent pas non plus et on a $\pi(p+1)<\pi(p+2)<\ldots<\pi(n)<n$.
Qu'un brin du premier groupe passe sous un brin du deuxi\`eme découle du plongement du patron montr\'e sur la figure~\ref{F:patron}. On a alors bien une tresse de permutation v\'erifiant~\eqref{Eq:PermutationLorenz}. 
Comme on s'int\'eresse \`a un n\oe ud et non \`a un entrelacs, la permutation a un unique cycle.

R\'eciproquement, une tresse $b$ de permutation, associ\'ee \`a une permutation~$\pi_b$ v\'erifiant~\eqref{Eq:PermutationLorenz} et n'ayant qu'un cycle,  d\'etermine un mot de Lyndon $w_b$ \`a $n$ lettres comme suit: la $i$-\`eme lettre de $w_b$ est un $\xx$ si on a $\pi^{i-1}(1)\le p$ et un $\yy$ sinon. 
On v\'erifie alors que la cl\^oture de $b$ est le n\oe ud de Lorenz admettant $w_b$ pour code de Lyndon. 
\end{likeproof}

\begin{example}
\label{X:TresseLorenz}
La tresse de Lorenz associ\'ee \`a l'orbite issue du point d'abscisse~$1/7$ est la tresse $\sigma_2\sigma_1$ dont la cl\^oture est un n\oe ud trivial, et la permutation de Lorenz associ\'ee est le cycle $(123)$. 
En revanche, la tresse de Lorenz associ\'ee au point ~$5/31$ est $\sigma_3 \sigma_2 \sigma_1 \sigma_4 \sigma_3 \sigma_2$ dont la cl\^oture est un n\oe ud de tr\`efle, la permutation de Lorenz associ\'ee est le cycle $(13524)$.

\begin{figure}[h]
\begin{center}
\includegraphics{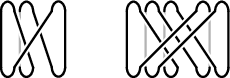}
\end{center}
\end{figure}
\end{example}


\subsection{Stabilisation}
\label{S:Stabilisation}

La proposition~\ref{T:Lyndon} fournit une \'enum\'eration de toutes les orbites p\'eriodiques de l'attracteur de Lorenz. 
Cependant, des orbites diff\'erentes peuvent \^etre isotopes, c'est-\`a-dire repr\'esenter le m\^eme n{\oe}ud, et donc plusieurs mots de Lyndon peuvent  coder le m\^eme n\oe ud de Lorenz. 
Pour r\'eduire la redondance du codage, on introduit maintenant la notion de stabilisation d'une orbite et du mot de Lyndon associ\'e.

\begin{figure}[htb]
	\begin{center}
	\includegraphics[scale=0.6]{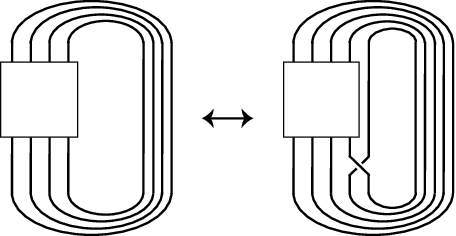}
	\end{center}
	\caption{\small \sf Une transformation de Markov\: la tresse $b$ de $\braid_n$  est remplac\'{e}e par la tresse $b\sigma_n^{\pm1}$ de $\braid_{n+1}$.}
	\label{F:Markov}
\end{figure}

Une {\it transformation de Markov} \`a droite sur une tresse $b$ \`{a} $n$ brins  consiste \`{a} ajouter un $n+1$-i\`eme brin et un croisement $\sigma_n^{\pm 1}$ (voir figure~\ref{F:Markov}). 
Elle ne change pas la classe d'isotopie de la cl\^oture de la tresse. 
On remarque que, si $b$ est une tresse de Lorenz \`{a} $n$ brins, alors $b\sigma_n$ est encore une tresse de Lorenz, cette fois-ci \`{a} $n+1$ brins. 
On peut de m\^{e}me introduire une transformation de Markov \`{a} gauche\: si on note $\tilde b$ la tresse $b$ o\`{u} tous les indices des g\'en\'erateurs~$\sigma_i$ ont \'{e}t\'{e} augment\'{e}s de 1, alors la tresse~$\tilde b\sigma_1$ est une tresse \`a $n+1$ brins dont la cl\^oture est la m\^{e}me que celle de $b$. 
Comme dans le cas de la transformation \`a droite, si $b$ est de Lorenz, alors $\tilde b\sigma_1$ est aussi de Lorenz. 
En termes d'orbites sur le patron de Lorenz, les transformations de Markov reviennent \`{a} ajouter une petite boucle faisant le tour d'un des deux points critiques de l'attracteur de Lorenz \`{a} la partie de l'orbite passant le plus pr\`es de celui-ci. 
La p\'{e}riode est alors augment\'ee de un
\footnote{En termes de permutations de Lorenz, une transformation de Markov \`a  droite correspond \`{a} modifier le motif $n \mapsto i$ en $n \mapsto n+1 \mapsto i$, et une transformation de Markov \`a gauche \`{a} augmenter d'une unit\'e tous les indices sup\'erieurs \`a~$1$ et \`a remplacer le motif $1\mapsto i+1$ en t\^{e}te par~$1\mapsto 2\mapsto i+1$.
}. 

Il est alors facile de traduire l'invariance des n\oe uds par transformations de Markov en termes de mots de Lyndon.

\begin{proposition}
\label{P:isotopiexy}
Soit $w$ un mot de Lyndon,
$w_{\xx}$ le mot obtenu en ajoutant une lettre $\xx$  en t\^ete de $w$ et
$w_{\yy}$ le mot obtenu en ins\'{e}rant une lettre $\yy$ apr\`es la
$i$-\`eme lettre de~$w$, o\`u $i$ est maximal pour l'ordre de Lyndon $<_w$. 
Alors les orbites p\'eriodiques du flot de Lorenz ayant $w$, $w_{\xx}$, et $w_{\yy}$ pour codes de Lyndon respectifs sont isotopes.
\end{proposition}

\begin{proof}
Soit $\gamma, \gamma_{\xx}$ et $\gamma_{\yy}$ les orbites ayant pour codes de Lyndon les mots $w, w_{\xx}$ et $w_{\yy}$. 
Alors, par d\'efinition du code, ajouter un $\xx$ en t\^ete de $w$ revient \`a ajouter une boucle faisant le tour de point critique gauche \`a l'orbite $\gamma$, soit, en termes de la tresse de Lorenz, \`a effectuer une transformation de Markov \`a gauche. 
Donc $\gamma_\xx$ est isotope \`a $\gamma$.

De m\^eme, si $i$ est maximal pour l'ordre de Lyndon associ\'e \`a~$w$, ajouter un $\yy$ apr\`es la $i$-\`eme lettre de~$w$ revient \`a effectuer une transformation de Markov, \`a droite cette fois, car, d'apr\`es le lemme~\ref{L:ordreLyndon}, la $i$-\`eme lettre de~$w$ correspond \`a l'intersection de~$\gamma$ avec l'axe~$A$ situ\'ee au plus pr\`es de l'attracteur droit. 
Donc $\gamma_\yy$ est isotope \`a~$\gamma$.
\end{proof}

Avec les notations de la proposition~\ref{P:isotopiexy}, on dira que $\gamma_\xx$ est obtenue par \emph{stabilisation} \`a gauche \`a partir de~$\gamma$, et que $\gamma_\yy$ est obtenue par \emph{stabilisation} \`a droite. 
De m\^eme, on dira que $w_\xx$ est obtenu par \emph{stabilisation} \`a gauche \`a partir de~$w$, et que $w_\yy$ est obtenu par \emph{stabilisation} \`a droite.

Par exemple, les mots $\xx\xx\xx\yy\xx\yy, \xx\xx\xx\xx\yy\xx\yy, \xx\xx\yy\yy\xx\yy$, et $\xx\xx\xx\xx\yy\yy\yy\xx\yy$ sont tous obtenus par stabilisation \`a partir du mot $\xx\xx\yy\xx\yy$. 
Les orbites associ\'ees sont donc toutes isotopes (en l'occurrence, elles repr\'esentent toutes des n\oe uds de tr\`efle).

Comme on peut stabiliser ind\'efiniment un mot de Lyndon  sans changer le n{\oe}ud correspondant, nous d\'eduisons:

\begin{corollary}
\label{C:Infinite}
Tout n{\oe}ud de Lorenz appara\^\i t comme classe d'isotopie d'une infinit\'e d'orbites p\'eriodiques du flot de Lorenz. 
\end{corollary}


\subsection{Diagrammes de Young}
\label{S:Young}

Lorsqu'on le repr\'esente comme sur la figure~\ref{F:permutation}, le diagramme d'une tresse de Lorenz \'evoque un diagramme de Young. 
On peut rendre cette correspondance formelle et, alors, les transformations de Markov correspondent \`a une op\'eration simple sur les diagrammes associ\'es.

\begin{figure}[htb]
	\begin{center}
	\includegraphics[width=\textwidth]{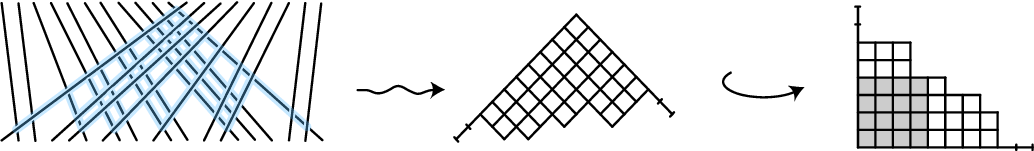}
	\end{center}
	\caption{\small \sf \`A une tresse de Lorenz on associe un diagramme de Young, augment\'e ici de deux segments \`a gauche et \`a droite, en ne gardant que la partie \'epaissie de la tresse de Lorenz. 
	Le \trip~$t$ (d\'efinition~\ref{D:trip}) est \'egal \`a la longueur du c\^ot\'e du carr\'e gris\'e, plus 1. 
	Pour tout $i>0$, il y a alors $n_i$ colonnes de hauteur~$i$ et $m_i$ lignes de largeur~$i$.} 
	\label{F:Young}
\end{figure}

On d\'efinit une correspondance entre tresses de Lorenz et diagrammes de Young comme suit. 
Soit $b$ une tresse de Lorenz. 
Dans le diagramme en segments de~$b$, on efface des fragments des brins de fa\c con \`a ne garder que les portions qui se trouvent entre deux croisements, voir figure~\ref{F:Young}. 
En tournant et d\'eformant la figure pour donner \`a chaque segment une orientation horizontale ou verticale, on obtient un diagramme de Young, \'eventuellement compl\'et\'e par deux segments, l'un vers la droite sur l'axe horizontal, l'autre vers le haut sur l'axe vertical.

\begin{definition}\label{D:Young}
Soit $\gamma$ une orbite p\'eriodique de Lorenz, $b$ la tresse de Lorenz \`a $n$ brins associ\'ee, $\pi$ la permutation de Lorenz et $p$ le nombre de brins de $b$ allant vers la droite, c'est-\`a-dire l'unique $p$ vérifiant~$\pi(p)=n$.
Pour $0\le i\le n$, on d\'efinit
	\begin{eqnarray}
	\label{E:ni}
	n_i &=& \mathrm{card}\{j\,\vert\,\pi(j)-j=i+1\,\,\mathrm{et}\,\,\pi(j)<\pi^2(j)\}\\
	\label{E:mi}
	m_i &=& \mathrm{card}\{j \,\vert\, j - \pi(j)=i+1\,\,\mathrm{et}\,\,\pi(j)>\pi^2(j)\}.
	\end{eqnarray}
Le \emph{diagramme de Young augment\'e associ\'e \`a~$b$} est le diagramme de Young ayant  un m\^at gauche de hauteur $n_0$, une base horizontale de longueur $m_0$, et $n_i$ colonnes de hauteur $i$.
\footnote{On dit que le diagramme est \emph{augment\'e} en raison du m\^at et de la base non standards qu'on lui a ajout\'es.}
\end{definition}

Les entiers $n_i$ ({\it resp.} $m_i$) correspondent au nombre de brins du premier ({\it resp.} second) groupe passant par-dessous  ({\it resp.} par-dessus) $i+1$ brins de l'autre groupe (voir figure~\ref{F:permutation}). 

Par exemple, le diagramme de Young compl\'et\'e associ\'e \`a l'orbite de code de Lyndon $\xx\xx\yy$ est un segment horizontal d'une unit\'e de long. Le diagramme de Young associ\'e au code de Lyndon $\xx\xx\xx\yy\xx\yy$ est un diagramme \`a deux cases avec un segment additionnel 
\begin{picture}(6,0)
	\put(0,-0.2){\includegraphics[scale=.7]{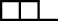}}
\end{picture}
.

\begin{proposition}
\label{P:Young}
La correspondance ci-dessus \'etablit une bijection entre les tresses de Lorenz et les diagrammes de Young augment\'es. 
\end{proposition}

\begin{proof}
Partons d'un diagramme de Young augment\'e de largeur~$p-1$ et de hauteur~$q-1$. 
On oriente le diagramme convenablement, on prolonge les segments, et on fait de chaque intersection un croisement de tresse positif. 
La tresse ainsi obtenue est une tresse de permutation dont ni les $p$ brins de gauche, ni les $q$ brins de droite ne se coupent entre eux. 
D'apr\`es la r\'eciproque dans le th\'eor\`eme~\ref{T:tresseLorenz}, c'est une tresse de Lorenz.
Il est clair que la correspondance ainsi d\'efinie est la r\'eciproque de celle d\'ecrite dans le chapeau de cette section.
\end{proof}

\begin{proposition}
\label{P:YoungStandard}
Tout n\oe ud de Lorenz peut \^etre obtenu comme cl\^oture d'une tresse de Lorenz associ\'ee \`a un diagramme de Young standard
(c'est-\`a-dire non augment\'e).
\end{proposition}

\begin{proof}
Soit $K$ un n\oe ud de Lorenz, et $b$ une tresse de Lorenz dont la cl\^oture est~$K$. 
D'apr\`es la proposition~\ref{P:Young}, la tresse~$b$ correspond \`a un diagramme~$T$, {\it a priori} augment\'e. 
Vue la description de la correspondance entre tresses de Lorenz et diagrammes de Young, supprimer un segment additionnel horizontal dans un diagramme compl\'et\'e revient \`a effectuer une d\'estabilisation gauche sur la tresse associ\'ee, tandis que supprimer un segment vertical revient \`a effectuer une d\'estabilisation droite. 
Partant du diagramme~$T$, on arrive en un nombre fini d'\'etapes du type pr\'ec\'edent \`a un diagramme standard. 
Par construction, celui-ci est associ\'e au n\oe ud initial~$K$.
\end{proof}

Par exemple, le n\oe ud de tr\`efle est associ\'e au tableau 
\begin{picture}(4,0)
	\put(0,-0.2){\includegraphics[scale=.7]{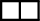}}
\end{picture}
. 
Z\'erologiquement, le n\oe ud trivial est associ\'e au tableau n'ayant aucune case. 

\begin{remark}
En termes de mot de Lyndon, un diagramme de Young non augment\'e correspond \`a un mot~$w$ qui est \emph{minimal} au sens o\`u, quand on \^ote \`a~$w$ la premi\`ere lettre~$\xx$, alors le mot~$w'$ obtenu n'est pas un mot de Lyndon et, quand on \^ote \`a~$w$ la lettre~$\yy$ dont la position~$i$ est maximale dans l'ordre de Lyndon $<_w$, alors $i$ n'est pas maximal dans $<_{w'}$. 
On d\'eduit alors directement de la proposition~\ref{P:YoungStandard} que tout n\oe ud de Lorenz peut \^etre repr\'esent\'e par une tresse de Lorenz cod\'ee par un mot de Lyndon minimal.

Par exemple, le mot de Lyndon $\xx\xx\xx\yy\xx\yy$ n'est pas minimal, puisque, si on lui \^ote le premier $\xx$, on obtient le mot $\xx\xx\yy\xx\yy$, qui est un mot de Lyndon. 
De m\^eme, le mot de Lyndon $\xx\xx\yy\yy\xx\yy$ n'est pas minimal. 
En effet, l'entier maximal dans l'ordre de Lyndon associ\'e est~$3$ et, quand on \^ote la lettre~$\yy$ en position~$3$, on obtient $\xx\xx\yy\xx\yy$, et $3$ reste maximal dans l'ordre de Lyndon associ\'e \`a ce dernier mot. 
Par contre, $\xx\xx\yy\xx\yy$ est minimal, puisque, d'une part,  $\xx\yy\xx\yy=(\xx\yy)^2$ n'est pas un
mot de Lyndon et, d'autre part, l'entier maximal dans l'ordre de Lyndon associ\'e \`a~$\xx\xx\yy\xx\yy$ est~$3$, alors que l'entier maximal dans l'ordre de Lyndon associ\'e \`a~$\xx\xx\xx\yy$ est~$4$.
\end{remark}

\`A ce point, on a donc une application qui associe \`a tout diagramme de Young un n\oe ud de Lorenz, et tout n\oe ud de Lorenz est repr\'esent\'e par au moins un diagramme de Young (standard, c'est-\`a-dire non augment\'e). 
La question de la redondance du codage ainsi obtenu se repose naturellement\: 

\begin{question}
Un n\oe ud de Lorenz donn\'e peut-il \^etre repr\'esent\'e par plusieurs diagrammes standards\? par une infinit\'e de diagrammes standards\? 
\end{question}

Il est facile de voir que la r\'eponse \`a la premi\`ere question est positive\: le n\oe ud de tr\`efle est repr\'esent\'e par chacun des deux diagrammes de Young	
\begin{picture}(2,0)
	\put(0,-1){\includegraphics[scale=.7]{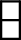}}
\end{picture}
et 
\begin{picture}(5,0)
	\put(0,-0.5){\includegraphics[scale=.7]{youngtrefle2.eps}}
\end{picture}, 
correspondant aux deux mots de Lyndon minimaux $\xx\xx\yy\xx\yy$ et $\xx\yy\xx\yy\yy$. 
Ce fait est un cas particulier d'un r\'esultat plus g\'en\'eral.
On appelle \emph{transpos\'e} d'un diagramme le diagramme obtenu par sym\'etrie autour de la premi\`ere diagonale, c'est-\`a-dire en
\'echangeant lignes horizontales et verticales.

\begin{proposition}
\label{T:transpose}
Les orbites de Lorenz associ\'ees \`a un diagramme de Young et \`a son transpos\'e sont isotopes.
\end{proposition}

\begin{proof}
Le patron de Lorenz est invariant par la rotation autour d'une droite passant par le milieu~$M$ de l'axe et \'echangeant les deux points critiques. Une orbite cod\'ee par un diagramme de Young est alors envoy\'ee sur l'orbite cod\'ee par le diagramme transpos\'e.
\end{proof}

Par contre, on \'etablira une r\'eponse n\'egative \`a la seconde question avec le corollaire~\ref{T:finitude}\: un n\oe ud de Lorenz ne peut \^etre repr\'esent\'e que par un nombre fini de diagrammes de Young non augment\'es. 


\section{Propri\'et\'es topologiques des n\oe uds de Lorenz}
\label{S:Topologie}

Dans cette partie, nous \'etudions l'\'etendue de la famille des n{\oe}uds de Lorenz, ainsi que leurs principales propri\'et\'es topologiques. 
Dans une br\`eve section~\ref{S:Primalite}, nous observons que tout n\oe ud de Lorenz est premier.
Dans la section~\ref{S:Torique}, nous montrons que tout n{\oe}ud torique est un n{\oe}ud de Lorenz, et que tout n{\oe}ud satellite d'un n{\oe}ud de Lorenz est de Lorenz, d'o\`u il r\'esulte que tout n{\oe}ud alg\'ebrique est de Lorenz. 
Puis, dans la section~\ref{S:Fibration}, nous montrons que tout n{\oe}ud de Lorenz est fibr\'e. 
Enfin, dans la section~\ref{S:Genre}, nous en d\'eduisons les genres des n\oe uds de Lorenz, avec une application au codage par les diagrammes de Young.

Pour la suite de l'article et pour plus de commodit\'es, nous identifions la sph\`ere~$\Sph^3$ avec $\reels^3 \cup \{\infty\}$, et plongeons tous les n\oe uds consid\'er\'es dans~$\Sph^3$.


\subsection{Primalit\'e}
\label{S:Primalite}

\begin{definition}
	\label{D:premier}
	Un n\oe ud $K$ est dit {\it premier} s'il n'est pas somme connexe de deux n\oe uds non~triviaux. 
\end{definition}

Autrement dit, un n\oe ud $K$ est premier si, pour toute sph\`ere~$S$ de dimension~2 coupant~$K$ en deux points, l'une des deux boules ainsi d\'elimit\'ees contient une corde non nou\'ee de~$K$.

\begin{theorem}[\cite{A:Williams}]
	\label{T:premier}
	Tout n\oe ud de Lorenz est premier.
\end{theorem}

\begin{likeproof}[Sch\'ema de la d\'emonstration]
	Soit~$K$ un n\oe ud de Lorenz et $S$ une sph\`ere coupant $K$ en deux points. 
	Le fait que $K$ soit une orbite p\'eriodique sur le patron de Lorenz donne des contraintes sur la position de $S$ dans~$\Sph^3$ et en particulier sur l'intersection de $S$ avec le patron de Lorenz, lesquelles se trouvent \^etre suffisamment fortes pour montrer que l'une des deux cordes d\'elimit\'ees par $S$ est non nou\'ee.
\end{likeproof}


\subsection{N\oe uds toriques, satellites et alg\'ebriques}
\label{S:Torique}

\begin{figure}[htb]
	\begin{center}
	\includegraphics[scale=0.6]{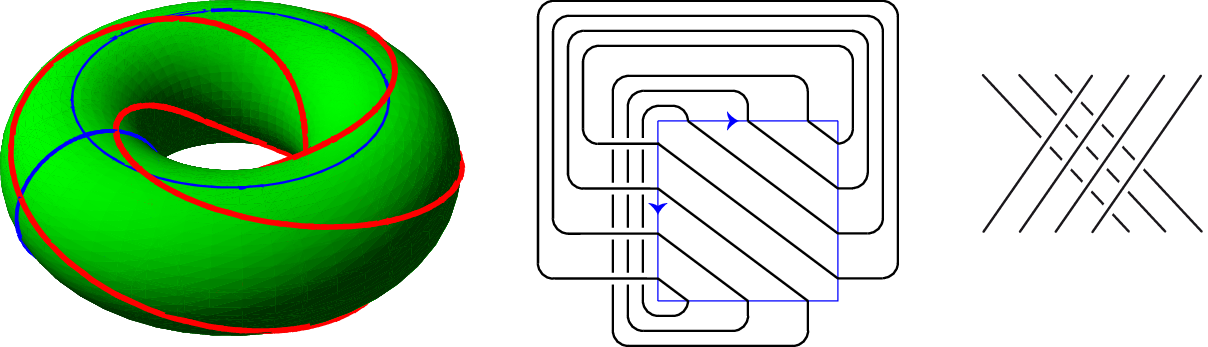}
	\end{center}
	\caption{\small \sf \`A gauche, le n\oe ud $T(4,3)$ dessin\'e \`a la surface d'un tore. 
	Au centre, la surface du tore est d\'ecoup\'ee le long d'un parall\`ele et d'un m\'eridien, puis d\'epli\'ee, faisant appara\^itre un diagramme avec $12$ croisements.
	\`A droite, une tresse positive dont la cl\^oture redonne le m\^eme diagramme.
	.} 
	\label{F:torique}
\end{figure}

\begin{definition}[voir figure~\ref{F:torique}]
	Soit $p$ et $q$ deux entiers relatifs premiers entre eux. 
	Un n\oe ud est dit {\it torique de type $(p,q)$} s'il peut \^etre trac\'e sur la surface d'un tore orient\'e plong\'e de mani\`ere standard dans $\Sph^3$, de sorte qu'il coupe $p$ fois chaque m\'eridien (orient\'e) du tore, et $q$ fois chaque parall\`ele (orient\'e). 
	On le note $T(p,q)$\footnote{\`A l'aide d'une isotopie de $\Sph^3$ \'echangeant les faces externes et internes du tore, on montre que les n\oe uds $T(p,q)$ et $T(-p,-q)$ sont isotopes. Par contre, les n\oe uds $T(p,q)$ et $T(p, -q)$ ne sont jamais isotopes.}.
\end{definition}

\vspace{-0.7cm}

\begin{proposition}[\cite{A:B-W}]
	Tout n\oe ud torique de type $(p,q)$ avec $p$ et $q$ positifs est un n\oe ud de Lorenz.
\end{proposition}

\begin{proof}
	Par construction, le n\oe ud~$T(p,q)$ est cl\^oture d'une tresse o\`u $p$ brins passent parall\`ellement par-dessus $q$~autres brins (voir figure~\ref{F:torique}). 
	D'apr\`es le th\'eor\`eme~\ref{T:tresseLorenz}, une telle tresse est une tresse de Lorenz.
\end{proof}

Le n\oe ud torique $T(p, -q)$ est l'image-miroir de~$T(p, q)$. 
Il est \'egalement obtenu comme cl\^oture d'une tresse o\`u $p$ brins venant de gauche passent parall\`ellement par-dessous $q$ autres venant de droite. 
Par contre, ce n'est pas un n\oe ud de Lorenz, l'orentation des croisements n'\'etant pas la bonne.

\begin{figure}[htb]
	\begin{center}
	\includegraphics[width=\textwidth]{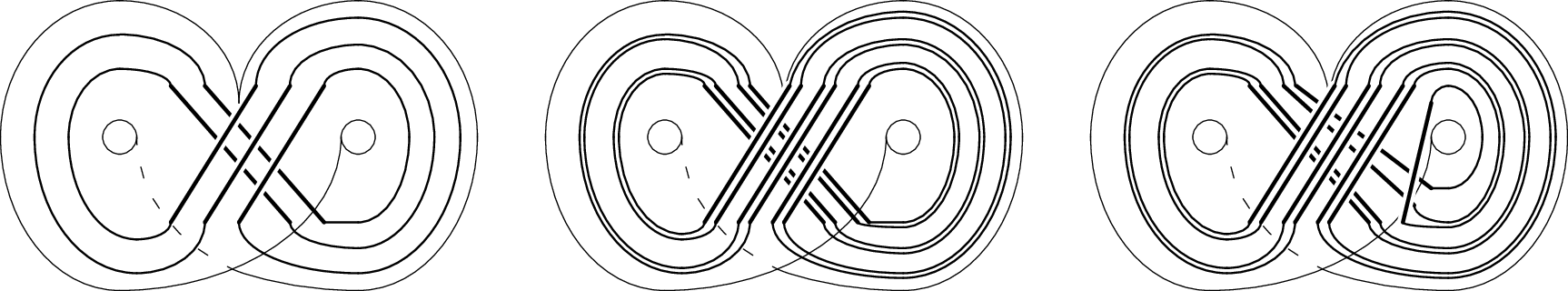}
	\end{center}
	\caption{\small \sf 
	Construction d'un n\oe ud satellite d'un n\oe ud de Lorenz. 
	On part d'un n\oe ud de Lorenz $K_c$  plong\'e dans le patron de Lorenz, ici \`a gauche. 
	On \'epaissit les brins de~$K_c$ en $p$ brins parall\`eles, ici au milieu avec $p=2$. 
	Au voisinage d'un quelconque des deux points critiques du patron, ici celui de droite, on fait passer le brin le plus proche du point critique par-dessous ses $p-1$ brins voisins. 
	Cette op\'eration peut \^etre r\'ep\'et\'ee un nombre $q$~premier avec $p$ quelconque de fois, ici $q=1$.
	D'apr\`es  la r\'eciproque dans le th\'eor\`eme~\ref{T:tresseLorenz}, la tresse obtenue est de Lorenz, et donc on obtient un n\oe ud de Lorenz, \`a droite.
	Ce n\oe ud est le satellite de $T(p,q)$ sur $K_c$ associ\'e \`a la trivialisation donn\'ee par le patron de Lorenz.
	} 
	\label{F:satellite}
\end{figure}

\begin{definition}[voir~\ref{F:satellite}]
	\label{D:satellite}
	Soit $K_c$ un n\oe ud dans~$\Sph^3$ et $K_i$ un n\oe ud dans~$\Sph^1\times \disque^2$ non isotope \`a~$\Sph^1\times\{0\}$.
	\'Etant donn\'es un voisinnage tubulaire $N_c$ de $K_c$ et une identification $f$ de $\Sph^1\times \disque^2$ avec $N_c$, le \emph{n\oe ud satellite} de $K_i$ sur $K_c$ associ\'e \`a $f$ est le n\oe ud $f(K_i)$. 
	Si $K_i$ est un n\oe ud torique, on dit aussi que le n\oe ud obtenu est un \emph{c\^ablage} de $K_c$.
\end{definition}

\begin{remark}
L'identification $f$ est importante, puisque le satellite est d\'etermin\'e \`a la fois par $K_c$, $K_i$ et par la classe d'isotopie de~$f$, c'est-\`a-dire par le choix d'un parall\`ele sur le tore~$f(\Sph^1\times\Sph^1)$, soit encore par le nombre d'enlacement entre $K_c=f(\Sph^1\times\{0\})$ et $f(\Sph^1\times\{1/2\})$.
\end{remark}

\vspace{-0.5cm}
\begin{figure}[htb]
	\begin{center}
	\begin{picture}(86,40)(0,0)
	\put(0,0){\includegraphics[width=0.6\textwidth]{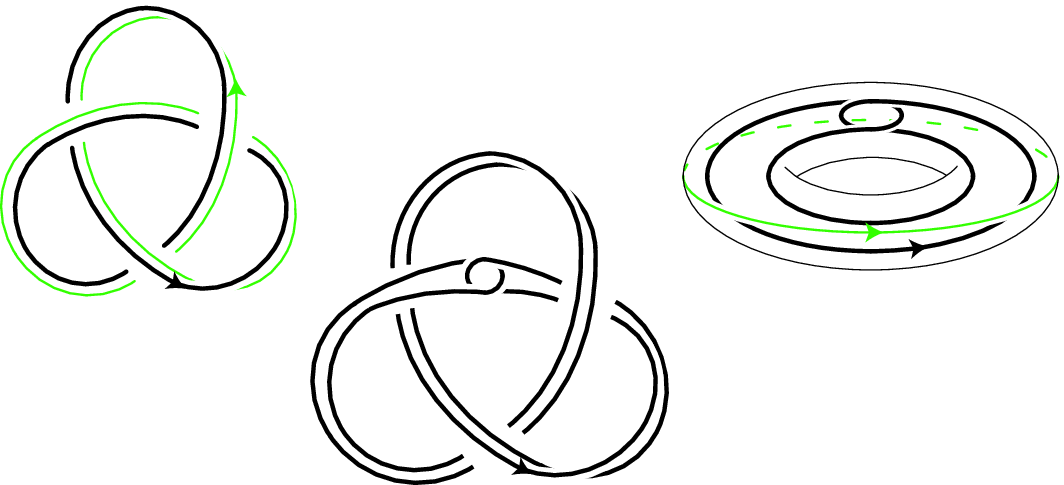}}
	\put(0,33){$K_c$}
	\put(86,33){$K_i$}
	\end{picture}
	\end{center}
	\caption{\small \sf Exemple de construction d'un n\oe ud satellite. 
	Le n\oe ud $K_c$ est un n\oe ud de tr\`efle (en haut \`a gauche), un choix de parall\`ele est \'egalement figur\'e. 
	Le n\oe ud $K_i$ est un n\oe ud trivial dans $\Sph^3$, mais pas dans le tore solide (en haut \`a droite) dont on a \'egalement marqu\'e un parall\`ele.
	Le satellite (au centre) est obtenu en envoyant le tore sur un voisinnage tubulaire de $K_c$ et en identifiant les deux parall\`eles marqu\'es.
	} 
	\label{F:ExempleSatellite}
\end{figure}

La notion de satellite g\'en\'eralise celle de somme connexe\: s'il existe $x$ dans~$\Sph^1$ et un disque dans~$\disque^2\times \Sph^1$ isotope \`a $\disque^2\times \{x\}$ coup\'e une fois exactement par $K_i$, alors tout satellite de $K_i$ sur $K_c$ est la somme connexe de ces deux n\oe uds. 
 
\begin{proposition}[\cite{A:B-W}]
	\label{T:cablage}
	Si $K_c$ est un n\oe ud de Lorenz et $T(p,q)$ un n\oe ud torique, avec $p, q > 0$, plong\'e naturellement dans un tore plein, alors le satellite de $T(p,q)$ sur $K_c$ associ\'e \`a la trivialisation donn\'ee par le patron de Lorenz est encore un n\oe ud de Lorenz.
\end{proposition}

La d\'emonstration est esquiss\'ee sur la figure~\ref{F:satellite}. 

\'Etant donn\'ee une courbe alg\'ebrique complexe~$H$ \`a singularit\'es isol\'ees dans~$\Compl^2$, identifi\'e avec~$\reels^4$, et un point $P$ sur $H$, l'intersection~$K_{H, \varepsilon}$ entre $H$ et la sph\`ere r\'eelle de rayon~$\varepsilon$ centr\'ee en $P$ est un entrelacs, avec \'eventuellement des points singuliers.
Dans~\cite{A:Milnor}, J.\,Milnor montre que, pour $\varepsilon$ assez petit, l'entrelacs~$K_{H, \varepsilon}$ n'a pas de points singuliers, et par cons\'equent sa classe d'isotopie est ind\'ependante du rayon $\varepsilon$. 
Si $P$ est un point r\'egulier de la courbe, alors l'entrelacs en question est le n\oe ud trivial.
Par contre, si $P$ est un point singulier irr\'eductible, c'est un n\oe ud, et il n'est jamais trivial.

\begin{definition}[\cite{A:Milnor}]
	Un n\oe ud $K$ est dit {\it alg\'ebrique} s'il est la classe d'isotopie de l'intersection d'une courbe alg\'ebrique de~$\Compl^2$, identifi\'e avec~$\reels^4$, avec une sph\`ere centr\'ee sur un point de la surface et de rayon suffisamment petit.
\end{definition}

\begin{example}\label{Ex:TrefleAlgebrique}
	Consid\'erons le polyn\^ome $P(x,y)=y^3-x^2$. Pour tout $\varepsilon$, l'intersection $K_\varepsilon$ du lieu d'annulation de $P$ et de la sph\`ere de rayon $\varepsilon$ est d\'efinie par $x^2=y^3$ et $\vert x \vert^2 + \vert y \vert^2 = \varepsilon$. 
	\'Etant donn\'e un point $(\lambda, \mu)$ du n\oe ud $K_\varepsilon$, ce dernier est alors param\'etr\'e par $t\mapsto (\lambda e^{2i\pi t/2}, \mu e^{2i\pi t/3})$.
	Il s'agit d'un n\oe ud de tr\`efle.
\end{example}

\begin{proposition}[\cite{A:B-W}]
	\label{T:noeudsalgebriques}
	Tout n\oe ud alg\'ebrique est un n\oe ud de Lorenz.
\end{proposition}

\begin{likeproof}[Sch\'ema de la d\'emonstration]
	Il est d\'emontr\'e dans~\cite{A:E-N} et~\cite{A:Milnor} que tout n\oe ud alg\'ebrique peut \^etre obtenu par c\^ablages successifs \`a partir d'un n\oe ud torique. 
	La proposition~\ref{T:cablage} assure alors qu'\`a chaque \'etape de c\^ablage, le n\oe ud obtenu est de Lorenz.
\end{likeproof}

La r\'eciproque de la proposition~\ref{T:noeudsalgebriques} est fausse.
Par exemple, le n\oe ud associ\'e au mot de Lyndon $\xx\yy\xx\yy^3\xx\yy^3$ n'est pas alg\'ebrique, puisque par exemple son polyn\^ome d'Alexander a des racines de module diff\'erent de~$1$, au contraire de tous les n\oe uds alg\'ebriques.
Le r\'esultat suivant, d\^u \`a M.\,El\,Rifai, est une r\'eciproque partielle de la proposition~\ref{T:cablage}.

\begin{theorem}[\cite{A:ElRifai2}]
	\label{T:satellite}
	Tout n\oe ud de Lorenz qui est un satellite d'un n\oe ud de Lorenz est satellite selon le sch\'ema de la proposition~\ref{T:cablage}, c'est-\`a-dire un c\^ablage sur un n\oe ud de Lorenz.
\end{theorem}

\begin{likeproof}[Sch\'ema de la d\'emonstration]
	L'id\'ee de base est la m\^eme que pour d\'emontrer la primalit\'e\: si~$K$ est un n\oe ud de Lorenz qui est un satellite, alors il existe un tore $\tore$ plong\'e de mani\`ere nou\'ee dans~$\Sph^3$, contenant $K$ en son int\'erieur mais ne contenant pas de boule contenant~$K$. 
	Encore une fois, le fait que $K$ soit une orbite p\'eriodique du flot de Lorenz donne suffisamment de contraintes sur $\tore$ et sur son intersection avec le patron de Lorenz pour montrer qu'il est un voisinage tubulaire d'un n\oe ud de Lorenz, et qu'en son int\'erieur le n\oe ud $K$ ne fait que tourner comme un n\oe ud torique. 
\end{likeproof}

Ce dernier r\'esultat est sp\'ecialement int\'eressant, puisqu'il d\'etermine quels n\oe uds de Lorenz sont alg\'ebriques, et, lorsque c'est le cas, comment ils apparaissent comme orbites p\'eriodiques du flot de Lorenz. 


\subsection{Caract\`ere fibr\'e}
\label{S:Fibration}

Le th\'eor\`eme~\ref{T:tresseLorenz} \'enonce une propri\'et\'e capitale des n\oe uds de Lorenz, \`a savoir qu'ils peuvent \^etre r\'ealis\'es comme cl\^otures de tresses positives. 
Une cons\'equence notable est le fait qu'ils sont fibr\'es.
D'une part, ceci permet calculer le genre des n\oe uds de Lorenz (proposition~\ref{T:genre}), et,  d'autre part, cela corrobore l'intuition que les n\oe uds de Lorenz sont des n\oe uds simples\: la structure du compl\'ementaire d'un n\oe ud fibr\'e est bien mieux comprise que celle d'un n\oe ud quelconque.

\begin{definition}
	\label{D:Seifert}
	Soit $K$ un entrelacs orient\'e dans $\Sph^3$. 
	On appelle {\it surface de Seifert}  pour $K$ une surface orient\'ee de $\Sph^3$ dont le bord orient\'e est~$K$. 
	On appelle {\it genre} de~$K$ le genre minimal d'une surface de Seifert pour~$K$. 
\end{definition}

Notons qu'un n\oe ud est trivial si et seulement si il borde un disque dans~$\Sph^3$, donc si et et seulement si son genre est nul. 

\begin{figure}[htbp]
	\begin{center}
	\includegraphics[scale=0.4]{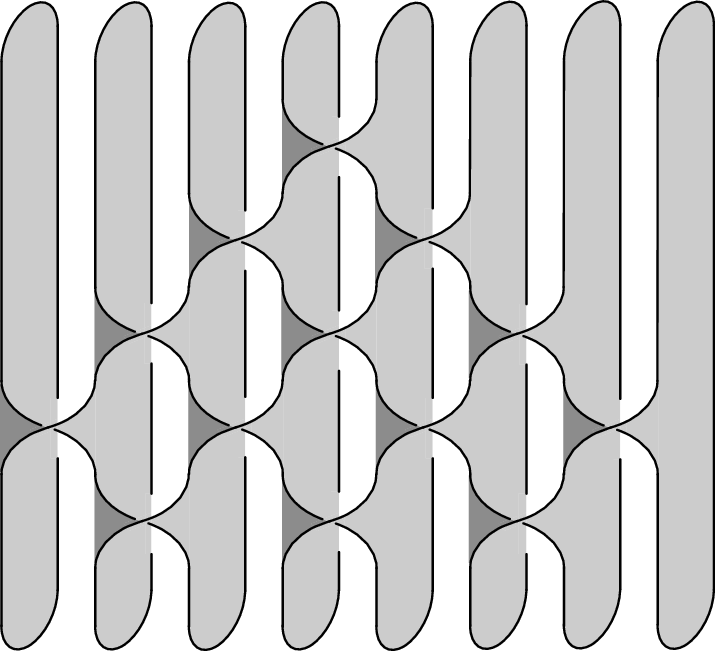}
	\end{center}
	\caption{\small \sf Une surface de Seifert pour le n\oe ud de Lorenz de code $\xx\xx\yy\xx\yy\yy\xx\yy$. 
	La m\^eme construction fournit une surface de Seifert, dite \emph{standard}, pour toute cl\^oture de tresse.
	D'apr\`es la proposition~\ref{T:genreminimal}, cette surface est de genre minimal dans le cas des tresses positives.}
	\label{F:seifertpositive}
\end{figure}

\begin{definition}
	\label{D:noeudfibre}
	Un n\oe ud~$K$ dans~$\Sph^3$ est dit \emph{fibr\'e} si son compl\'ementaire \emph{fibre} sur le cercle~$\Sph^1$, c'est-\`a-dire s'il existe une surface de Seifert~$\surf_K$ pour~$K$ et un diff\'eomorphisme~$h$  de~$\surf_K$ pr\'eservant le bord tels que le compl\'ementaire de $K$ dans $\Sph^3$ est isomorphe \`a $\surf_K \times [-\pi, \pi] /_{(h(x),-\pi)\sim (x,\pi)}$.
\end{definition}

On remarque que, via l'isomorphisme de la d\'efinition~\ref{D:noeudfibre}, la projection~$\theta$ sur le deuxi\`eme facteur d\'efinit une fonction de~$\Sph^3 \smallsetminus K$ sur $\Sph^1=[-\pi,\pi]/_{-\pi\sim\pi}$ qui poss\`ede les propri\'et\'es suivantes\footnote{On dit alors que couple $(K, \theta)$ est un \emph{livre ouvert} dont $K$ est la \emph{reliure} et les fibres $\theta^{-1}(t)$ sont les \emph{pages}.}\:

	$(i)$ $\theta$ est lisse et sans point critique\;

	$(ii)$ sur un voisinage tubulaire~$N(K)$ du n\oe ud hom\'eomorphe \`a un tore plein $\Sph^1\times\disque^2$, la fonction~$\theta$ s'identifie \`a la fonction argument sur le facteur~$\disque^2$. 
	
Le n\oe ud trivial est un exemple de n\oe ud fibr\'e. On part d'un cercle en m\'etal que l'on trempe dans du savon liquide\; quand on ressort le cercle, le savon forme un disque bord\'e par le n\oe ud en m\'etal, lequel est topologiquement trivial. 
En soufflant sur la pellicule de savon pour la gonfler jusqu'\`a son \'eclatement, on d\'efinit une fibration sur l'intervalle~$[0, \pi[$. 
En refaisant la m\^eme op\'eration en soufflant de l'autre c\^ot\'e, on d\'efinit une fibration sur~$]-\pi,0]$. 

Regardons $\Sph^3$ comme la sph\`ere unit\'e de~$\Compl^2$, c'est-\`a-dire comme l'ensemble des couples~$(z_1,z_2)$ de nombres complexes v\'erifiant $\module{z_1}^2 + \module{z_2}^2=1$. 
Pour tout $\lambda$ dans $\Compl\cup\{\infty\}$, l'ensemble $\{(z_1, z_2) \mid z_1=\lambda z_2\}$ est un n\oe ud trivial dans~$\Sph^3$. 
Un entrelacs constitu\'e par une union quelconque de tels n\oe uds triviaux est un entrelacs dit \emph{de Hopf}\: deux composantes distinctes ont un nombre d'enlacement \'egal \`a~$+1$. 
Par exemple, l'ensemble $\{(z_1, z_2) \in \Sph^3 \mid z_1=\pm z_2\}$ est un entrelacs de Hopf \`a deux composantes. 

\begin{proposition}[\cite{A:Milnor}]
	\label{T:Hopffibre}
	Tout entrelacs de Hopf est fibr\'e.
\end{proposition}

\begin{proof}
	Soit $H_n$ un entrelacs de Hopf \`a $n$ composantes. 
	Par définition, il existe des nombres complexes $k_1, \dots, k_n$ différents de~$1$ tels que $H_n$ est l'intersection dans~$\reels^4$ de~$\Sph^3$ et du lieu d'annulation du polyn\^ome 
	\[P_n(z_1,z_2)=\prod_{k=1}^{n}(z_1 - k_i z_2).\]
	Pour tout~$\theta$ dans $[-\pi,\pi]$, l'ensemble $\{(z_1, z_2) \in\Sph^3 \mid \arg(P_n(z_1,z_2))=\theta\}$ est une surface dont le bord est~$H_n$. 
	L'argument du polyn\^ome~$P$ est alors une fibration de~$\Sph^3\smallsetminus H_n$ sur~$\Sph^1$.
\end{proof}

La d\'emonstration pr\'ec\'edente peut \^etre adapt\'ee pour montrer que tout entrelacs alg\'ebrique est fibr\'e~\cite{A:Milnor}. 
Le r\'esultat principal de cette section est une g\'en\'eralisation qui implique en particulier que les n\oe uds de Lorenz sont fibr\'es.

\begin{theorem}[\cite{A:B-W, A:Stallings}]
	\label{T:fibrationpositive}
	Tout entrelacs qui est cl\^oture d'une tresse positive est fibr\'e.
\end{theorem}

L'id\'ee de la d\'emonstration que nous pr\'esentons est de construire de proche en proche une fibration pour chaque tresse positive en ajoutant un \`a  un les croisements. 
Comme l'entrelacs de Hopf \`a deux composantes est fibr\'e, on l'utilise comme brique \'el\'ementaire. 
D'autre part, on peut recoller certaines fibrations dans~$\Sph^3$ si elles sont non triviales sur des ensembles suffisamment disjoints. 
La notion de somme de Murasugi formalise cette id\'ee, et on va l'utiliser comme ciment de la construction.

\begin{figure}[htbp]
	\begin{center}
	\includegraphics[scale=0.3]{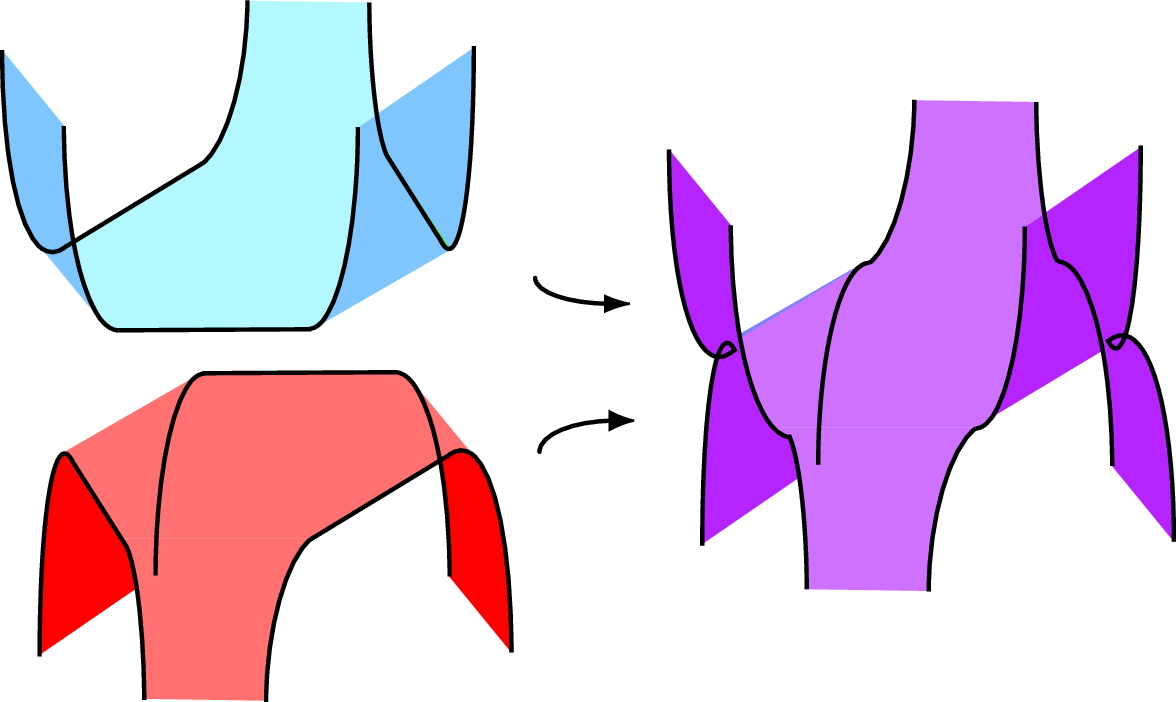}
	\end{center}
	\caption{\small \sf La somme de Murasugi de deux surfaces \`a bord.}
	\label{F:sommeMurasugi}
\end{figure}

\begin{definition}(voir figure~\ref{F:sommeMurasugi})
	\label{D:sommemurasugi}
	Soit $\surf_1$ et $\surf_2$ deux surfaces orient\'ees plong\'ees dans $\Sph^3$ de  bords respectifs~$K_1$ et $K_2$. 
	Soit~$\Pi$ une sph\`ere (que l'on voit comme un plan horizontal dans $\reels^3\cup\{\infty\}$) s\'eparant~$\Sph^3$ en deux boules ouvertes~$B_1$ et $B_2$. 
	On suppose que 

$(i)$ $\surf_1$ est incluse dans l'adh\'erence de~$B_1$ et $\surf_2$ dans l'adh\'erence de~$B_2$\;

$(ii)$ $\surf_1\cap\surf_2$ est un $2n$-gone $P$ contenu dans $\Pi$\;

$(iii)$ $K_1$ et $K_2$ sont deux entrelacs se coupant uniquement aux sommets~$x_1, \dots, x_{2n}$ de $P$.

\noindent On d\'efinit alors la \emph{somme de Murasugi $\surf_1\Murasomme_P\surf_2$ des surfaces $\surf_1$ et $\surf_2$ le long de $P$} comme leur r\'eunion $\surf_1\cup\surf_2$. 
	Par extension, on d\'efinit la \emph{somme de Murasugi $K_1\Murasomme_P K_2$ des n\oe uds $K_1$ et $K_2$ le long de $P$} comme le n\oe ud $K_1\cup K_2 \smallsetminus \bigcup ] x_i , x_{i+1} [$
\footnote{Plus g\'en\'eralement, on d\'efinit la somme de Murasugi de surfaces $\surf_i, i=1,2$ le long de deux polygones $P_i$ dont un c\^ot\'e sur deux est dans le bord $\partial \surf_i$ comme la classe d'isotopie dans~$\Sph^3$ de la somme de deux surfaces $\widehat{P_i}$ isotopes \`a $P_i$ et v\'erifiant les crit\`eres $(i), (ii)$ et $(iii)$, le r\'esultat ne d\'ependant alors pas du choix des $\widehat{P_i}$ dans leur classe d'isotopie. On la note~$\surf_1\Murasomme_{P_1\sim P_2}\surf_2$}.
\end{definition}

La somme de Murasugi g\'en\'eralise la somme connexe, qui correspond au cas $n=1$ dans la d\'efinition. 
Elle g\'en\'eralise \'egalement le plombage, qui correspond au cas $n=2$. 
La somme de Murasugi est une op\'eration g\'eom\'etrique naturelle pour les surfaces (et pour les entrelacs qu'elles bordent), \`a savoir qu'elle conserve plusieurs propri\'et\'es importantes~\cite{A:Gabai1, A:Gabai2}, comme par exemple le fait d'\^etre une surface incompressible, le fait d'\^etre une surface de Seifert de genre minimal, ou encore le fait d'\^etre un entrelacs fibr\'e, comme nous allons le voir. 

L'id\'ee fondamentale pour montrer que la somme de Murasugi de deux entrelacs fibr\'es est fibr\'ee consiste \`a d\'eformer pr\'ealablement les fibrations associ\'ees pour qu'elle soient non triviales sur des ensembles disjoints, d'o\`u la d\'efinition suivante.

\begin{definition}\label{D:fibrationsuperieure}
	Soit $K$ un entrelacs fibr\'e, $\theta:\Sph^3\smallsetminus K \mapsto \Sph^1$ d\'efinissant la fibration associ\'ee et $\Pi$ une sph\`ere divisant $\Sph^3$ en deux boules $B_1$ et $B_2$. 
	Soit $N(K)$ un voisinnage tubulaire de $K$. 
	On dit que la fibration $\theta$ est \emph{\tangente~inf\'erieurement} (resp. \emph{sup\'erieurement}) \`a~$\Pi$ si 

$(i)$ il existe un polygone $P$ de sommets $x_1,x_2,\dots,x_{2n}$ sur~$\Pi$ tel que $K$ rencontre $\Pi$ le long des c\^ot\'es $[x_1 x_2], [x_3x_4], \dots, [x_{2n-1}x_{2n}]$ (\emph{resp.} $[x_2 x_3], \dots, [x_{2n}x_{1}]$) de $P$ et $K$ ne rencontre pas $B_2$ (\emph{resp.}~$B_1$)\;

$(ii)$ $N(K)\cap \Pi$ consiste en un $\varepsilon$-voisinage des c\^ot\'es $[x_1 x_2], \dots, [x_{2n-1}x_{2n}]$ de $P$ (\emph{resp.} $[x_2 x_3], \dots$)\;

$(iii)$ \`a l'ext\'erieur du $\varepsilon$-voisinage $N(\partial P)$ du bord du polygone $P$ dans~$\Pi$, la fonction $\theta$ vaut 0 \`a l'int\'erieur de $P$ et $\pi$ \`a l'ext\'erieur\;

$(iv)$ \`a l'int\'erieur de $N(\partial P)\smallsetminus N(K)$, la fonction~$\theta$ prend ses valeurs dans l'intervalle~$]0, \pi[$ (\emph{resp.}~$]-\pi, 0[$).
\end{definition}  

Si on voit une fibration comme le d\'eplacement continu d'une fibre qui remplit tout le compl\'ementaire d'un n\oe ud, alors une fibration est tangente sup\'erieurement \`a un plan si, sous ce plan, la fibre se d\'eplace comme la bulle de savon d\'ecrite pour la fibration du n\oe ud trivial.
Le r\'esultat suivant assure qu'\'etant donn\'e un entrelacs fibr\'e et un polyg\^one $P$ inclus dans une fibre, on peut supposer la fibration tangente \`a une petite sph\`ere \'epaississant~$P$.

\begin{lemma}
	\label{L:positionfibre}
	Soit $K$ un entrelacs fibr\'e et $\theta:\Sph^3\smallsetminus K \to \Sph^1$ la fibration associ\'ee. 
	Soit $P$ un polygone \`a $2n$ c\^ot\'es inclus dans une fibre $\theta^{-1}(t)$, et rencontrant $K$ le long ses c\^ot\'es $[x_1 x_2], \dots, [x_{2n-1}x_{2n}]$ (\emph{resp.}$[x_2 x_3], \dots, [x_{2n}x_{1}]$). 
	Alors il existe une sph\`ere $\Pi$ dans $\Sph^3$ contenant $P$, rencontrant $K$ uniquement le long des arcs $[x_1 x_2], \dots, [x_{2n-1}x_{2n}]$ (\emph{resp.}$[x_2 x_3], \dots, [x_{2n}x_{1}]$), et telle que $K$ soit contenu dans une seule des deux boules ferm\'ees d\'elimit\'ees par~$\Pi$.
	De plus, on peut reparam\'etrer $\theta$ de sorte que $\theta$ soit tangente inf\'erieurement (\emph{resp.} sup\'erieurement) \`a~$\Pi$.
\end{lemma}

\begin{figure}[htbp]
	\begin{center}
	\begin{picture}(100,40)(0,0)
	\put(-4,-8){\includegraphics[width=0.8\textwidth]{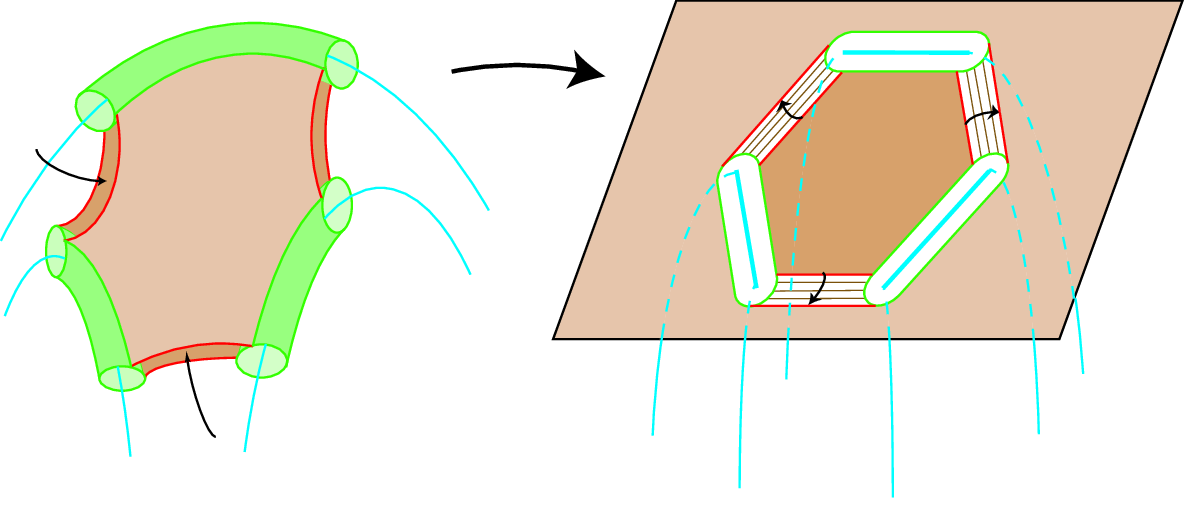}}
	\put(-2,29){$t$} 
	\put(18,-2){$t$}
	\put(13,23){$t+\varepsilon$}
	\put(82,24){$0$} 
	\put(93,15){$\pi$}
	\end{picture}
	\end{center}
	\caption{\small \sf D\'eformation de la fibration pour obtenir une fibration triviale dans le demi-plan sup\'erieur.}
	\label{F:fibrationreparametree}
\end{figure}

\begin{proof}(voir figure~\ref{F:fibrationreparametree})
	Comme $P$ est contractile, au voisinnage de celui-ci, la fibration $\theta$ est isomorphe \`a la projection sur le second facteur $P\times[t-\varepsilon,t+\varepsilon] \to [t-\varepsilon,t+\varepsilon]$. 
	Soit $B_2$ l'image r\'eciproque de  $P\times[0,t+\varepsilon]$ et $\Pi$ son bord. 
	La sph\`ere $\Pi$ consiste alors en deux copies de $P$, l'une dans la fibre~$\theta^{-1}(t)$ et l'autre dans une fibre $\theta^{-1}(t+\varepsilon)$, et un anneau $\partial P\times ]0,\varepsilon[$, comme sur la partie gauche de la figure~\ref{F:fibrationreparametree}.
	
	Composons $\theta$ \`a gauche par un diff\'eomorphisme du cercle envoyant l'intervalle $[t, t+\varepsilon]$ sur $[0, \pi]$. 
	Si on pense au param\`etre de la fibration comme \`a un temps, cela revient \`a ralentir le temps pour passer un temps $\pi$ dans la (petite) boule $B_2$. 
	La fibration obtenue est tangente inf\'erieurement \`a~$\Pi$. 
	En effet, si on zoome sur $B_2$ pour lui faire remplir tout le demi-espace sup\'erieur, la fibration est alors semblable \`a celle qui est repr\'esent\'ee sur la partie droite de la figure~\ref{F:fibrationreparametree}, laquelle est tangente inf\'erieurement au plan horizontal.
\end{proof}

Ces pr\'eliminaires nous permettent de montrer le r\'esultat-cl\'e suivant.

\begin{theorem}[\cite{A:Gabai1}]
	\label{T:sommefibre}
	Soit $K_1$ et $K_2$ deux entrelacs fibr\'es dans~$\Sph^3$ admettant pour fibres respectives deux surfaces $\surf_1$ et $\surf_2$.
	Soit $P_1$ ({\it resp.} $P_2$) un polygone \`a $2n$~c\^ot\'es sur $\surf_1$ ({\it resp.} $\surf_2$) dont les c\^ot\'es pairs ({\it resp.} impairs) sont inclus dans le bord $K_1$ ({\it resp.} $K_2$) de $\surf_1$ ({\it resp.} $\surf_2$).
	Alors la somme de Murasugi $K_1\Murasomme_{P_1\sim P_2} K_2$ est fibr\'ee de fibre $\surf_1\Murasomme_{P_1\sim P_2}\surf_2$.
\end{theorem}

\begin{figure}[htbp]
	\begin{center}
	\includegraphics[width=0.7\textwidth]{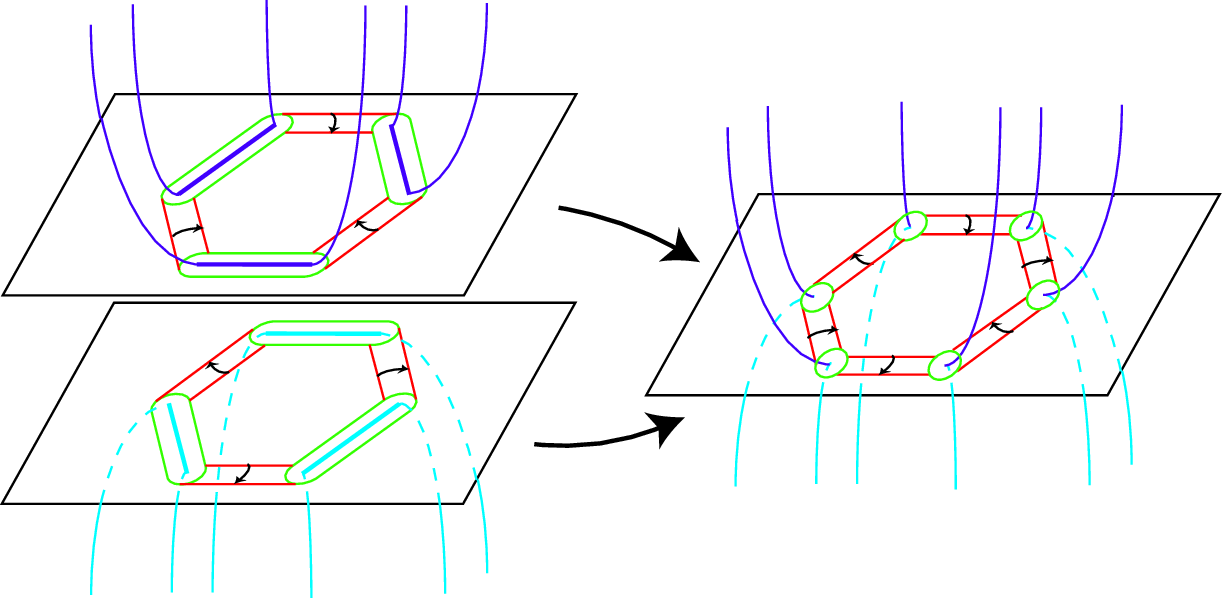}
	\end{center}
	\caption{\small \sf Le collage de deux fibrations pour obtenir une fibration de la somme de Murasugi de deux n\oe uds. 
	}
	\label{F:collagefibration}
\end{figure}

\begin{proof}(voir figure~\ref{F:collagefibration})
Soit $\Pi$ le plan horizontal dans $\reels^3$. 
Pour $i=1,2$, soit $\theta_i$ une fibration associ\'ee \`a $K_i$ avec $\surf_i$ pour fibre.
D'apr\`es le lemme~\ref{L:positionfibre}, on peut d\'eformer $K_1$ et $\theta_1$ de sorte que $\theta_1$ soit tangente inf\'erieurement \`a~$\Pi$, et de m\^eme on peut d\'eformer $K_2$ et $\theta_2$ de sorte que $\theta_2$ soit tangente sup\'erieurement \`a~$\Pi$. 
La somme de Murasugi $K_1\Murasomme_{P_1\sim P_2} K_2$ est alors obtenue en recollant  le demi-espace inf\'erieur associ\'e \`a $K_1$ et le demi-espace sup\'erieur associ\'e \`a $K_2$, et en supprimant les arcs formant le bord du polyg\^one~$P=x_1\ldots x_{2n}$ dans~$\Pi$. 

On d\'efinit une nouvelle fonction $\theta$, dont nous montrerons par la suite qu'elle induit une fibration. 
Sur le demi-espace inf\'erieur, $\theta$ co\"\i  ncide avec~$\theta_1$, sauf sur un voisinage tubulaire des arcs $[x_1, x_2], \ldots, [x_{2n-1}x_{2n}]$ (o\`u $\theta_1$ est suppos\'ee du type fonction-argument). 
De m\^eme, sur le demi-espace sup\'erieur, $\theta$ co\"\i  ncide avec $\theta_2$, sauf sur un voisinage tubulaire des arcs~$[x_2, x_3],\ldots, [x_{2n}x_1]$. 

Au voisinage des arcs~$[x_1, x_2], [x_3, x_4], \ldots, [x_{2n-1}x_{2n}]$, la fonction $\theta$ est d\'efinie sur le demi-espace sup\'erieur, et varie de~$\pi$ \`a $2\pi$ quand on se d\'eplace de l'ext\'erieur vers l'int\'erieur du polyg\^one. 
Elle est \'egalement d\'efinie sur le demi-espace inf\'erieur jusqu'\`a la fronti\`ere d'un voisinage tubulaire de ces arcs, et varie \'egalement de~$\pi$ \`a $2\pi$ quand on se d\'eplace de l'ext\'erieur vers l'int\'erieur sur cette fronti\`ere. 
Il reste donc \`a prolonger~$\theta$ sur $n$ cylindres disjoints, sachant que, sur le bord de ces cylindres, $\theta$ varie de~$\pi$ \`a~$2\pi$ comme indiqu\'e sur la gauche de la figure~\ref{F:cylindre}. 
On le fait \`a l'aide de la fonction hauteur indiqu\'ee \`a droite de la figure~\ref{F:cylindre}. 
De m\^eme, on prolonge~$\theta$ avec des valeurs dans l'intervalle~$[0, \pi]$ dans le demi-plan sup\'erieur au voisinage des arcs~$[x_2, x_3],\ldots, [x_{2n}x_1]$. 

\begin{figure}[htb]
	\begin{center}
	\includegraphics[scale=0.6]{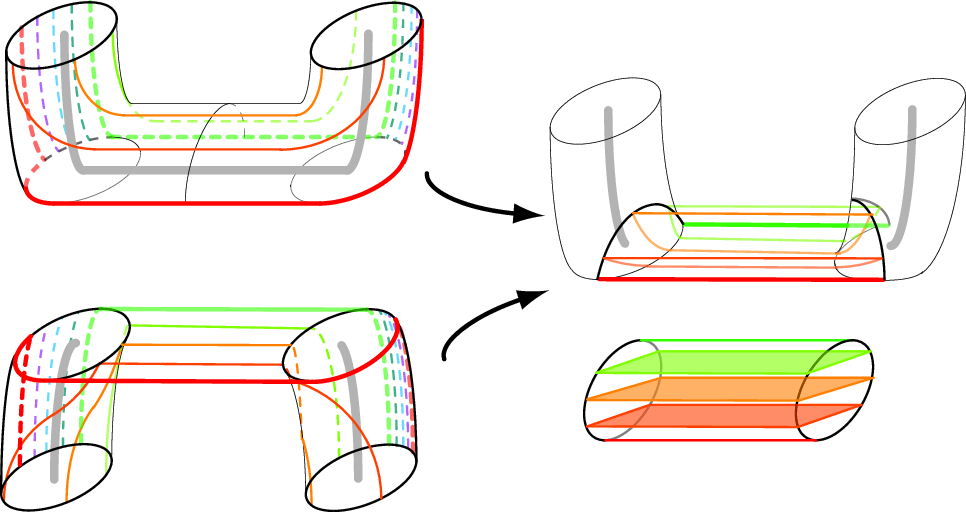}
	\end{center}
	\caption{\small \sf Compl\'etion de la fibration d'un somme de Murasugi \`a partir de la r\'eunion des fibrations associ\'ees \`a chacun des demi-espaces. 
\`A gauche, la fibration au voisinage d'un arc de type $]x_i, x_{i+1}[$ pour chacun des deux n\oe uds $K_1$ (en bas) et $K_2$ (en haut). 
Les arcs repr\'esentent les niveaux de la fonction $\theta$, et plus pr\'ecis\'ement l'intersection des fibres avec un voisinage tubulaire de $K_1$ et $K_2$ respectivement,  sur la partie du bas, on a \'egalement ajout\'e les niveaux de $\theta$ dans le plan~$\Pi$ de recollement. 
En haut \`a droite, la partie sur laquelle il faut d\'efinir de nouvelles valeurs pour la fonction $\theta$, qui est un cylindre. 
Les conditions au bord sont alors fix\'ees. 
En bas \`a droite, un feuilletage d'un cylindre remplissant ces contraintes de bord.}
	\label{F:cylindre}
\end{figure}

La fonction $\theta$ est alors d\'efinie sur tout $\Sph^3\smallsetminus N(K_1\Murasomme K_2)$, et \`a valeurs dans $\Sph^1$.

Pour v\'erifier que $\theta$ est une fibration, \'etablissons que la fonction $\theta$ n'a pas de singularit\'e. 
Comme les fonctions $\theta_1$ et $\theta_2$ sont tangentes au plan~$\Pi$ de recollement hors de $N(\partial P)$, et sans singularit\'e, et comme $\theta$ co\"\i  ncide partout hors de $N(P)$ avec l'une de ces deux fonctions, la fonction $\theta$ y est sans singularit\'e. 
Au voisinnage des arcs $]x_ix_{i+1}[$, c'est-\`a-dire \`a l'int\'erieur de la diff\'erence sym\'etrique~$N(K_1)\Delta N(K_2)$, le prolongement de $\theta$ a \'et\'e justement choisi pour ne pas faire appara\^\i tre de singularit\'e. 

Pour compl\'eter la d\'emonstration, il reste \`a v\'erifier que $\theta$ est bien un argument au voisinage des sommets~$x_1, \ldots x-{2n}$ du polygone~$P$. 
La figure~\ref{F:collagefibration} montre que la fonction $\theta$  y est bien de la forme d\'esir\'ee. 
\end{proof}

On applique maintenant aux tresses les r\'esultats g\'en\'eraux qui pr\'ec\`edent.

\begin{lemma}\label{T:ColonneFibre}
Pour tout entier $n$ strictement positif, la cl\^oture de la tresse $\sigma_1^n$ est un entrelacs fibr\'e.
\end{lemma}

\begin{figure}[htbp]
	\begin{center}
	\includegraphics[width=0.7\textwidth]{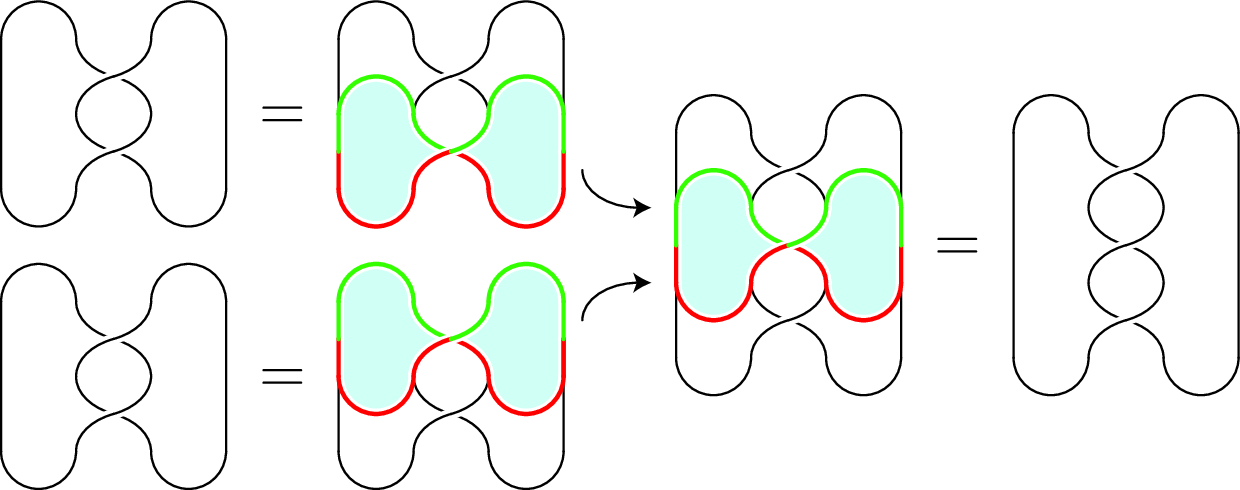}
	\end{center}
	\caption{\small \sf Construction it\'erative d'une surface de Seifert pour la cl\^oture de $\sigma_1^n$.}
	\label{F:TrefleFibre}
\end{figure}

\begin{proof}
Le cas $n=1$ correspond au n\oe ud trivial, tandis que le cas $n=2$ correspond \`a l'entrelacs de Hopf \`a deux composantes.
Sinon, on obtient la surface de Seifert standard $\surf^n$ pour la cl\^oture de $\sigma_1^n$ en collant les surfaces standards $\surf^{n-1}$ et $\surf^2$ le long du quadrilat\`ere gris\'e sur la figure~\ref{F:TrefleFibre}. 
Le lemme~\ref{T:sommefibre} assure alors qu'\`a chaque it\'eration, l'entrelacs est fibr\'e avec $\surf^n$ pour fibre.
\end{proof}

On est maintenant pr\^et pour d\'emontrer que la cl\^oture de toute tresse positive est un entrelacs fibr\'e.

\begin{figure}[htbp]
	\begin{center}
	\includegraphics[width=0.45\textwidth]{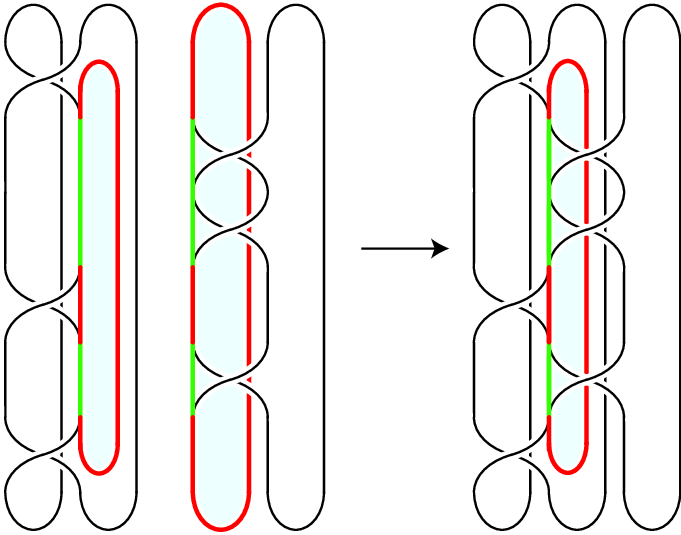}
	\end{center}
	\caption{\small \sf Collage de deux tresses l'une \`a c\^ot\'e de l'autre. En bleu, le polygone de collage.}
	\label{F:CollageLateral}
\end{figure}

\begin{proof}[du th\'eor\`eme~\ref{T:fibrationpositive}] 
D'apr\`es le lemme~\ref{T:ColonneFibre}, la cl\^oture d'une tresse-colonne $\sigma_i^n$ est un entrelacs fibr\'e. 
La figure~\ref{F:CollageLateral} montre comment juxtaposer une \`a une les surfaces standards associ\'ees \`a chaque colonne par somme de Murasugi, afin d'ainsi obtenir la surface standard associ\'ee \`a la cl\^oture de n'importe quelle tresse positive.
D'apr\`es le lemme~\ref{T:sommefibre}, l'entrelacs obtenu est fibr\'e, avec la surface standard pour fibre.
\end{proof}

\begin{figure}[htb]
	\begin{center}
	\begin{picture}(64,74)(3,0)
	\put(4,0){\includegraphics*[scale=0.4]{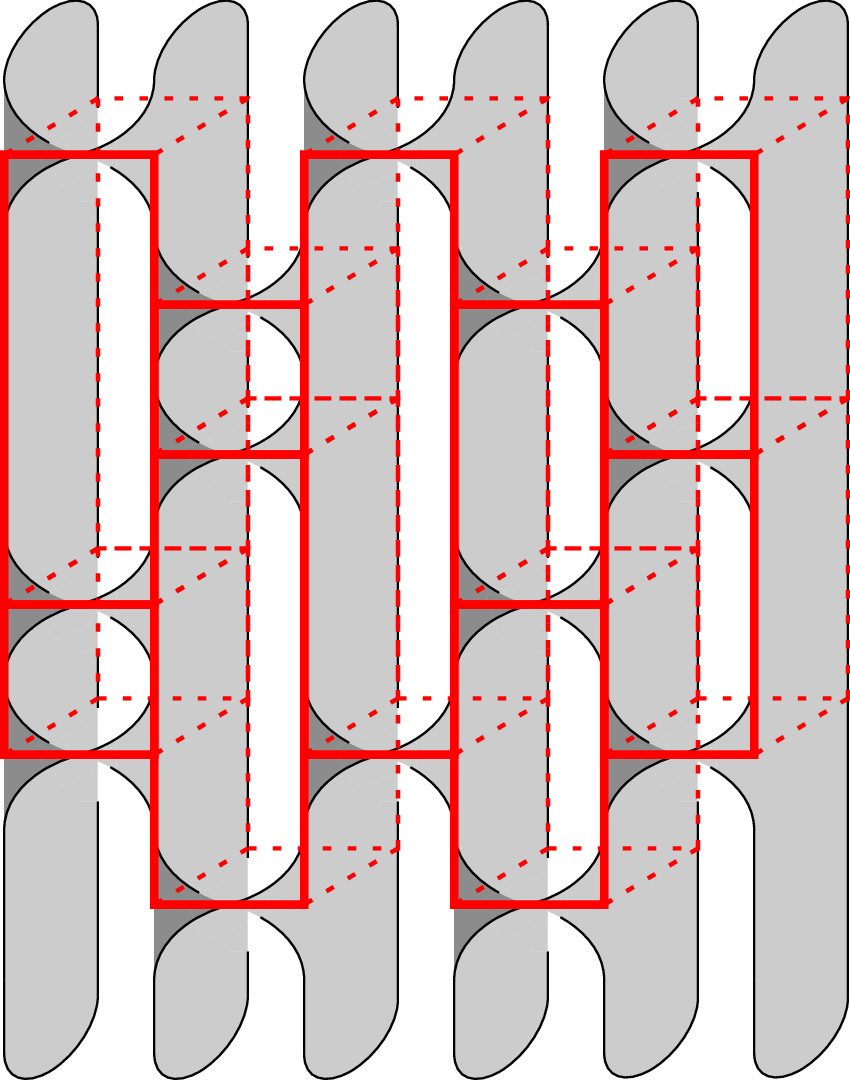}}
	\put(11,47){$1$}
	\put(11,26){$2$} 
	\put(21.5,47){$3$}
	\put(21.5,26){$4$} 
	\put(31.5,42){$5$}
	\put(41.8,42){$6$} 
	\put(41.8,21){$7$} 
	\put(52,52){$8$} 
	\put(52,32){$9$} 
	\end{picture}
	\caption{\small \sf Construction it\'erative de la fibration associ\'ee \`a la cl\^oture d'une tresse positive. 
	Chaque brique correspond \`a un entrelacs de Hopf \`a deux composantes. 
	En les empilant dans l'ordre prescrit, on s'assure que la somme de Murasugi fibre \`a chaque \'etape de la construction.}
	\label{F:constructionbrique}
	\end{center}
\end{figure}

\begin{remark}
On a utilis\'e la positivit\'e de $b$ de mani\`ere cruciale en affirmant que chaque colonne \'el\'ementaire est fibr\'ee avec la surface standard pour fibre. 
En effet, si la tresse $b$ a un croisement positif et un croisement n\'egatif cons\'ecutifs dans une m\^eme colonne, celle-ci peut n\'eanmoins \^etre fibr\'ee, mais la surface standard ne sera pas une fibre. 
Le th\'eor\`eme~\ref{T:sommefibre} ne s'applique alors plus dans la construction par r\'ecurrence, et donc la preuve s'effondre.

Par contre, l'argument peut \^etre \'etendu au cas d'une tresse {\it homog\`ene}, d\'efinie comme une tresse poss\'edant une d\'ecomposition dans laquelle, pour chaque $i$, l'un au plus des g\'en\'erateurs $\sigma_i$ ou $\sigma_i^{-1}$ appara\^\i t.
Ceci montre par exemple que le n\oe ud de huit, cl\^oture de la tresse $\sigma_1\sigma_2^{-1}\sigma_1\sigma_2^{-1}$, est fibr\'e.
\end{remark}


\subsection{Genre}
\label{S:Genre}

Dans la section pr\'ec\'edente, on a montr\'e qu'un n\oe ud qui est cl\^oture d'une tresse positive est fibr\'e. 
La technique utilis\'ee a \'egalement permis de construire une fibre pour la fibration associ\'ee. 
Dans cette section, nous utilisons ces informations pour calculer  le genre des n\oe uds de Lorenz et en d\'eduire un crit\`ere simple permettant de montrer qu'un n{\oe}ud donn\'e n'est pas un n{\oe}ud de Lorenz.

\begin{proposition}
\label{T:genreminimal}
Soit $b$ une tresse positive, $K$ sa cl\^oture et $\surf$ la surface de Seifert standard.
Alors $\surf$ est une surface de Seifert de genre minimal pour~$K$.
\end{proposition}

\begin{proof}
Soit $g$ le genre de la surface~$\surf$. 
Reprenons les notations de la section pr\'ec\'edente.
On a vu que le n\oe ud $K$ admet une fibration d\'efinie par une fonction $\theta$ telle que $\surf = \theta^{-1}(0)$. 
Comme $\theta$ d\'efinit une fibration de $\Sph^3\smallsetminus K$, le rev\^etement cyclique infini du compl\'ementaire de $K$ est isomorphe \`a $\theta^{-1}(0) \times \reels$, et donc son groupe fondamental est un groupe libre de rang~$2g$. 
De plus, le lacet horizontal $\partial\Sigma(0)\times \{0\}$ dans $\Sigma\times \reels$ est le produit de $g$ commutateurs dans~$\pi_1(\Sph^3\smallsetminus K)$, et pas moins.

Soit $\Sigma'$ une surface de Seifert quelconque pour $K$ de genre $g'$. 
En passant au rev\^etement universel, on plonge $\Sigma'$ dans $\Sigma \times \reels$, et en particulier son groupe fondamental $\pi_1(\Sigma')$ se plonge dans~$\pi_1(\Sigma \times \reels)$. 
Le lacet $\partial\Sigma'$ est, dans $\pi_1(\Sigma')$, le produit de $g'$ commutateurs.
Or il est isotope dans $\Sigma \times \reels$ au lacet $\partial\Sigma\times \{0\}$, qui est le produit de $g$ commutateurs et pas moins.
On en d\'eduit l'in\'egalit\'e $g\le g'$, comme~voulu.
\end{proof}

\begin{proposition}
\label{T:genre}
Soit~$\gamma$ une orbite p\'eriodique de p\'eriode~$n$ du flot de Lorenz, $K$ le n\oe ud associ\'e et $c$~le nombre de croisements de la tresse de Lorenz associ\'ee. 
Alors le genre~$g$ de~$K$ est donn\'e par la formule $2g=n-c$.
\end{proposition}

\begin{proof}
La fibre $\theta^{-1}(0)$ exhib\'ee dans la preuve du th\'eor\`eme~\ref{T:fibrationpositive} est une surface de Seifert de genre minimal. 
Or, il s'agit de la surface de Seifert standard, qui par cons\'equent est de genre minimal. 
On v\'erifie facilement que sa caract\'eristique d'Euler est $1+c-n$ et qu'elle n'a qu'une composante de bord\; par cons\'equent son genre est~$(n-c)/2$.
\end{proof}

On a vu dans la proposition~\ref{C:Infinite} que tout n{\oe}ud de Lorenz est r\'ealis\'e par une infinit\'e d'orbites du flot de Lorenz, et, par cons\'equent, est associ\'e \`a une infinit\'e de tresses de Lorenz. 
La proposition~\ref{T:genre} implique que, pour toutes ces tresses, la quantit\'e~$n - c$, diff\'erence entre le nombre de brins et le nombre de croisements, prend la m\^eme valeur.
La proposition~\ref{T:genre} se reformule directement en termes de diagrammes de Young.

\begin{corollary}
\label{C:GenreTableau}
Le genre d'un n\oe ud de Lorenz associ\'e \`a un diagramme de Young est la moiti\'e du nombre de cases de ce diagramme.
\end{corollary}

Nous avons vu avec la proposition~\ref{T:transpose} que deux diagrammes de Young distincts peuvent coder le m\^eme n\oe ud de Lorenz. 
Nous pouvons maintenant borner sup\'erieurement le d\'efaut d'injectivit\'e de ce codage.

\begin{proposition}
\label{T:finitude}
Pour tout n\oe ud de Lorenz~$K$, il n'y a qu'un nombre fini de diagrammes de Young standards associ\'es \`a~$K$.
\end{proposition}

\begin{proof}
Si $K$ a pour genre~$g$, tout diagramme de Young associ\'e~$K$ a $2g$ cases en vertu du corollaire~\ref{C:GenreTableau}. 
Or il n'existe qu'un nombre fini de diagrammes de Young ayant $2g$ cases.
\end{proof}

Ce r\'esultat fournit un moyen de d\'emontrer qu'un n\oe ud donn\'e n'est pas un n\oe ud de Lorenz. 
Par exemple, le n\oe ud de huit n'est pas un n\oe ud de Lorenz, car son genre est~1, et les seuls diagrammes de Young \`a deux cases sont tous deux associ\'es au n\oe ud de tr\`efle. 

\begin{proposition}
\label{T:16croisements}
Parmi tous les n\oe uds admettant un repr\'esentant planaire ayant au plus seize croisements, exactement vingt sont des n\oe uds de Lorenz.
\end{proposition}

\begin{proof}
Dans~\cite{A:atlas}, nous avons d\'etermin\'e les polyn\^omes d'Alexander et de Jones de tous les n\oe uds de Lorenz de p\'eriode au plus~$21$. 
Cette liste inclut en particulier tous les n\oe uds de Lorenz de genre au plus~$9$. 
Or elle ne compte que vingt n\oe uds ayant au plus seize croisements. 
Si on s'int\'eresse \`a la question de savoir quels n\oe uds ayant au plus seize croisements sont de Lorenz, cette liste est suffisante puisqu'un n\oe ud ayant au plus $2g+1$ croisements est de genre au plus~$g$. 
Par cons\'equent, les n\oe uds ayant au plus seize croisements sont de genre au plus~$8$, et donc ont \'et\'e recens\'es dans la liste de~\cite{A:atlas}.
\end{proof}

Signalons une derni\`ere propri\'et\'e g\'eom\'etrique des n\oe uds de Lorenz.

\begin{proposition}[\cite{A:B-W}]
	La signature de tout n\oe ud de Lorenz est strictement positive.
\end{proposition}

\begin{proof}
	Un th\'eor\`eme de L.\,Rudolph~\cite{A:Rudolph} affirme que la signature de la cl\^oture de toute tresse positive est positive. 
	Il s'applique aux tresses de Lorenz.
\end{proof}

Comme la signature de l'image-miroir d'un n\oe ud est l'oppos\'e de la signature du n\oe ud, on d\'eduit:

\begin{corollary}
	\label{T:nonamphichiraux}
	Aucun n\oe ud de Lorenz n'est isotope \`a son image-miroir.
\end{corollary}


\section{L'indice de tresse}
\label{S:IndiceTresse}

Nous poursuivons maintenant l'\'etude topologique des n{\oe}uds de Lorenz  en d\'eterminant de fa\c{c}on d\'etaill\'ee leur indice de tresse, c'est-\`{a}-dire le plus petit nombre de brins d'une tresse repr\'{e}sentant ce n\oe ud. 
Pour ce faire, nous introduisons dans la section~\ref{S:TresseBW} une nouvelle famille de tresses, dite de Birman-Williams, de sorte qu'une orbite de Lorenz donn\'ee est cl\^oture de la tresse de Birman-Williams associ\'ee. 
Le nombre de brins d'une tresse de Birman-Williams est, par d\'efinition, le \trip\ de l'orbite de Lorenz associ\'ee, une donn\'ee combinatoire qui se lit par exemple sur les mots de Lyndon. 
Nous montrons ensuite que cette tresse minimise le nombre de brins parmi les tresses ayant la m\^eme cl\^oture.
Pour ce faire, nous utilisons un r\'esultat g\'en\'eral sur l'indice de tresse d'un n\oe ud qui est cl\^oture d'une tresse positive, d\^u \`a J.\,Franks, H.\,Morton et R.\,Williams (th\'eor\`eme~\ref{T:F-W})\footnote{Le nom de H.\,Morton n'appara\^it pas dans les auteurs de l'article~\cite{A:F-W}, mais il y est mentionn\'e que ce dernier est aussi l'auteur d'une d\'emonstration du th\'eor\`eme~\ref{T:F-W}.}. 
Nous donnons une d\'emonstration d\'etaill\'ee de ce r\'esultat, \`a la fois pour son int\'er\^et propre et pour offrir une version plus d\'etaill\'ee que celle de l'article original~\cite{A:F-W}. 
La section~\ref{S:Homfly} rappelle la construction r\'ecursive du polyn\^ome HOMFLY \`a partir des relations d'\'echeveau et introduit le formalisme des arbres de calcul positifs. 
Dans la section~\ref{S:Variante}, nous introduisons une variante du polyn\^ome HOMFLY qui est bien adapt\'ee au cas des n{\oe}uds qui sont cl\^otures de tresses positives. 
Ceci nous m\`ene dans la section~\ref{S:Application} \`a la d\'emonstration du th\'eor\`eme~\ref{T:F-W} et \`a son application au calcul de l'indice de tresse des n{\oe}uds de Lorenz.
Enfin, dans la section~\ref{S:B-K}, nous d\'etaillons le lien entre entrelacs de Lorenz et $T$-entrelacs, une autre famille d'entrelacs admettant une description combinatoire simple.


\subsection{Tresse de Birman--Williams et indice de tresse}
\label{S:TresseBW}

\begin{figure}[htbp]
	\begin{center}
	\includegraphics[height=0.8\textheight]{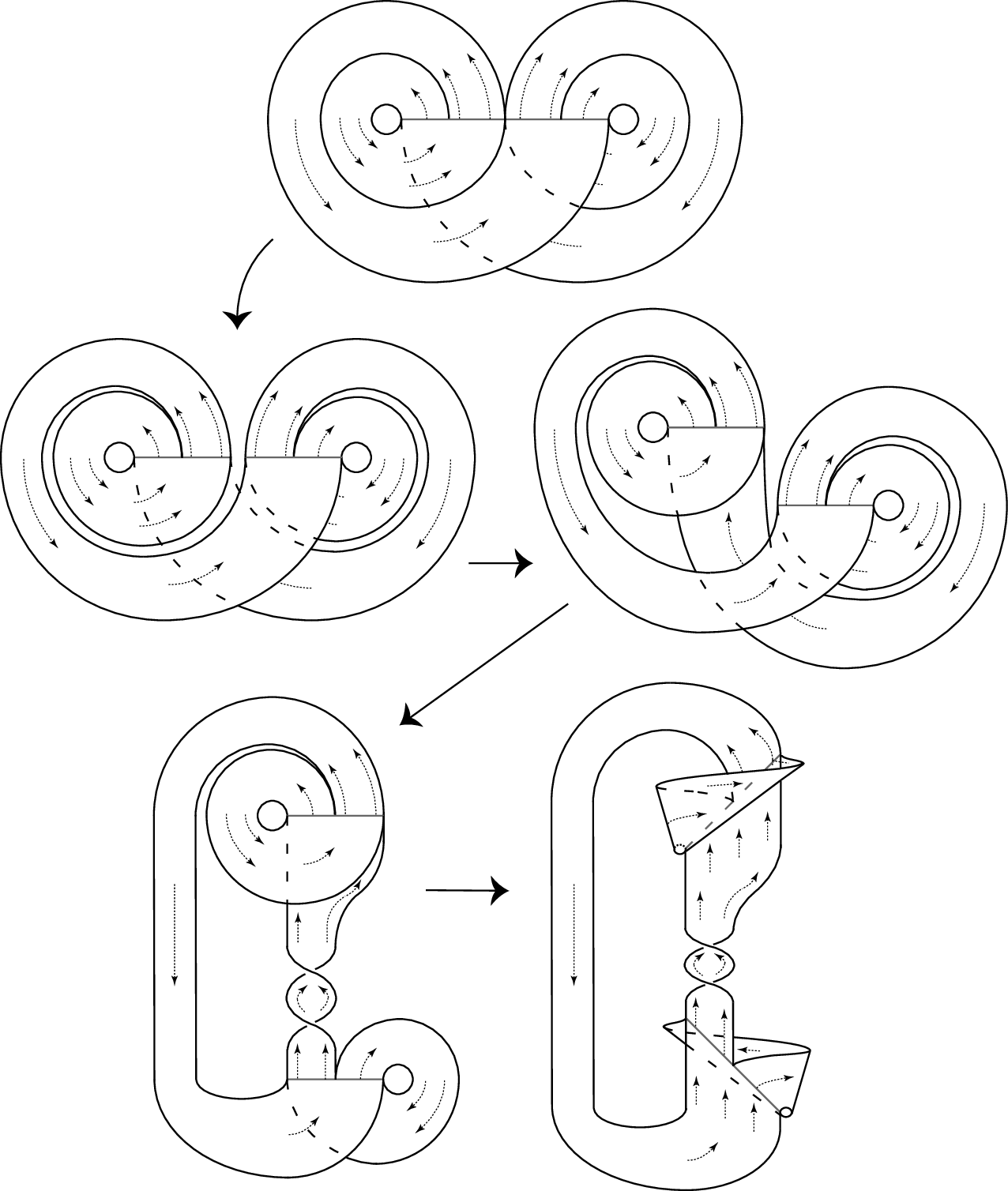}
	\end{center}
	\caption{\small \sf Le patron de Lorenz d\'{e}form\'{e} de sorte que le flot soit ascendant sur la branche droite et descendant sur la branche gauche.}
	\label{F:deformation}
\end{figure}

Dans le th\'eor\`eme~\ref{T:tresseLorenz}, on a associ\'e \`a toute orbite~$\gamma$ du flot de Lorenz une tresse privil\'egi\'ee, dite tresse de Lorenz, dont la cl\^oture repr\'esente le n\oe ud~$K$ associ\'e \`a~$\gamma$. 
En d\'{e}formant le patron comme sur la figure~\ref{F:deformation}, on fait appara\^\i tre une autre tresse dont la cl\^oture est~$K$. 
Cette tresse, dite \emph{de Birman--Williams}, est, comme la tresse de Lorenz, une tresse positive. 
Par contre, ce n'est pas une tresse de permutation. 
Elle fournit un repr\'esentant plus \'econome que la tresse de Lorenz  dans la mesure o\`u elle met en jeu moins de brins et moins de croisements. 

Pour d\'{e}crire les tresses de Birman--Williams, on d\'{e}finit d'abord le \trip\ d'une orbite de Lorenz.

\begin{definition}
\label{D:trip}
On appelle \emph{\trip} d'une orbite p\'eriodique du flot de Lorenz le nombre d'occurrences du mot~$\xx\yy$ dans le code de Lyndon de cette orbite.
\end{definition}

En termes de tresses de Lorenz (voir la figure~\ref{F:permutation}), le \trip\ correspond au nombre de brins qui traversent de gauche \`a droite la ligne imaginaire s\'{e}parant les $p$ premiers brins des $n-p$ derniers. 
En termes de diagrammes de Young (voir la figure~\ref{F:Young}), le \trip\ est le plus grand entier~$t$ tel qu'on puisse inscrire un carr\'e $(t-1)\times (t-1)$ dans le diagramme de Young correspondant \`a l'orbite.

On rappelle que, pour~$t$ entier positif, $\Delta_t^2$ d\'esigne la tresse $(\sigma_1\dots\sigma_{t-1})^t$ correspondant \`{a} un tour complet de $t$~brins sur eux-m\^{e}mes. 
Par ailleurs, une tresse positive~$b$ est dite \emph{multiple} de~$\Delta_t^2$ s'il existe une tresse positive~$b'$ v\'erifiant $b = \Delta_t^2 b'$.

\begin{theorem}[\cite{A:B-W} et consulter l'erratum]
	\label{T:tresseB-W}
	Soit $\gamma$ une orbite p\'{e}riodique du flot de Lorenz de \trip~$t$, et $n_i$ et $m_i$ d\'efinis par les \'equations~\eqref{E:ni} et~\eqref{E:mi}.
	Soit $K$ le n{\oe}ud associ\'{e}, et $b$ la tresse \`{a} $t$ brins
	\begin{equation}
	\label{E:B-W}
	\Delta_t^2 \ \prod_{i=1}^{t-1}(\sigma_1\sigma_2\dots\sigma_i)^{n_i} 
	\prod_{i=t-1}^{1}(\sigma_{t-1}\dots\sigma_{i})^{m_{t-i}}.
	\end{equation}
	Alors la cl\^oture de la tresse $b$ est un repr\'esentant de $K$.
\end{theorem}

\begin{likeproof}[Principe de la d\'emonstration]
	On d\'ecrit l'image de la tresse de Lorenz de la figure~\ref{F:permutation}, quand on la
d\'{e}forme selon le sch\'{e}ma indiqu\'{e} sur la figure~\ref{F:deformation}.
\end{likeproof}

La tresse d\'efinie par la formule~\eqref{E:B-W} est appel\'{e}e {\it tresse de Birman-Williams} pour l'orbite~$\gamma$. 

\begin{example}
	\label{X:BirmanWilliams}
	Revenons \`{a} l'orbite partant du point $5/31$. Le mot de Lyndon associ\'{e} est $\xx\xx\yy\xx\yy$, qui contient deux fois $\xx\yy$. 
	Le \trip\ vaut donc~2. 
	La permutation de Lorenz est $(13524)$, le seul indice non nul dans la tresse de Birman-Williams est alors $n_1$ qui vaut 1,  et donc la tresse de Birman-Williams associ\'{e}e est $\sigma_1^3$, dont la cl\^oture est bien un n{\oe}ud de tr\`efle. 
	Pour comparaison, la tresse de Lorenz correspondante est $\sigma_2\sigma_1\sigma_3\sigma_2\sigma_4\sigma_3$.
\end{example}

Notre but dans la suite de cette partie est de d\'emontrer le r\'esultat suivant.

\begin{theorem}
	\label{T:Trip}
	L'indice de tresse d'un n\oe ud de Lorenz est \'egal au \trip~de toute orbite du flot de Lorenz le r\'ealisant.
\end{theorem}

Eu \'egard au th\'eor\`eme~\ref{T:tresseB-W}, le th\'eor\`eme~\ref{T:Trip} est un corollaire direct du r\'esultat suivant, qui est celui que nous allons d\'emontrer dans la suite.

\begin{theorem}[\cite{A:F-W}]
	\label{T:F-W}
	Soit $K$ un n{\oe}ud qui est la cl\^oture d'une tresse positive $b$ \`{a} $t$ brins multiple de $\Delta_t^2$. 
	Alors l'indice de tresse de~$K$ est \'{e}gal \`{a}~$t$.
\end{theorem}


\subsection{Le polyn\^ome HOMFLY et les arbres de calcul}
\label{S:Homfly}

La d\'emonstration du th\'eor\`eme~\ref{T:F-W} va faire appel au polyn\^ome HOMFLY, aussi appel\'e HOMFLY-PT. 
Celui-ci est une extension \`a deux variables du polyn\^ome de Jones, et, comme ce dernier, il admet une d\'efinition simple en termes de relations d'\'echeveau pour les diagrammes planaires.

\begin{figure}[htb]
	\begin{center}
	\begin{picture}(45,16)(-23,0)
	\put(-20,2){\includegraphics[scale=0.5]{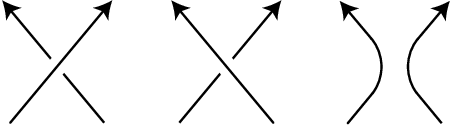}}
	\put(-18,14){$K_+$} 
	\put(-3,14){$K_-$} 
	\put(12,14){$K_0$}
	\put(-27,-2){$x~\HOM(K_+) + x^{-1}\HOM(K_-) + y~\HOM(K_0) = 0$}
	\end{picture}
	\end{center}
	\caption{\small \sf La relation d'\'{e}cheveau d\'{e}finissant le polyn\^ome HOMFLY.} 
	\label{F:skein}
\end{figure}

Soit $K$ un entrelacs orient\'e et $\mathcal K$ une projection planaire r\'{e}guli\`ere de $K$, c'est-\`{a}-dire ne poss\'edant qu'un nombre fini de points multiples et telle que ceux-ci soient tous des points doubles. 
Choisissons un point double $M$ de $\mathcal K$. 
D\'{e}finissons $K_+$ et $K_-$ comme les entrelacs admettant $\mathcal K$ pour projection, sauf \'eventuellement au point $M$, o\`{u} $K_+$ est projet\'e sur un croisement positif, et $K_-$ sur un croisement n\'egatif \footnote{Remarquons que $K$ co\"\i ncide avec $K_+$ ou avec $K_-$}. 
D\'efinissons $K_0$ comme l'entrelacs orient\'e o\`{u} le croisement en $M$ a \'{e}t\'{e} supprim\'{e} (voir figure~\ref{F:skein}). 
Alors, par d\'efinition, les polyn\^omes HOMFLY des trois entrelacs associ\'{e}s sont reli\'{e}s par la relation
\begin{equation}
	\label{E:Skein}
	x~\HOM_{K_+} + x^{-1}\HOM_{K_-} + y~\HOM_{K_0} = 0,
\end{equation}
dite \emph{relation d'\'{e}cheveau}. 
D\'{e}signons par $\rond$ le n{\oe}ud trivial. 
En imposant $\HOM_\rond=1$, le polyn\^ome $\HOM_K$ est alors d\'{e}fini de mani\`ere unique pour tout n\oe ud $K$. 
On montre que la relation d'\'{e}cheveau impose
$\HOM_{\rond^k} = \left(-(x+x^{-1})/y\right)^{k-1}$, 
o\`{u} $\rond^k$ d\'{e}signe la r\'eunion de $k$ n{\oe}uds triviaux deux \`{a} deux non enlac\'{e}s. 
Le polyn\^ome de Jones correspond \`a la sp\'ecialisation $y=1$.
\footnote{Une description plus pr\'ecise du polyn\^ome HOMFLY peut \^etre trouv\'ee par exemple dans~\cite[chapitre 15]{A:Lickorish}.}

Le calcul effectif du polyn\^ome HOMFLY d'un n{\oe}ud peut s'effectuer en introduisant un arbre de calcul dont les bifurcations successives correspondent aux points doubles o\`u la relation d'\'echeveau~\eqref{E:Skein} est invoqu\'ee. 
Un choix arbitraire des points doubles peut mener \`{a} une boucle, par exemple si on choisit deux fois de suite le m\^{e}me point. 
Par contre, on va voir que, pour autant qu'on change de point d'application \`{a} chaque \'{e}tape, le processus de calcul converge en un nombre fini d'\'{e}tapes.

Comme nous consid\'erons dans la suite des tresses qui ne poss\`edent pas toutes le m\^eme nombre de brins, il est commode de d\'efinir une \emph{tresse marqu\'ee} comme un couple~$(b, n)$ o\`u $b$ est une tresse \`a au plus $n$~brins, et  une \emph{tresse positive marqu\'ee} comme une tresse marqu\'ee $(b,n)$ o\`u $b$ est une tresse positive \`a au plus $n$~brins. 

\begin{definition} 
	\label{D:TransformationMarkov}
	On dit que deux tresses marqu\'ees~$(b, n)$ et $(b', n')$ sont reli\'ees par une \emph{transformation de Markov positive} si 

-- ou bien on a $n' = n$ et $b'$ est conjugu\'ee \`a~$b$ dans le groupe des tresses \`a $n$~brins, 

-- ou bien on a $n' = n+1$ et $b' = b \sigma_n$, 

-- ou bien, sym\'etriquement, on a $n = n'+1$ et $b = b' \sigma_n$.
\end{definition}

Remarquer que, si on note $c(b,n)$ la diff\'erence entre le nombre de croisements positifs dans~$(b,n)$ et le nombre de croisements n\'egatifs, alors la quantit\'e $c(b,n) - n$ est invariante par transformation de Markov positive.
La propri\'et\'e que nous utiliserons dans la suite est que, si deux tresses marqu\'ees sont reli\'ees par une suite de transformations de Markov positives, alors leur cl\^otures sont isotopes.\footnote{Dans la section~\ref{S:Stabilisation}, on a mentionn\'e les transformations de Markov g\'en\'erales, o\`u on peut
aussi passer de $(b,n)$ a $(b\sigma_n^{-1}, n+1)$. 
Rappelons que le th\'eor\`eme de Markov affirme que les cl\^otures de deux tresses marqu\'ees sont des entrelacs isotopes si et seulement si on peut passer de~$(b, n)$ \`a~$(b', n')$ par une suite finie de transformations de Markov (g\'en\'erales).}

Afin de suivre le sch\'{e}ma de calcul sugg\'{e}r\'{e} par la relation  d'\'{e}cheveau~\eqref{E:Skein}, J.\,Conway a introduit pour les tresses marqu\'ees une relation qui en est la contrepartie.

\begin{definition}
	\label{D:ScindementConway}
	Soit $(b, n)$ une tresse marqu\'ee. On dit que deux tresses marqu\'ees $(b_0, n)$ et $(b_-, n)$ sont obtenues par \emph{scindement de Conway} \`a partir de~$(b, n)$ s'il existe des tresses~$b_1, b_2$ et v\'erifiant
	\begin{equation}
		b = b_1 \sigma_i b_2, \quad
		b_0 = b_1 b_2, \quad\mbox{et}\quad
		b_- = b_1 \sigma_i^{-1} b_2.
	\end{equation}
\end{definition}

Si $(b,n)$ est une tresse marqu\'ee, on note $\HOM_{(b,n)}$ le polyn\^ome HOMFLY de sa cl\^oture. 
La relation d'\'echeveau~\eqref{E:Skein} implique que, si $(b_0, n)$ et $(b_-, n)$ sont obtenues par scindement de Conway \`a partir de~$(b, n)$, alors on a
\begin{equation}
	\label{E:SkeinTresse}
	\HOM_{(b,n)}(x,y) = -x^{-1}y \HOM_{(b_0,n)}(x,y) - x^{-2}\HOM_{(b_-,n)}(x,y).
\end{equation}
On peut alors d\'{e}finir  un arbre de calcul pour le polyn\^ome $\HOM$, qui mime en termes de tresses marqu\'ees le calcul du polyn\^ome HOMFLY de la cl\^oture.

\begin{definition}
	Soit $(b,n)$ une tresse marqu\'ee. 
	Un \emph{arbre de calcul positif pour $(b,n)$} est un arbre orient\'e, binaire, fini, dont chaque sommet est \'{e}tiquet\'{e} par une tresse marqu\'ee, dont chaque ar\^{e}te est \'{e}tiquet\'{e}e par un mon\^ome de $\relatifs[x, x^{-1}, y, y^{-1}]$,  et qui satisfait aux propri\'{e}t\'{e}s suivantes\:

$(i)$ La racine est \'etiquet\'ee $(b,n)$.

$(ii)$ Si $S$ est un sommet \'{e}tiquet\'{e} $(b_s, n_s)$ qui n'est pas une feuille, et si $S_g$ et $S_d$ sont ses fils, \'{e}tiquet\'{e}s respectivement $(b_g, n_g)$ et $(b_d, n_d)$, alors $(b_g, n_g)$ et $(b_d, n_d)$ s'obtiennent  \`{a} partir de $(b_s, n_s)$ en faisant d'abord une suite de transformations de Markov positives, puis un scindement de Conway.

$(iii)$ Si, de plus, les ar\^{e}tes reliant $S$ \`{a} $S_g$ et $S_d$ sont  \'{e}tiquet\'{e}es par les mon\^omes $\monome_g$ et $\monome_d$ respectivement, alors on a $\HOM_{(b_s,n_s)} = \monome_g\HOM_{(b_g,n_g)} + \monome_d\HOM_{(b_d,n_d)}$.

$(iv)$ Si $f$ est une feuille \'{e}tiquet\'{e}e $(b_f,n_f)$, alors $b_f$ est la tresse triviale \`a $n_f$ brins.\end{definition}

	S'il existe un arbre de calcul positif $\arbre$ pour une tresse marqu\'ee $(b,n)$, il devrait \^etre clair que cet arbre fournit une strat\'{e}gie pour d\'{e}terminer le polyn\^ome HOMFLY de la cl\^oture de $(b,n)$, en it\'{e}rant la relation d'\'{e}cheveau~\eqref{E:SkeinTresse}. 
	En effet, notons $\feuilles$ l'ensemble des feuilles de $\arbre$, et, pour une feuille~$f$ dans
$\feuilles$, notons $\monome_f$ le produit des mon\^omes rencontr\'{e}s sur les ar\^{e}tes de $\arbre$ le long du chemin allant de la racine \`{a} $f$. 
	Alors, par construction, on a
	\[\HOM_{(b,n)}(x,y) = \sum_{f \in \feuilles}\monome_f \cdot \left(-\frac{x+x^{-1}}{y}\right)^{n_f-1}.\]
	Par contre, il n'est pas \'{e}vident {\it a priori} qu'un tel arbre de calcul existe pour toute tresse. 
	Nous allons montrer maintenant que c'est bien le cas. 
	\`A partir de maintenant, on ne travaille plus qu'avec des tresses positives
\footnote{Si on autorise \'egalement des scindements de Conway {\it n\'egatifs}, c'est-\`a-dire o\`u les r\^oles de $\sigma_i$ et $\sigma_i^{-1}$ sont \'echang\'es, les arbres de calculs existent en fait pour toute tresse marqu\'ee. 
La d\'emonstration est un peu plus longue, mais repose sur les m\^emes arguments~\cite{A:F-W}.}.
	Commen\c cons par un r\'{e}sultat technique sur les transformations de Markov positives. 

\begin{lemma}
	\label{L:TressePositive}
	Soit $(b,n)$ une tresse positive marqu\'ee. 
	Alors, par une suite de transformations de Markov positives, on peut transformer $(b,n)$ en $(b',n')$, o\`{u} $b'$ est 
soit la tresse triviale sur $n'$ brins, soit une tresse positive qui est multiple de~$\sigma_i^2$ pour au moins un entier~$i$.
\end{lemma}

\begin{proof}
	D\'efinissons le \emph{poids} $\poids$ d'un mot de tresse positif par $\poids(\sigma_i)=i$ et $\poids(w_1w_2)=\poids(w_1)+\poids(w_2)$.
	Nous allons prouver le r\'{e}sultat par r\'{e}currence sur le poids d'un mot de tresse $\mot$ repr\'esentant~$b$. 
	Le seul mot de tresse de poids nul est le mot vide, correspondant \`a la tresse triviale, pour laquelle le r\'{e}sultat est vrai. 
	Soit $m$ le plus grand entier tel que $\sigma_m$ appara\^\i t dans le mot $\mot$. 

	Si $\sigma_m$ n'appara\^\i t qu'une fois dans $\mot$, alors $\mot$ est de la forme $\mot_0 \sigma_m \mot_4$, et donc $b$ est conjugu\'{e}e \`{a} la tresse d\'ecrite par le mot $\mot_4 \mot_0 \sigma_m$.
	Une transformation de Markov positive supprime alors le $\sigma_m$ final. 
	Comme $\poids(\mot_4 \mot_0) = \poids(\mot) - m$, par hypoth\`ese de r\'{e}currence, le r\'{e}sultat est vrai pour la tresse repr\'esent\'ee par le mot $\mot_4 \mot_0$, et donc pour la tresse~$b$.

	Sinon, on peut supposer que $\sigma_m$ appara\^it au moins deux fois dans $\mot$.
	Alors le mot $\mot$ se d\'ecompose en $\mot_0 \sigma_m \mot_2 \sigma_m \mot_4$ de sorte que $\mot_2$ ne contienne aucun $\sigma_m$. 
	Il y a alors trois cas possibles:

{\it Premier cas\:} $\sigma_{m-1}$ n'appara\^it pas dans $\mot_2$. 
	La tresse repr\'esent\'ee par le mot $\sigma_m \mot_2 \sigma_m$ est \'egalement repr\'esent\'ee par $\mot_2 \sigma_m^2$, et par cons\'equent $b$ est conjugu\'{e}e \`{a} la tresse repr\'esent\'ee par~$\mot_4 \mot_0 \mot_2 \sigma_m^2$.

{\it Deuxi\`eme cas\:} $\sigma_{m-1}$ appara\^\i t une fois exactement dans $\mot_2$. 
	Alors $\mot_2$ se d\'ecompose en $\mot_1 \sigma_{m-1} \mot_3$, et la tresse $b$, repr\'esent\'ee par le mot $\mot_0 \sigma_m \mot_1 \sigma_{m-1} \mot_3 \sigma_m \mot_4$, est \'egalement repr\'esent\'ee par le mot  $\mot_0 \mot_1 \sigma_m \sigma_{m-1} \sigma_m \mot_3 \mot_4$, dont le poids est aussi $\poids(\mot)$. 
	Mais, d'apr\`es la relation de tresse $\sigma_{m-1}\sigma_{m}\sigma_{m-1}=\sigma_{m}\sigma_{m-1}\sigma_{m}$, $b$ est aussi repr\'esent\'ee par le mot $\mot_0 \mot_1 \sigma_{m-1} \sigma_m \sigma_{m-1} \mot_3 \mot_4$, dont le poids est $\poids(\mot) - 1$.
	L'hypoth\`ese de r\'{e}currence s'applique alors.

{\it Troisi\`eme cas\:} $\sigma_{m-1}$ appara\^\i t au moins deux fois dans $\mot_2$. 
	Alors, $\mot_2$ se d\'ecompose sous la forme $\mot_1 \sigma_{m-1} \mot_2' \sigma_{m-1} \mot_3$, avec $\mot_2'$ ne contenant aucun $\sigma_{m-1}$.
	On r\'eit\`ere le processus \`{a} partir du mot $\mot_2'$, \`{a} la recherche de la lettre $\sigma_{m-2}$. 
	\`A chaque it\'eration, soit on trouve un mot repr\'esentant la tresse $b$ et contenant un facteur $\sigma_i^2$, soit le poids d\'ecro\^\i t, soit l'indice $m$ du g\'en\'erateur recherch\'e d\'ecro\^\i t.  Cette derni\`ere \'eventualit\'e ne peut se produire une infinit\'e de fois, par cons\'{e}quent le processus s'ach\`eve.
\end{proof}

Le r\'esultat suivant justifie l'introduction des arbres de calcul positifs, puisqu'il donne un moyen de calculer le polyn\^ome HOMFLY des tresses positives. 
On d\'efinit la longueur~$\longueur(b)$ d'une tresse positive~$b$ comme la longueur de n'importe quel mot de tresse positif repr\'esentant~$b$.

\begin{theorem}[\cite{A:F-W}]
	\label{T:ArbreCalcul}
	Il existe un arbre de calcul positif pour toute tresse positive marqu\'ee.
\end{theorem}

\begin{proof}
	Soit~$(b,n)$ une tresse positive marqu\'ee.
	La preuve se fait par r\'{e}currence sur la longueur $\longueur(b)$. 
	Si $b$ est de longueur nulle, alors la tresse~$b$ est triviale, et l'arbre n'ayant qu'un sommet \'etiquet\'e $(b,n)$ convient.

	Sinon, on utilise le lemme~\ref{L:TressePositive}\: la tresse positive marqu\'ee~$(b, n)$ peut \^{e}tre convertie par une suite de transformations de Markov positives en $(b',n')$, o\`{u} $b'$ est soit la tresse triviale, soit de la forme $b'' \sigma_i^2$ avec $b''$ positive. 
	Dans le premier cas, la tresse $b'$ est elle-m\^eme triviale et l'arbre n'ayant qu'un sommet \'etiquet\'e $(b,n)$ convient. 
	Dans le second cas, on utilise un scindement de Conway sur un des deux croisements $\sigma_i$, et les deux fils sont alors \'{e}tiquet\'{e}s $(b''\sigma_i, n')$ et $(b'', n')$. 
	Par construction, les tresses positives $b'' \sigma_i$ et $b''$ ont une longueur strictement inf\'erieure \`a celle de~$b$, donc, par hypoth\`ese de r\'ecurrence, on peut compl\'{e}ter l'arbre de calcul sous les sommets \'etiquet\'es $(b''\sigma_i, n')$ et $(b'', n')$, et, de l\`a, sous $(b, n)$.
\end{proof}

\subsection{Une variante du polyn\^ome HOMFLY}
\label{S:Variante}

Dans la d\'emonstration du th\'eor\`eme~\ref{T:ArbreCalcul}, toutes les branches gauches de l'arbre construit sont \'{e}tiquet\'{e}es $-x^{-1}y$ et toutes les branches droites $-x^{-2}$.
Introduisons alors une variante du polyn\^ome HOMFLY qui se comporte bien vis-\`{a}-vis des tresses marqu\'ees et des arbres de calcul positifs. 

\begin{proposition}
	\label{T:RCT}
	$(i)$ Pour toute tresse positive marqu\'ee~$(b, n)$, il existe un polyn\^ome $J_{(b,n)}(R,C,T)$ \`{a} coefficients entiers positifs v\'erifiant
	\begin{equation}
		\HOM_{(b,n)}(x, y) = J_{(b,n)}\left(-x^{-2}, -x^{-1}y, -\frac{x+x^{-1}}{y}\right).
	\end{equation}

	$(ii)$ Si on donne \`{a} $R$, $C$ et $T$ les degr\'{e}s $1$, $2$ et $-1$ respectivement, alors $J_{(b,n)}(R,C,T)$ est homog\`ene de degr\'{e} $\longueur(b) - n + 1$,

	$(iii)$ Le polyn\^ome $J_{(b,n)}$ ne d\'{e}pend que de la classe d'isotopie de la cl\^oture de $(b,n)$.

	$(iv)$ Si $q$ est le nombre de composantes de la cl\^oture de $(b,n)$, alors on a $J_{(b,n)}(0, C,T) = C^pT^{q-1}$, o\`{u} $p$ est un entier.

	$(v)$ Le degr\'e en~$T$ de~$J_{(b,n)}(R,C,T)$ est strictement inf\'erieur \`a l'indice de tresse de~$b$.
\end{proposition}

\begin{proof}
	Nous reprenons les arbres de calcul positifs construits plus haut, mais en modifiant l'\'{e}tiquetage des ar\^{e}tes. 
	Soit $(b,n)$ une tresse positive marqu\'ee. 
	Soit $\arbre$ un arbre de calcul positif pour~$(b,n)$, construit selon le proc\'ed\'e d\'ecrit dans la d\'{e}monstration du th\'{e}or\`eme~\ref{T:ArbreCalcul}.
	Par construction, toutes les ar\^{e}tes allant vers un fils gauche de $\arbre$ sont \'{e}tiquet\'{e}es $-x^{-1}y$ et toutes les ar\^etes allant vers un fils droit sont \'etiquet\'ees $-x^{-2}$. 
	Rempla\c{c}ons ces \'{e}tiquettes par $C$ et $R$ respectivement. 
	Notons~$\feuilles$ l'ensemble des feuilles de $\arbre$.
	D\'{e}finissons alors $\monj_f(R,C)$ comme le mon\^ome produit de toutes les \'{e}tiquettes $R$ et $C$ rencontr\'{e}es le long du chemin allant de la racine \`{a} $f$. 
	Enfin, posons 
	\begin{equation}
		\label{E:defJ} 
		J_{(b,n)}(R,C,T)=\sum_{f \in \feuilles}\monj_f(R,C)\cdot T^{n_f-1}.
	\end{equation}

	{\it A priori}, la valeur de $J_{(b,n)}$ d\'{e}pend de l'arbre de calcul choisi. 
	C'est seulement lorsque nous aurons d\'emontr\'e le point $(iii)$ que nous en d\'{e}duirons que $J_{(b,n)}$ ne d\'epend que de la cl\^oture de~$(b, n)$. 

	Le point $(i)$ se d\'{e}duit imm\'{e}diatement de la d\'{e}finition de $J$, puisque, pour le d\'efinir, on a repris les arbres de calcul positifs introduits pour calculer le polyn\^ome HOMFLY, en rempla\c cant les $-x^{-2}$ par~$R$, les $-x^{-1}y$ par $C$, et les $-\frac{x+x^{-1}}{y}$ par $T$. 
	Or, pour toute feuille $f$ dans~$\feuilles$, le nombre~$n_f$ de brins de l'\'etiquette est \'egalement le nombre de composantes de l'entrelacs (trivial) associ\'{e}. 

	D\'emontrons le point $(ii)$. 
	Si on a une ar\^{e}te \'{e}tiquet\'{e}e $C$ et reliant des sommets de~$\arbre$ \'etiquet\'es respectivement $(b_0, n_0)$ et $(b_1, n_1)$, alors, par d\'{e}finition de l'arbre de calcul et des scindements de Conway, on a
	\[\longueur(b_1) - n_1 = \longueur(b_0) - n_0 - 2 = \longueur(b_0) - n_0 - \deg(C).\] 
	De m\^{e}me, si on a une ar\^{e}te \'{e}tiquet\'{e}e $R$, on a
	\[\longueur(b_1) - n_1 = \longueur(b_0) - n_0 - 1  = \longueur(b_0) - n_0 - \deg(R).\]
	Par r\'ecurrence, si la racine de $\arbre$ est \'{e}tiquet\'{e}e $(b,n)$ et si $f$ est une feuille dont l'\'{e}tiquette est une tresse triviale \`a $n_f$ brins, on a
	\[\longueur(b) - n = \deg(\monj_f) - n_f = \deg\left(\monj_f \cdot T^{n_f-1}\right) - 1.\] 
	Somme de mon\^omes tous de m\^eme degr\'e, $J_{(b,n)}(R,C,T)$ est donc homog\`ene de degr\'{e} $\longueur(b) - n + 1$.

	D\'emontrons maintenant le point $(iii)$. 
	Revenons aux arbres de calcul positifs calculant le polyn\^ome HOMFLY. 
	On remarque que, si $f$ est une feuille de $\arbre$, on a l'\'{e}galit\'{e} $\deg(\monj_f)=\deg_x (\monome_f)$. 
	Or, si $(b,n)$ est l'\'{e}tiquette de la racine de $\arbre$ et si $f$ est une feuille, on a 
		\[\longueur(b) - n = \deg(\monj_f) - n_f = \deg_x(\monome_f) - n_f.\]
	Par cons\'{e}quent, il existe un exposant $s(f)$ v\'erifiant 
	$\monome_f(x,y) = \pm y^{s(f)}x^{n - \longueur(b) - n_f}$, d'o\`{u}
		\[\monome_f(x,y) \cdot \left(-\frac{x+x^{-1}}{y}\right)^{n_f-1} 
	   	= \pm y^{s(f) - n_f + 1} x^{n - \longueur(b) + 1}(1 + x^{-2})^{n_f - 1}.\] 
	Posons $k = n - \longueur(b) + 1$. 
	Par d\'{e}finition d'un arbre de calcul positif, on a alors
		\[\HOM_{(b,n)}(x,y) = x^k \sum_{f \in \feuilles}\pm y^{s(f) - n_f - 1}(1+x^{-2})^{n_f-1}.\] 
	Il existe par cons\'{e}quent des entiers $a_{ij}$ uniques v\'erifiant
	\begin{equation*}
		\HOM_{(b,n)}(x,y) = \sum_{i,j}a_{ij}x^{k - j}y^{i + j}((x+x^{-1})/y)^j
		= \sum_{i,j}(-1)^{\frac{i-k}2}a_{ij} (-x^{-2})^{-\frac{i+k}2}(-y/x)^{i+j}(-(x+x^{-1})/y)^j.
	\end{equation*}
	En substituant $C=-x^{-2}, R=-y/x$ et $T=-(x+x^{-1})/y$, on obtient
	\begin{equation}
		\label{E:Jtilde}
		\HOM_{(b,n)}(x,y)
		=\sum_{i,j}(-1)^{\frac{i-k}2}a_{ij}R^{i+j}C^{-\frac{i+k}2}T^j.
	\end{equation} 
	Appelons $\widetilde J(R,C,T)$ le polyn\^ome du second membre de~\eqref{E:Jtilde}. 
	Alors $\widetilde J$ est homog\`ene de degr\'{e} $-k$ si on donne \`{a} $R,C$ et $T$ les poids $1,2$ et $-1$ respectivement. 
	Comme les $a_{ij}$ ne d\'ependent que de~$\HOM_{(b,n)}$, le polyn\^ome $\widetilde J$ est l'unique polyn\^ome homog\`ene se sp\'{e}cialisant en $\HOM_{(b,n)}$ lorsqu'on substitue $-x^{-2}$ \`a~$C$, $-y/x$ \`a $R$ et $-(x+x^{-1})/y$ \`a~$T$. 
	Or $J_{(b,n)}$ a les m\^{e}mes propri\'{e}t\'{e}s, donc les polyn\^omes~$J_{(b,n)}$ et~$\widetilde J$ co\"\i ncident. 
	Par cons\'equent, $J_{(b,n)}$ ne d\'{e}pend que du polyn\^ome $\HOM_{(b,n)}$, et donc {\it a fortiori} que de la cl\^oture de~$(b, n)$.

	Pour ce qui est de $(iv)$, remarquons que la seule branche de l'arbre~$\arbre$ ne comportant aucune \'etiquette~$R$ est la branche la plus \`{a} droite. 
	Soit $(b_0, n_0)$ l'\'etiquette d'un sommet de cette branche, et $(b_1, n_1)$ l'\'etiquette de son fils droit. 
	Alors, par construction de~$\arbre$, il existe un entier~$i$ tel qu'on passe de $b_0$ \`a $b_1 \sigma_i^2$ par des transformations de Markov positives. 
	Ces transformations ne changent pas le nombre de composantes de la cl\^oture, et, d'autre part, en supprimant un facteur $\sigma_i^2$ d'une tresse, on ne change pas non plus le nombre de composantes de sa cl\^oture. 
	Par cons\'equent, les cl\^otures de $(b_0, n_0)$ et $(b_1,n_1)$ ont le m\^eme nombre de composantes. 
	De proche en proche, on en d\'eduit que les n\oe uds associ\'es \`a la racine et \`{a} la feuille la plus \`a droite de $\arbre$ ont le m\^{e}me nombre de composantes, d'o\`u le point~$(iv)$.

	Enfin, pour $(v)$, on revient \`{a} la formule~\eqref{E:defJ}.
	Observons que l'exposant maximal de $T$ dans $J_{(b,n)}$ est de la forme $n_f-1$ pour une certaine feuille $f$ de $\arbre$. 
	Or, le nombre maximal de composantes d'un entrelacs associ\'{e} \`{a} une feuille est born\'{e} par le nombre de brins d'une tresse ayant cet entrelacs pour cl\^oture. 
	Donc, partant d'une tresse marqu\'ee \`{a} $n$ brins, notre construction d'arbre de calcul positif ne fait pas appara\^itre de tresse \`{a} plus de $n$ brins, d'o\`{u} le point~$(v)$.
\end{proof}

\subsection{Application au calcul de l'indice}
\label{S:Application}

Les r\'esultats pr\'ec\'edents permettent de d\'eterminer avec pr\'ecision l'indice 
de tresse de la cl\^oture de certaines tresses positives, en particulier des
tresses de Birman-Williams.

\begin{proposition}
	\label{T:NbComposants}
	Supposons que $(b,n)$ est une tresse positive marqu\'ee et que $\arbre$ est un arbre de calcul positif dont la racine est \'etiquet\'ee par $(b,n)$ et tous les sommets par des tresses positives. 
	Alors, si $\arbre$ a un sommet \'{e}tiquet\'{e} par une tresse dont la cl\^oture est un entrelacs \`a $n$~composantes, l'indice de tresse de la cl\^oture de $(b,n)$ est~$n$.
\end{proposition}

\begin{proof}
	Soit $(b_0, n)$ le sommet de $\arbre$ dont l'entrelacs associ\'{e} a $n$ composantes. 
	D'apr\`es le point $(iv)$ du th\'{e}or\`eme~\ref{T:RCT}, le polyn\^ome \`a coefficients positifs $J_{(b_0,n)}$ a au moins un mon\^ome dont l'exposant de $T$ est \'{e}gal \`{a} $n-1$. 
	Alors, en utilisant l'arbre de calcul $\arbre$, on d\'{e}duit $J_{(b,n)} = P_1J_{(b_0,n)} + P_2$, o\`{u} $P_1$ et $P_2$ sont des polyn\^omes \`{a} coefficients positifs. 
	Par cons\'{e}quent, $J_{(b,n)}$ a au moins un terme dont le degr\'e en~$T$ est au moins \'{e}gal \`{a} $n-1$.
	D'apr\`es le point $(v)$ de la proposition~\ref{T:RCT}, l'indice de tresse de la cl\^oture de $(b,n)$ est au moins~$n$. 
	D'un autre c\^ot\'{e}, $b$ est une tresse \`a $n$ brins, donc l'indice de tresse de sa cl\^oture est au plus~$n$. 
	De ces deux in\'{e}galit\'{e}s, on d\'{e}duit que l'indice de tresse de la cl\^oture de $(b,n)$ est $n$ exactement.
\end{proof}

Nous sommes maintenant pr\^ets \`a d\'emontrer le th\'eor\`eme~\ref{T:F-W}.

\begin{likeproof}[D\'emonstration du th\'eor\`eme~\ref{T:F-W}]
	D'apr\`es la proposition~\ref{T:NbComposants}, il suffit de montrer  que, pour toute tresse positive~$b$ \`a au plus $t$~brins, il existe un arbre de calcul positif pour~$(\Delta_t^2 b, t)$ contenant un sommet \'etiquet\'e~$(\Delta_t^2, t)$. 
	On proc\`ede par r\'ecurrence sur la longueur de~$b$. Pour $b = 1$, le r\'esultat est trivial. Supposons $b = \sigma_i b'$.
	On sait, voir par exemple~\cite{A:Garside}, qu'il existe une tresse positive~$\Delta'$ v\'erifiant $\Delta_t^2 = \Delta' \sigma_i$.
	On a alors $\Delta_t^2 b = \Delta' \sigma_i^2 b'$. 
	Il existe donc un scindement de Conway de la tresse positive marqu\'ee~$(b, t)$ en les tresses positives marqu\'ees $(\Delta' \sigma_i b', t)$, c'est-\`a-dire $(\Delta_t^2 b', t)$, et $(\Delta' b', t)$.
	Comme la longueur de $b'$ est strictement inf\'erieure \`a celle de~$b$, l'hypoth\`ese de r\'ecurrence garantit qu'il existe un arbre de calcul positif~$\arbre'$ pour~$(\Delta_t^2 b', t)$ contenant un sommet \'etiquet\'e~$(\Delta_t^2, t)$. 
	Comme $\Delta' b'$ est une tresse positive, on peut compl\'eter~$\arbre'$ en un arbre de calcul positif pour~$(\Delta_t^2 b, t)$, lequel contient un sommet \'etiquet\'e~$(\Delta_t^2, t)$ puisqu'il inclut~$\arbre'$.
\end{likeproof}

\subsection{N\oe uds de Lorenz et $T$-entrelacs}
\label{S:B-K}

Le th\'eor\`eme~\ref{T:tresseB-W} fournit, pour tout n\oe ud de Lorenz, une tresse ayant un nombre de brins minimal. 
Cependant, il existe diff\'erentes orbites du flot de Lorenz qui sont isotopes, mais fournissent des tresses de Birman-Williams diff\'erentes. 
Ainsi, la repr\'esentation n'est pas injective. 
R\'eduire la redondance de ce codage, ou trouver des mouvements passant d'une tresse de Birman-Williams \`a une autre d\'ecrivant le m\^eme n\oe ud de Lorenz est donc un objectif int\'eressant.

R\'ecemment, J.\,Birman et I.\,Kofman \cite{A:B-K} ont montr\'e que les n\oe uds de Lorenz, et plus g\'en\'eralement les entrelacs de Lorenz g\'en\'eralis\'es (pour lesquels on autorise diff\'erentes composantes \`a suivre des trajectoires parall\`eles sur le patron de Lorenz), admettent une autre description combinatoire.

\begin{definition}\label{D:Tlinks}
	Soit $2\le r_1\le r_2 \le \ldots \le r_k$ et $s_1, \ldots, s_k>0$ des entiers. 
	On d\'efinit le {\it T-entrelacs} $T\left((r_1, s_1), \ldots, (r_k, s_k)\right)$ comme la cl\^oture de la tresse
\[(\sigma_1\sigma_2\ldots\sigma_{r_1-1})^{s_1} (\sigma_1\sigma_2\ldots\sigma_{r_2-1})^{s_2}\dots (\sigma_1\sigma_2\ldots\sigma_{r_k-1})^{s_k}.\]
\end{definition}

La famille des $T$-entrelacs g\'en\'eralise celle des entrelacs alg\'ebriques, puisque le $T$-entrelacs $T\left((r_1, s_1), \ldots, (r_k, s_k)\right)$ est alg\'ebrique si et seulement si $r_1\vert r_2 \vert \ldots \vert r_k$ et $s_i\ge r_i$ pour tout $i$. 
Le r\'esultat suivant affirme qu'on passe des entrelacs alg\'ebriques \`a ceux de Lorenz justement en levant ces conditions.

\begin{theorem}[\cite{A:B-K}]
	La famille des entrelacs de Lorenz g\'en\'eralis\'es et celle des $T$-entrelacs co\"incident. 
\end{theorem}

La d\'emonstration se base sur des manipulations alg\'ebriques dans le groupe des tresses \`a partir de la tresse de Birman-Williams permettant de la transformer en une $T$-tresse. 
Elle donne \'egalement des formules pour les param\`etres $r_i$ et $s_i$ en fonction des param\`etres $n_i$ et $m_i$ d\'efinis en~\eqref{E:ni} et~\eqref{E:mi}. 
Bien que ce th\'eor\`eme permette de r\'eduire encore la redondance du codage des n\oe uds de Lorenz, il ne permet pas de r\'esoudre directement le probl\`eme d'isotopie, puisque, par exemple, la sym\'etrie du patron de Lorenz montre que les entrelacs $T\left((r_1, s_1), (r_2, s_2)\right)$ et $T\left((s_2, r_2-r_1), (s_1+s_2, r_1)\right)$ sont toujours isotopes.


\section{N{\oe}uds de Lorenz, n\oe uds modulaires et th\'eorie des nombres}
\label{S:Arithmetique}

Cette partie est consacr\'ee au lien entre n{\oe}uds de Lorenz et n{\oe}uds modulaires. 
On y explique l'origine arithm\'etique de ces derniers \`a partir des classes de conjugaison dans~$\SLZ$, on d\'ecrit la correspondance de Ghys, et on \'etablit les r\'esultats nouveaux annonc\'es dans l'introduction. 
Dans la section~\ref{S:Conjugaison}, nous d\'eterminons des repr\'esentants des classes de conjugaison dans~$\SLZ$. 
La section~\ref{S:Ideaux} est centr\'ee sur la correspondance classique entre les classes de conjugaison pr\'ec\'edentes et les classes d'id\'eaux dans un ordre d'un corps quadratique. 
Nous introduisons la surface modulaire et les n{\oe}uds modulaires dans les sections \ref{S:SurfaceModulaire} et~\ref{S:NoeudsModulaires}. 
C'est dans la section~\ref{S:ModulairesLorenz} que nous d\'ecrivons la correspondance de Ghys proprement dite. 
Puis, dans la section~\ref{S:GroupeClasses}, nous rappelons la d\'efinition du groupe des classes. 
Enfin, c'est dans la section~\ref{S:Nouveaux} que nous \'etablissons les r\'esultats sur les orbites triviales et inverses.

\begin{figure}[htb]
	\begin{center}
	\begin{picture}(121,36)(0,0)
	\put(15,2){\includegraphics[scale=1]{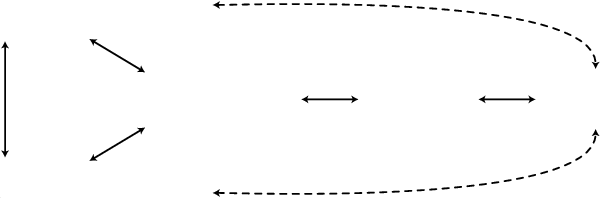}}
	\put(0,35){classes de formes quadratiques} 
	\put(7,31.5){\`a coefficients entiers} 
	\put(8,3.5){classes d'id\'eaux dans} 
	\put(7,0){des corps quadratiques} 
	\put(31,20){classes de conjugaison} 
	\put(39,16){dans $\SLZ$}
	\put(80,20){n\oe uds} 
	\put(77,16){modulaires}
	\put(110,20){n\oe uds} 
	\put(107,16){de Lorenz}
	\put(3,18){({\it Gauss})} 
	\put(95,21){({\it Ghys})}
	\put(110,32){?}
	\put(110,3){?}
	\end{picture}
	\end{center}
	\caption{\small \sf Liens entre les n\oe uds de Lorenz et la th\'eorie des nombres.
	L'existence d'un lien direct entre les n\oe uds de Lorenz et les classes de formes quadratiques d'une part, et les classes d'id\'eaux d'autre part, reste hypoth\'etique.} 
	\label{F:trianglenombres}
\end{figure}


\subsection{Classes de conjugaison dans $\SLZ$ sous l'action de~$\GLZ$.}
\label{S:Conjugaison}

Nous commen\c{c}ons par un r\'esultat pr\'eparatoire sur les \'el\'ements de~$\PSL$, \`a savoir la d\'etermination des classes de conjugaison de matrices hyperboliques. 
On rappelle qu'une matrice est dite \emph{hyperbolique} si la valeur absolue de sa trace est strictement plus grande que~$2$. 

Dans~$\PSLR$, les classes de conjugaison sont faciles \`a d\'eterminer, puisque, \`a trace~$t$ fix\'ee, il y a une classe de conjugaison si $\vert t \vert > 2$ et deux classes si $\vert t \vert \le 2$ (les matrices de rotation $\mathrm{rot}(\theta)$ et $\mathrm{rot}({-\theta})$ ne sont pas conjugu\'ees dans~$\PSLR$).

\begin{proposition}[voir \cite{A:Katok}, chap. 3.5]
	\label{T:presentationSL}
	Dans $\PSL$, on note $a$ la classe contenant les deux matrices~$\pm \mama 0 1 {-1} 0$ et $b$ la classe contenant $\pm \mama 0 1 {-1} 1$.
	Alors le groupe $\PSL$ est engendr\'e par $a$ et $b$, avec les relations $a^2=b^3=1$.
\end{proposition}

Appelons {\it mot r\'{e}duit} un mot en les lettres $a$ et $b$ ne contenant aucun sous-mot $a^2$ ou $b^3$. 
La proposition~\ref{T:presentationSL} implique que tout \'{e}l\'{e}ment de~$\PSL$ est repr\'{e}sent\'{e} par un unique mot r\'{e}duit\footnote{Autrement dit, le groupe $\PSL$ est isomorphe au produit libre $\ZZ/2\ZZ * \ZZ/3\ZZ$.}.

\begin{lemma}
	\label{L:motreduit}
	Tout \'el\'ement de~$\PSL$ qui n'est conjugu\'e ni \`a~$a$, ni \`a~$b$, ni \`a~$b^2$, est conjugu\'e, dans $\PSL$, \`a un \'el\'ement dont le repr\'esentant r\'eduit commence par~$b$ et finit par~$a$.
\end{lemma}

\begin{proof}
Soit $g$ dans $\PSL$. 
Par conjugaison, on peut permuter cycliquement l'ordre des lettres du mot r\'eduit associ\'e. 
Si $g$ a un repr\'{e}sentant r\'{e}duit contenant des lettres $a$ et $b$, alors on peut en permuter les facteurs afin que $b$ soit en t\^{e}te, et $a$ en queue. 
Sinon, les seuls mots r\'{e}duits qui ne comportent pas \`a la fois des $a$ et des $b$ sont $a$, $b$ et~$b^2$.
\end{proof}

Passons au groupe $\SLZ$, et posons $X=\mama 1 0 1 1$ et $Y=\mama 1 1 0 1$. On d\'eduit du lemme~\ref{L:motreduit} une caract\'{e}risation des classes de conjugaison hyperboliques de~$\SLZ$.

\begin{proposition}
	\label{T:conjugaisonSL}
	Toute matrice hyperbolique de $\SLZ$ est conjugu\'{e}e, dans~$\SLZ$, \`{a} un  produit de matrices~$X$ et $Y$ contenant au moins un facteur~$X$ et au moins un facteur~$Y$, et ce produit est unique \`{a} permutation circulaire des facteurs pr\`es.
\end{proposition}

\begin{proof}
Soit $M$ une matrice de~$\SLZ$ de trace $> 2$. 
La classe de $X$ dans $\PSL$ est~$ba$ et celle de $Y$ est $b^2a$. 
Par cons\'equent, d'apr\`es le lemme~\ref{L:motreduit}, il existe un produit~$P$ de matrices $X$ et~$Y$ tel que $M$ est conjugu\'ee \`a~$\pm P$. 
Or les traces de~$M$ et $P$ sont positives, donc $M$ est conjugu\'ee \`a~$P$.
D'autre part, pour tous $m$ et $n$ entiers, on a $tr(X^m)=tr(Y^n)=2$. 
Or on a suppos\'e $tr(M) > 2$. 
Donc $P$ ne peut \^etre ni de la forme~$X^m$, ni de la forme~$Y^n$, et il contient donc au moins un facteur~$X$ et au moins un facteur~$Y$. 

Si $P'$ est obtenu \`a partir de~$P$ par permutation cyclique des lettres, alors $P'$ est conjugu\'e \`a~$P$, et donc \`a $M$. 
Enfin, l'unicit\'e de~$P$ d\'ecoule de l'unicit\'e du repr\'esentant r\'eduit du lemme~\ref{L:motreduit}.
\end{proof}

En rapprochant la proposition~\ref{T:conjugaisonSL} de la proposition~\ref{T:Lyndon} affirmant que les orbites p\'eriodiques du flot de Lorenz sont repr\'esent\'ees par les mots en deux lettres $\xx$ et $\yy$, uniques \`a permutation cyclique des lettres pr\`es, on peut entrevoir un lien entre classes de conjugaison dans~$\SLZ$ et n{\oe}uds de Lorenz. 
C'est ce lien qui sera explicit\'e et approfondi dans la section~\ref{S:ModulairesLorenz}.
Dans un groupe, un \'el\'ement~$g$ est dit \emph{primitif} si l'\'egalit\'e $g=h^k$ implique $k=\pm 1$.
C'est une propri\'et\'e invariante par conjugaison.

\begin{corollary}
\label{C:ConjLyndon}
Les classes de conjugaison de matrices de $\SLZ$ hyperboliques et primitives sont index\'ees par les mots de Lyndon.
\end{corollary}

\begin{exemple}\label{Ex:tracedix1}
Des tests rapides permettent de v\'erifier qu'il n'y a que trois mots de Lyndon associ\'es \`a des classes de conjugaison de matrices de trace~10 et de d\'eterminant~1, \`a savoir les mots $X^8Y, X^4Y^2$ et $X^2YXY$. D'apr\`es le r\'esultat pr\'ecedent, toute matrice primitive de $\SLZ$ de trace~10 est donc conjugu\'ee \`a l'un de ces produits.

Consid\'erons par exemple la matrice $M=\mama 5 {-8} {-3} 5$. 
L'algorithme d'Euclide appliqu\'e \`a la premi\`ere colonne de $M$ permet de d\'ecomposer celle-ci en le produit $Y^{-2}X^{3}Y^{-2}$. 
Dans $\PSL$ la classe de ce produit est $(b^2a)^{-2}(ba)^3(b^2a)^{-2}=(ab)^2(ba)^3b(ab)^2=abab^2abab^2ab$. 
Une permutation cyclique des \'el\'ements de ce produit montre qu'il est conjugu\'e \`a $(ba)^2(ba^2)(ba)(ba^2)$, qui est la classe du produit $X^2YXY$. 
On en d\'eduit que la matrice $M$ est conjugu\'ee \`a~$X^2YXY$. 
Le m\^eme proc\'ed\'e permet de d\'eterminer le mot de Lyndon repr\'esentant la classe de conjugaison de n'importe quelle matrice primitive hyperbolique de~$\SLZ$.
\end{exemple}

Dans la suite, nous allons consid\'erer des classes un peu plus grosses que les classes de conjugaison par des matrices de~$\SLZ$. 
On note $\SLpm$ l'ensemble des matrices \`a coefficients dans $\relatifs$ de d\'{e}terminant $\pm 1$. 
On s'int\'{e}resse maintenant \`{a} la conjugaison par des matrices de $\SLpm$. 

\begin{proposition}\label{T:conjugaisonGL}
Toute matrice hyperbolique et de trace positive de $\SLZ$ est conjugu\'{e}e, dans~$\SLpm$, \`{a} un produit de $X$ et de $Y$ contenant au moins un facteur~$X$ et au moins un facteur~$Y$, et ce produit est unique \`{a} permutation circulaire des facteurs pr\`es et \`{a} interversion des caract\`eres $X$ et $Y$ pr\`es.
\end{proposition}

\begin{proof}
Comme $\SLZ$ est un sous-groupe (d'indice~$2$) de $\SLpm$, si deux matrices sont conjugu\'{e}es via des matrices de $\SLZ$, elles le sont via des matrices de $\SLpm$. 
D'autre part, toute matrice de d\'{e}terminant $-1$ est le produit d'une matrice de $\SLZ$ par la matrice $T=\mama 0 1 1 0$. 
Donc il suffit de conna\^\i tre le r\'{e}sultat de l'action de la matrice $T$ par conjugaison sur $X$ et $Y$ pour adapter la proposition~\ref{T:conjugaisonSL}. 
Or on a $TXT^{-1}=Y$ et $TYT^{-1}=X$, donc conjuguer par $T$ revient \`{a} intervertir les lettres $X$ et $Y$.
\end{proof} 


\subsection{Classes d'id\'{e}aux et classes de conjugaison}
\label{S:Ideaux}

Nous allons maintenant \'{e}tablir une correspondance classique en th\'{e}orie des nombres entre classes d'id\'{e}aux d'un ordre dans un corps quadratique et  classes de conjugaison de d\'{e}terminant non nul dans $\MMZ$ sous l'action de $\SLpm$. 
Les r\'{e}sultats s'\'{e}tendent \`{a} des corps de degr\'{e} sup\'{e}rieur et des matrices de taille sup\'{e}rieure. 
Cependant, ils s'\'{e}noncent plus facilement dans le cas de la taille 2 qui nous int\'{e}resse ici. 

Un {\it corps quadratique} est un corps de la forme $\rationnels[\alpha]$ avec $\alpha$ nombre alg\'{e}brique de degr\'{e}~$2$. 
L'{\it anneau des entiers} $\grandordre$ de~$\rationnels[\alpha]$ est l'ensemble des \'{e}l\'{e}ments de $\rationnels[\alpha]$ admettant un polyn\^ome minimal unitaire \`{a} coefficients entiers. 
Un {\it ordre} de $\grandordre$ est un sous-anneau de~$\grandordre$ contenant l'\'{e}l\'{e}ment~$1$ et non isomorphe \`{a}~$\relatifs$.
Un r\'{e}sultat classique affirme que $\grandordre$ est un r\'{e}seau de $\rationnels [\alpha]$, c'est-\`{a}-dire un $\ZZ$-module libre de rang~$2$ dans $\rationnels[\alpha]$. 
Un ordre est alors un sous-r\'{e}seau de~$\grandordre$ contenant~$1$.

\begin{example}
Dans le corps $\rationnels[\!\sqrt5]$, l'anneau des entiers est $\relatifs[\frac{1+\!\sqrt5}{2}]$. 
Deux exemples d'ordres sont $\relatifs[\!\sqrt5]$, qui est d'indice 2, et $\relatifs[\frac{1+3\!\sqrt5}{2}]$, qui est d'indice~3.
\end{example}

Un {\it id\'{e}al} $\Ideal$ d'un ordre $\ordre$ est un sous-groupe additif de $\ordre$ stable par multiplication par les \'{e}l\'{e}ments de~$\ordre$. 
Si $\lbid\omega_1, \omega_2\rbid$ est une $\relatifs$-base de~$\ordre$, alors tout id\'{e}al~$\Ideal$ de~$\ordre$ admet une
$\relatifs$-base $\lbid\alpha_1, \alpha_2\rbid$ o\`{u} $\alpha_1$ et $\alpha_2$ sont deux \'{e}l\'{e}ments du r\'{e}seau engendr\'{e} par $\omega_1$ et $\omega_2$. 

Un id\'eal~$\Ideal$ d'un ordre~$\ordre$ est dit {\it principal} s'il existe un \'el\'ement $a$ dans $\Ideal$ tel que $\Ideal$ est exactement l'ensemble des multiples de $a$ par des \'el\'ements de~$\ordre$. 
On note alors $\lid a\rid$ cet id\'eal. 
Soit $\Ideal, \Ideal'$ deux id\'{e}aux de $\ordre$. 
On dit que $\Ideal$ et $\Ideal'$ sont dans la m\^{e}me {\it classe} s'il existe $a, a'$ non nuls dans $\ordre$ v\'erifiant $\lid a\rid\Ideal = \lid a'\rid\Ideal'$. 
Il s'agit d'une relation d'\'{e}quivalence entre id\'{e}aux. 
Par exemple, la classe de l'id\'{e}al~$\lid 1 \rid$ est l'ensemble des id\'{e}aux principaux de $\ordre$. 
Le r\'esultat principal de cette section est le suivant.

\begin{proposition}[voir \cite{A:Cohn}]
	\label{T:nombreclasses}
	Soit $t$ et $d$ deux entiers vérifiant $t^2 > 4d > 0$. 
	Soit $P(x)$ le polyn\^ome~$x^2 - tx + d$ et soit $\alpha$ une racine de~$P$. 
	Alors il existe une bijection entre les classes de conjugaison sous l'action de $\SLpm$ de matrices~$A$ de~$\MMZ$ v\'{e}rifiant $P(A)=0$ d'une part, et les classes d'id\'{e}aux de l'ordre~$\relatifs[\alpha]$ d'autre part.
\end{proposition}

\begin{proof}
Soit $A$ dans $\MMZ$ v\'erifiant $P(A)=0$. 
Alors $A$ admet $\alpha$ et $d/\alpha$ pour valeurs propres et $P$ pour polyn\^ome caract\'{e}ristique. 
Soit $v$ un vecteur propre de~$A$ associ\'e \`{a} la valeur propre~$\alpha$. 
Le vecteur~$v$ est unique \`{a} multiplication par un scalaire pr\`es et peut \^{e}tre choisi dans~$\relatifs[\alpha]^2$. 
On le note alors $\left(\begin{matrix} {\alpha_1} \\ {\alpha_2} \end{matrix} \right)$.
Par d\'{e}finition, on a 
   $A\left(\begin{matrix} {\alpha_1} \\ {\alpha_2} \end{matrix} \right) 
   =\alpha \left( \begin{matrix} {\alpha_1} \\ {\alpha_2} \end{matrix} \right)$. 
Par cons\'{e}quent, pour tout polyn\^ome $Q$ \`{a} coefficients entiers, on a 
   $Q(A)\left( \begin{matrix} {\alpha_1} \\ {\alpha_2} \end{matrix} \right) 
   = Q(\alpha) \left( \begin{matrix} {\alpha_1} \\ {\alpha_2} \end{matrix} \right)$. 
Donc $\lbid\alpha_1, \alpha_2\rbid$ est une $\relatifs$-base d'un certain id\'{e}al~$\Ideal_v$ de~$\relatifs[\alpha]$.

Si $v'$ est un autre vecteur propre associ\'e \`a la valeur propre~$\alpha$ dont les coordonn\'ees sont dans~$\relatifs[\alpha]$, alors on a $v'=\lambda v$ pour un certain $\lambda$ dans~$\relatifs[\alpha]$, et donc l'id\'eal associ\'e est $\lid\lambda\rid\Ideal_v$. 
La classe de $\Ideal_v$ ne d\'epend donc que de~$A$. 
On note cette classe~$\Cl(A)$.

Si $A'$ est conjugu\'ee \`a $A$ par une matrice~$M$ de~$\SLpm$, alors  le vecteur $M^{-1}v$ est un vecteur propre associ\'e \`a la valeur propre~$\alpha$ de $A'$. 
Comme $M$ est de d\'eterminant~$\pm 1$, l'id\'eal $\Ideal_{M^{-1}v}$ est \'egal \`a~$\Ideal_v$, mais la base associ\'ee est diff\'erente. 
La classe~$\Cl(A)$ ne d\'epend donc que de la classe de conjugaison de $A$ sous l'action de~$\SLpm$.
Donc $\Cl$ induit une application bien d\'efinie, \'egalement not\'ee~$\Cl$ dans la suite, de l'ensemble des classes de conjugaison de matrices~$A$ v\'erifiant $P(A) = 0$ vers les classes d'id\'eaux de l'ordre~$\relatifs[\alpha]$.

Dans l'autre direction, soit $\Ideal$ un id\'{e}al de $\relatifs[\alpha]$ et $\lbid\beta_1, \beta_2\rbid$ une base de~$\Ideal$.  
Comme $\Ideal$ est stable par multiplication par~$\alpha$, il existe une matrice $B$ dans~$\MMZ$ v\'erifiant $B \left(\begin{matrix} {\beta_1} \\ {\beta_2} \end{matrix} \right) 
= \alpha \left( \begin{matrix} {\beta_1} \\ {\beta_2} \end{matrix} \right)$. 
Par cons\'equent $B$ admet $\alpha$ pour valeur propre et $\left( \begin{matrix} {\beta_1} \\ {\beta_2} \end{matrix} \right)$ pour vecteur propre associ\'{e}. 
Comme la trace de $B$ est enti\`ere, l'autre valeur propre $\alpha'$ de~$B$ est \'{e}gale \`{a} $n-\alpha$, pour un certain $n$ dans $\relatifs$. 
On a alors  
$$\det B=\alpha(n-\alpha)=n\alpha-\alpha^2=n\alpha-(t\alpha-d).$$
Or $\det B$ est entier, ce qui entra\^ine $n-t=0$, et donc $\alpha'$ est le conjugu\'{e} de $\alpha$. 
Par cons\'{e}quent, la matrice $B$ v\'{e}rifie $P(B)=0$. 

Si on choisit une autre base $\lbid\beta_1', \beta_2'\rbid$ de $\Ideal$, alors il existe $T$ dans $\SLpm$  v\'erifiant $\left( \begin{matrix} {\beta'_1} \\ {\beta'_2} \end{matrix} \right) = T\left( \begin{matrix} {\beta_1} \\ {\beta_2} \end{matrix} \right)$. 
La matrice associ\'ee est alors $T\!\,BT^{-1}$. 
La classe de conjugaison de $B$ ne d\'epend donc que de~$\Ideal$. 
On note cette classe de conjugaison~$\Conj(\Ideal).$

Si on choisit un autre id\'eal $\Ideal'$ dans la classe de $\Ideal$, on v\'erifie \'egalement que la matrice associ\'ee \`a toute base de~$\Ideal'$ est dans~$\Conj(\Ideal)$.

On a ainsi construit deux applications $\Cl$ et $\Conj$ entre classes de conjugaison dans~$\MMZ$ sous l'action de~$\SLpm$ et classes d'id\'eaux dans des corps quadratiques. 
La d\'efinition de~$\Conj$ montre imm\'ediatement que $\Cl$ et $\Conj$ sont inverses l'une de l'autre.
\end{proof}

La question est maintenant de calculer effectivement la bijection~$\Phi$ de la proposition~\ref{T:nombreclasses}, \`{a} savoir, \`{a} partir des matrices, de calculer des \'{e}l\'{e}ments des classes d'id\'{e}aux correspondants. Une notion adapt\'ee \`a cette fin est la suivante.

\begin{definition}[voir~\cite{A:Cohn}]
\label{D:matideale}
Soit~$\ordre$ un ordre dans $\relatifs[\alpha]$ et $\Ideal$ un id\'eal de $\ordre$.
Une matrice $M$ dans $\MMZ$ est dite {\it id\'{e}ale} pour $\Ideal$ s'il existe une base $\{\omega_1, \omega_2\}$ de $\ordre$ et une base~$\{\alpha_1, \alpha_2\}$ de $\Ideal$ telles que
$M\left( \begin{matrix} {\omega_1} \\ {\omega_2} \end{matrix} \right) 
= \left( \begin{matrix} {\alpha_1} \\ {\alpha_2} \end{matrix} \right)$. 
\end{definition}

\begin{proposition}\label{T:bijectionclasses}
Soit $P(x)= x^2 - tx + d$ un polyn\^ome unitaire irr\'{e}ductible et $\alpha$ une racine de~$P$. 
Alors toute matrice $A$ dans $\MMZ$ v\'erifiant $P(A)=0$ est de la forme $X_\Ideal C_PX_\Ideal^{-1}$, o\`{u} $X_\Ideal$ est une matrice id\'{e}ale pour un id\'{e}al~$\Ideal$ appartenant \`a~$\Cl(A)$ et o\`u $C_P$ est la matrice compagnon~$\mama0 {1} {-d} t$ du polyn\^ome $P$.
R\'{e}ciproquement,  toute matrice $X$ dans~$\MMZ$ v\'erifiant $X^{-1}AX=C_P$ est la matrice id\'{e}ale d'un id\'{e}al appartenant \`a~$\Cl(A)$.
\end{proposition}

Noter que la relation $A=X_\Ideal C_PX_\Ideal^{-1}$ n'implique pas que la matrice~$A$ est conjugu\'ee \`a $C_P$ dans~$\SLpm$, car en g\'en\'eral, la matrice $X_\Ideal$ n'est pas de d\'eterminant~$\pm 1$.

\begin{proof}
Par d\'efinition de~$\alpha$, on a 
$\alpha\left( \begin{matrix} {1} \\ {\alpha} \end{matrix} \right)
=\left( \begin{matrix} {\alpha} \\ {t\alpha-d} \end{matrix} \right),$ 
donc, par la d\'emonstration de la proposition~\ref{T:nombreclasses}, la matrice $C_P$ est la matrice de $\Conj\left(\lid1\rid\right)$ associ\'ee \`a la base $\lbid1, \alpha\rbid$ de l'id\'eal~$\lid 1 \rid$.

Soit $A$ une matrice de~$\MMZ$ annulant~$P$. 
On lui a associ\'e dans la d\'{e}monstration de la proposition~\ref{T:nombreclasses} une base $\{\alpha_1, \alpha_2\}$ d'un id\'{e}al $\Ideal$ appartenant \`a~$\Cl(A)$. 
Par d\'efinition, la matrice id\'{e}ale $X_\Ideal$ correspondante, relativement aux bases
$\{1,\alpha\}$ de~$\relatifs[\alpha]$ et $\{\alpha_1, \alpha_2\}$ de $\Ideal$, v\'erifie
   $X_\Ideal\left( \begin{matrix} {1} \\ {\alpha} \end{matrix} \right) 
   = \left( \begin{matrix} {\alpha_1} \\ {\alpha_2} \end{matrix} \right).$ 
On a donc
   $X_\Ideal^{-1}\left( \begin{matrix} {\alpha_1} \\ {\alpha_2} \end{matrix} \right) 
   = \left(\begin{matrix} {1} \\ {\alpha} \end{matrix} \right),$
d'o\`{u}
   $C_PX_\Ideal^{-1}\left( \begin{matrix} {\alpha_1} \\ {\alpha_2} \end{matrix} \right) 
   = C_P\left( \begin{matrix} {1} \\ {\alpha} \end{matrix} \right)
   = \alpha\left( \begin{matrix} {1} \\ {\alpha} \end{matrix} \right),$
et donc
$$X_\Ideal C_PX_\Ideal^{-1}\left( \begin{matrix} {\alpha_1} \\ {\alpha_2}
\end{matrix}
\right) 
   = \alpha X_\Ideal\left( \begin{matrix} {1} \\ {\alpha} \end{matrix} \right)
   = \alpha\left( \begin{matrix} {\alpha_1} \\ {\alpha_2} \end{matrix}
\right).$$ 
Par cons\'{e}quent, on a $A=X_\Ideal C_PX_\Ideal^{-1}$, comme annonc\'e.

Pour la r\'{e}ciproque, soit $\{\alpha_1, \alpha_2\}$ une base de l'id\'{e}al correspondant \`{a} $A$. 
On a alors
    $A\left( \begin{matrix} {\alpha_1} \\ {\alpha_2} \end{matrix} \right) 
    = \alpha\left( \begin{matrix} {\alpha_1} \\ {\alpha_2} \end{matrix} \right).$
Soit $X$ v\'erifiant $X^{-1}AX=C$. 
On a donc
    $XCX^{-1}\left( \begin{matrix} {\alpha_1} \\ {\alpha_2} \end{matrix} \right) 
    = \alpha\left( \begin{matrix} {\alpha_1} \\ {\alpha_2} \end{matrix} \right),$
d'o\`{u}
    $CX^{-1}\left( \begin{matrix} {\alpha_1} \\ {\alpha_2} \end{matrix} \right) 
    = \alpha X^{-1}\left( \begin{matrix} {\alpha_1} \\ {\alpha_2} \end{matrix} \right).$
Comme $\alpha$ est une valeur propre simple de $C$, le vecteur
$X^{-1}\left( \begin{matrix} {\alpha_1} \\ {\alpha_2} \end{matrix} \right)$ est propre, et donc il existe $p, q$ dans $\rationnels(\alpha)$ v\'erifiant
    $X^{-1}\left( \begin{matrix} {\alpha_1} \\ {\alpha_2} \end{matrix} \right) 
    = \displaystyle\frac{p}{q}\left( \begin{matrix} {1} \\ {\alpha} \end{matrix} \right)$, 
soit
    $X^{-1}\left( \begin{matrix} {q \alpha_1} \\ {q \alpha_2} \end{matrix} \right) 
    = \left( \begin{matrix} {p} \\ {p \alpha} \end{matrix} \right)$,
d'o\`u
    $\left( \begin{matrix} {q \alpha_1} \\ {q \alpha_2} \end{matrix} \right) 
    = X\left( \begin{matrix} {p} \\ {p \alpha} \end{matrix} \right)$. 
Par cons\'equent, $X$ est la matrice id\'{e}ale de l'id\'eal $\lid q\alpha_1, q\alpha_2 \rid$, lequel appartient \`a la classe~$\Cl(A)$.
\end{proof}

La d\'{e}monstration de la proposition~\ref{T:bijectionclasses} montre en particulier que la classe des id\'{e}aux principaux de~$\relatifs[\alpha]$ correspond \`{a} la classe de matrices contenant la matrice compagnon~$C_P$ du polyn\^ome caract\'{e}ristique~$P$ de~$\alpha$. 
Dans le cas $P(x)=x^2 - tx + 1$, la classe des id\'{e}aux principaux correspond \'{e}galement \`{a} la classe de matrices engendr\'{e}e par la matrice~$X^{t-2}Y$, puisque $X^{t-2}Y$ est conjugu\'{e}e \`{a}~$C_P$ par la matrice~$\mama 1 0 {-1} 1$, qui est dans~$\SLpm$.

\begin{example}\label{Ex:tracedix2}
Revenons aux matrices de trace 10 et de d\'{e}terminant~1. 
Le polyn\^ome minimal correspondant est $P(x)=x^2-10x+1$, dont les racines sont $5+2\!\sqrt{6}$ et $5-2\!\sqrt{6}$. 
Le corps quadratique associ\'{e} est $\rationnels[\!\sqrt6]$, dont l'anneau des entiers est $\relatifs[\!\sqrt6]$; en revanche, nous nous int\'{e}ressons ici \`a l'ordre~$\relatifs[2 \sqrt6]$ et \`a ses classes d'id\'eaux. 
Nous avons vu dans l'exemple~\ref{Ex:tracedix1} que toute matrice de trace 10 et de d\'eterminant 1 est conjugu\'ee \`a l'un des trois produits $X^8Y, X^4Y^2$ ou $X^2YXY$, c'est-\`a-dire \`a l'une des matrices $\mama 1 1 8 9, \mama 1 2 4 9$ ou~$\mama 2 3 5 8$. 
Par la propri\'et\'e~\ref{T:nombreclasses}, il y a donc trois classes d'id\'eaux dans l'ordre~$\relatifs[2\sqrt 6]$.

On vient de voir que la classe de la matrice $X^8Y=\mama 1 1 8 9$ correspond \`{a} la classe des id\'{e}aux principaux, donc $\mama 1 2 4 9$ et $\mama 2 3 5 8$ correspondent aux deux classes d'id\'{e}aux non principaux. 
Comme les matrices conjuguant ces deux \'{e}l\'{e}ments \`{a} $C_P$ sont $\mama 2 0 {-1} 1$ et $\mama 3 0 {-1} 1$, et que la base de $\relatifs[2\sqrt 6]$ associ\'{e}e \`{a} $C_P$ est $\{1, \alpha\} = \{1, 5+2\sqrt6\}$, on en d\'{e}duit au passage que les id\'{e}aux $\lid 2, 4+2\sqrt 6 \rid$ et $\lid 3, 4+2\sqrt 6 \rid$ sont repr\'{e}sentants des deux classes d'id\'{e}aux non principaux de~$\relatifs[2\sqrt 6]$.
\end{example}


\subsection{La surface modulaire}
\label{S:SurfaceModulaire}

Dans~\cite{A:Ghys3} et~\cite{A:Ghys4}, \'E. Ghys d\'emontre, et illustre visuellement, que les n{\oe}uds de Lorenz apparaissent comme orbites p\'{e}riodiques du flot modulaire d\'{e}fini sur la vari\'{e}t\'{e}~$\PSLR/\PSL$. 
Afin d'expliquer ce r\'esultat, nous allons d'abord donner plusieurs descriptions de la vari\'et\'e modulaire, et d\'ecrire le flot en termes de r\'eseaux du plan.

Le groupe~$\PSLR$ agit sur le demi-plan de Poincar\'e~$\Hy$ par homographie:
\[ \mama a b c d.~z = \frac{az+b}{cz+d},\]
l'action \'etant transitive et fid\`ele. 
Un domaine fondamental~$\domfond$ pour cette action est donn\'e par 
\[\left\{z\in\Compl \mid \Im z >0,  -1/2 \le \Re z \le 1/2,  \vert z \vert \ge 1 \right\}.\]

La restriction de l'action \`a~$\PSL$ est discr\`ete, elle n'est pas libre car certains points, en fait les orbites des points $i$ et $j$, ont des stabilisateurs non triviaux. Le quotient~$\Hy/\PSL$ est alors un orbifold, c'est-\`a-dire une surface avec deux points singuliers (correspondant aux orbites des points $i$ et~$j$) au voisinage desquels l'angle total n'est pas~$2\pi$ mais respectivement $2\pi/2$ et $2\pi/3$.
Cet orbifold, qui n'est donc pas une surface \`a proprement parler, est traditionnellement appel\'e {\it surface modulaire} et not\'e $\surfmod$ dans la suite. 
Il s'identifie naturellement \`a $\domfond$ modulo les recollements $-\frac 1 2+i y\sim \frac 1 2 + i y$ pour $y \ge \sqrt 3 / 2$ et $e^{i(\frac \pi 2 - \alpha)} \sim e^{i(\frac \pi 2 + \alpha)}$ pour $0 < \alpha \le \pi/12$.
Il h\'erite de la m\'etrique hyperbolique de~$\Hy$. 

Comme le jacobien d'une homographie vaut toujours $1$, l'action s'\'etend au fibr\'e unitaire tangent~$\un\Hy \simeq \Hy\times\Sph^1$. 
Elle devient alors transitive et libre, par cons\'equent, les vari\'et\'es $\un\Hy$ et~$\PSLR$ sont isomorphes.
Comme l'action de $\PSL$ sur $\un\Hy$ est discr\`ete, les vari\'et\'es quotients~$\PSLR/\PSL$ et~$\un\Hy/\PSL=\USmod$ sont isomorphes.
Le fibr\'{e} unitaire tangent~$\USmod$ \`a la surface modulaire est alors une vari\'et\'e de dimension 3, naturellement munie d'une m\'etrique riemannienne~$g$. Il s'identifie de mani\`ere canonique \`a $\domfond \times \Sph^1$ modulo les recollements  sur les bords $(-\frac 1 2 + iy, \theta)\sim (\frac 1 2 + iy, \theta)$ d'une part, et $(e^{i(\frac \pi 2 - \alpha)}, \theta-\alpha) \sim (e^{i(\frac \pi 2 + \alpha)}, \pi + \theta+\alpha)$ d'autre part, qui se prolongent aux coins en $(i, \theta) \sim (i, \theta+\pi)$ et $(j, \theta) \sim (j, \theta+\frac{2\pi} 3)$.

\begin{figure}[htb]
	\begin{center}
	\begin{picture}(120,80)(0,0)
	\put(0,0){\includegraphics[width=.8\textwidth]{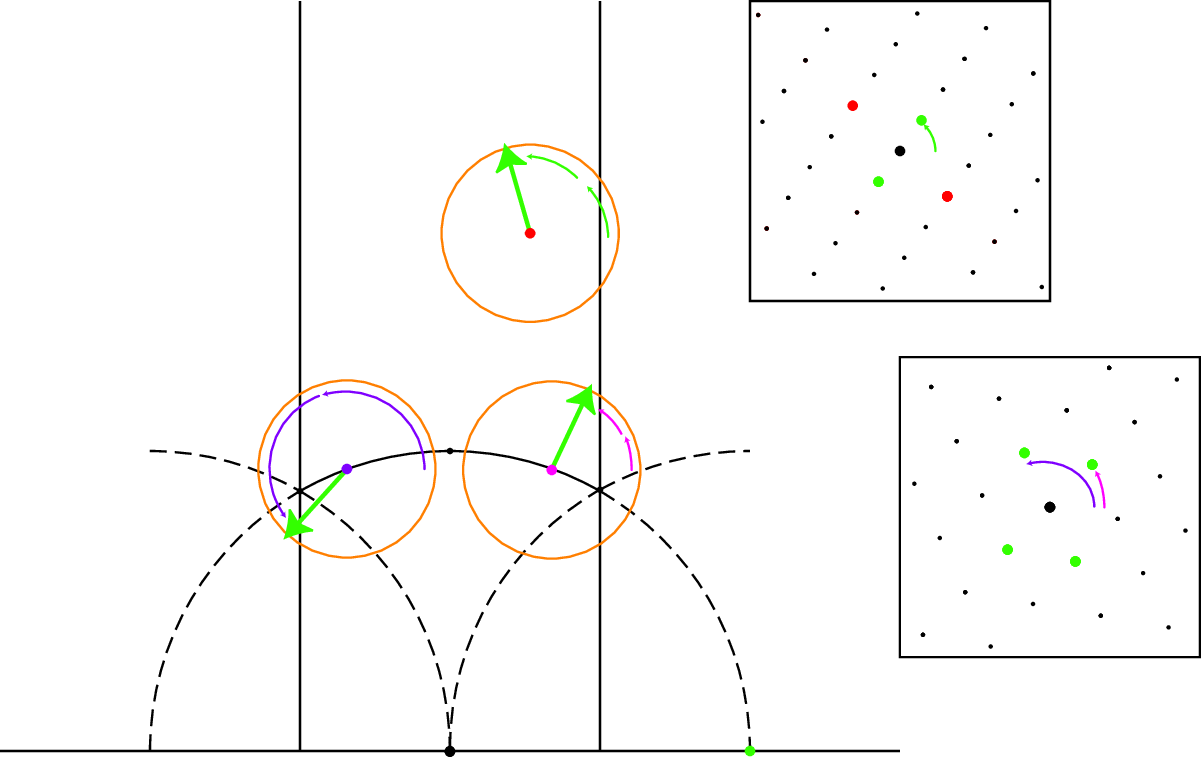}}
	\put(44,-3){0}
	\put(74,-3){1}
	\put(44,32){$i$}
	\put(51.5,49.5){$z$}
	\put(49,62){$2\theta$}
	\put(93.5,61.5){$\theta$}
	\put(91,64.5){$v$}
	\put(83.5,66){$v'$}
	\put(52.5,26){$z_1$}
	\put(60,37){$2\theta_1$}
	\put(109.5,29){$v_1$}
	\put(96,31){$v_2=v'_1$}
	\put(97,17){$v'_2$}
	\put(34,26){$z_2$}
	\put(24.5,19){$2\theta_2$}
	\put(36,64){$\domfond$}
	\end{picture}
	\end{center}
	\caption{\small \sf Comment associer un point de $\USmod$ \`a un r\'eseau de covolume~$1$. Les deux vecteurs tangents repr\'esent\'es en $z_1$ et $z_2$ correspondent au m\^eme r\'eseau, repr\'esent\'e \`a droite, et sont donc identifi\'es.} 
	\label{F:ReseauModulaire}
\end{figure}

Un r\'eseau du plan est dit \emph{de covolume~$1$} s'il admet une maille d'aire~$1$. 
Un tel r\'eseau est d\'etermin\'e par les coordonn\'ees des vecteurs d'une base, donc par une matrice~$A$ de~$\SLR$. 
De plus, changer la base (ordonn\'ee) du r\'eseau revient \`a multiplier \`a gauche~$A$ par un \'el\'ement de~$\SLZ$. 
On obtient ainsi un isomorphisme~$\resslr$ entre l'espace des r\'eseaux de covolume~1 et la vari\'et\'e~$\SLR/\SLZ \simeq \PSLR/\PSL$. 
Ainsi, les trois vari\'et\'es $\PSLR/\PSL, \USmod$ et l'espace des r\'eseaux de covolume~1 sont canoniquement isomorphes. 
D\'etaillons le lien entre les deux derni\`eres (voir \'egalement la figure~\ref{F:ReseauModulaire}).

Soit $\Lambda$ un r\'eseau de covolume 1 dans le plan, identifi\'e \`a~$\Compl$. 
Nous allons lui associer un point~$\resmod(\Lambda)$ de~$\USmod$, c'est-\`a-dire un point de $\surfmod$ et un vecteur unitaire tangent en ce point. 
L'id\'ee est de passer par  un r\'eseau semblable $\rho e^{i\theta}\Lambda$ ayant une base de la forme $\{1, z\}$, o\`u $z$ est un point de $\domfond$.

Supposons d'abord que $\Lambda$ n'a que deux plus courts vecteurs que nous nommons $v$ et $-v$.
Supposons ensuite que $\Lambda$ n'a que deux plus courts vecteurs non colin\'eaires \`a $v$, nous les nommons $v'$ et $-v'$ de sorte que la base $(v, v')$ soit directe. 
Il existe une unique simillitude directe envoyant $v$ en 1.
On d\'esigne par $z$ l'image de $v'$ et par $\theta$ l'angle de la simillitude.
Remarquons que, si $\theta$ d\'epend du choix de $v$ que nous avons fait, ni $z$, ni $2\theta$ n'en d\'ependent. 
Il se trouve que $z$ est \`a l'int\'erieur du domaine~$\domfond$, et on d\'efinit alors sans ambigu\"it\'e $\resmod(\Lambda)$ par la formule~$\resmod(\Lambda)=(z, 2\theta)$.

Dans la construction pr\'ec\'edente, s'il y a deux choix possibles pour $v'$, alors ces options fournissent deux points $z$ de $\domfond$ de la forme $-\frac 1 2 + iy$ et $\frac 1 2 + iy$, qui correspondent au m\^eme point de $\surfmod$. 
Comme la simillitude ne d\'epend pas du choix de $v'$, ces deux choix possibles de $z$ dans~$\domfond$ m\`enent au m\^eme angle $\theta$, et donc au m\^eme point de~$\USmod$. On peut alors d\'efinir $\resmod\Lambda)$ par la m\^eme formule sans ambigu\" it\'e.

S'il y a plus deux choix possibles pour~$v$, alors soit $\Lambda$ est \`a maille carr\'ee, soit \`a maille triangulaire \'equilat\'erale, soit \`a maille triangulaire isoc\`ele non \'equilat\'erale. 
Dans le dernier cas (illutr\'e sur la figure~\ref{F:ReseauModulaire}), il y a quatre choix possibles pour $v$, menant \`a deux points $z$ de $\domfond$ de la forme $e^{i(\frac \pi 2 - \alpha)}$ et $e^{i(\frac \pi 2 + \alpha)}$, qui correspondent encore au m\^eme point de $\surfmod$. 
Les angles des simillitudes correspondantes diff\`erent de $\frac \pi 2 + \alpha$, dont le double est exactement la diff\'erence d\^ue au recollement de $\domfond$. Par cons\'equent, on peut encore d\'efinir $\resmod(\Lambda)$ sans ambigu\"it\'e.

De m\^eme , dans le cas d'une maille carr\'ee, le point $z$ associ\'e est le point $i$, et les diff\'erents choix du plus court vecteur~$v$ sont compens\'es par l'identification $(i, \theta) \sim (i, \theta+\pi)$. Et dans le cas d'une maille triangulaire équilatérale, les diff\'erents choix possibles, \`a la fois pour $v$ et pour $v'$ sont compens\'es par les identifications $(j, \theta)\sim(j, \theta+\frac{2\pi} 3)\sim(j+1, \theta)$. On a alors

\begin{proposition}
\label{usmodegalreseau}
L'application $\resmod$ \'etablit un isomorphime entre l'espace des r\'eseaux de covolume~$1$ et le fibr\'e unitaire tangent~$\USmod$.
\end{proposition}


\subsection{Flot modulaire, n{\oe}uds modulaires et classes de conjugaison}
\label{S:NoeudsModulaires}

Soit $\Lambda$ un r\'{e}seau de~$\reels^2$ de covolume~1 et $(b_1, b_2)$ une base de~$\Lambda$. 
Pour tout $t$, les vecteurs $b_{1}^t=\mama{e^t}{0}{0}{e^{-t}}b_1$ et $b_{2}^t=\mama{e^t}{0}{0}{e^{-t}}b_2$ forment encore une base d'un r\'{e}seau de covolume~1, lequel ne d\'{e}pend que de $\Lambda$. 
On d\'{e}finit $\phi^t(\Lambda)$ comme le r\'{e}seau engendr\'{e} par les vecteurs $b_{1}^t$ et $b_{2}^t$. 
Le flot~$\phi$ est appel\'e {\it flot modulaire} sur l'espace des r\'eseaux de covolume~1. 
Par la proposition~\ref{usmodegalreseau}, il induit un flot sur la vari\'et\'e~$\USmod$.
Il se trouve que ce dernier est le flot g\'eod\'esique associ\'e \`a la m\'etrique~$g$ sur $\USmod$.

Si les vecteurs $b_1^{t_0}$ et $b_2^{t_0}$ appartiennent au r\'eseau~$\Lambda$ et sont minimaux, alors on a $\phi^{t_0}(\Lambda)=\Lambda$ et l'application $t\mapsto \phi^t(\Lambda)$ restreinte \`a l'intervalle $[0, t_0]$ d\'efinit un lacet ferm\'{e} dans~$\USmod$, et donc un n{\oe}ud dans ce m\^eme espace. 
Le r\'esultat suivant d\'ecrit les circonstances o\`u un tel \'ev\'enement se produit et annonce la correspondance de~Ghys. 

\begin{proposition}
\label{T:noeudsmodulaires}
Les orbites p\'eriodiques du flot modulaire sont index\'ees par les classes de conjugaison de $\PSL$.
\end{proposition}

\begin{likeproof}[Sch\'ema de la d\'emonstration]
Comme le flot modulaire co\"\i ncide avec le flot g\'eod\'esique sur le quotient de~$\Hy$ par $\PSL$, une orbite modulaire p\'eriodique provient d'une g\'eod\'esique~$\geod$ de~$\Hy$ qui est invariante sous l'action d'un certain \'el\'ement hyperbolique~$A$ de $\PSL$. 
Quand on passe au quotient~$\Hy/\PSL$, cette g\'eod\'esique~$\geod$ devient invariante sous l'action de toute la classe de conjugaison de~$A$ dans $\SLZ$. 
R\'eciproquement, les seuls \'el\'ements de $\PSL$ fixant $\geod$ sont les conjugu\'ees des puissances de $A$. 
Une orbite modulaire p\'eriodique est donc naturellement associ\'ee \`a une classe de conjugaison dans $\SLZ$ sous l'action de $\SLZ$, \`a savoir la classe de~$A$ dans~$\SLZ$, ou de sa plus petite racine si $A$ n'est pas primitive. 
Or le corollaire~\ref{C:ConjLyndon} fournit une bijection entre les mots de Lyndon et les classes de conjugaison de matrices primitives hyperboliques.
\end{likeproof}

\begin{corollary}
\label{T:LyndonIndexeModulaire}
Les orbites p\'eriodiques du flot modulaire sont index\'ees par les mots de Lyndon.
\end{corollary}

On a maintenant obtenu de mani\`ere naturelle des n{\oe}uds sur le fibr\'e unitaire tangent \`a la surface modulaire. 
Or celui-ci se trouve \^etre isomorphe \`a la sph\`ere~$\Sph^3$ priv\'ee d'un n\oe ud de tr\`efle, via aux coordonn\'ees suivantes, introduites par Gauss, pour les r\'eseaux de~$\reels^2$.
Identifions $\reels^2$ avec $\Compl$, et posons
$$g_2(\Lambda)=60\sum_{z\in\Lambda\setminus\{0\}}z^{-4}
\mbox{\quad et \quad} g_3(\Lambda)=140\sum_{z\in\Lambda\setminus\{0\}}z^{-6}.$$
Les sommes $g_2$ et $g_3$ convergent et il se trouve que l'application $\Lambda\mapsto(g_2(\Lambda), g_3(\Lambda))$ est injective, c'est-\`a-dire  que le r\'eseau~$\Lambda$ est enti\`erement d\'{e}termin\'e par les valeurs des fonctions $g_2$ et $g_3$.
De plus, la non-d\'{e}g\'{e}n\'{e}rescence de~$\Lambda$ correspond \`a la condition $g_2^3-  27g_3^2 \neq 0$. 
Enfin, le fait que $\Lambda$ est de covolume~$1$ impose une condition de renormalisation sur $g_2(\Lambda)$ et~$g_3(\Lambda)$. 
Quitte \`a changer de renormalisation, on peut supposer 
\begin{equation}
	\abs{g_2(\Lambda)}^2+\abs{g_3(\Lambda)}^2=1.
\end{equation} 
Cette derni\`ere condition permet d'identifier l'image de $(g_2, g_3)$ \`{a} une partie de la sph\`ere $\Sph^3$,  qui est le compl\'{e}mentaire d'un n{\oe}ud de tr\`efle, not\'e \trefoil, dans $\Sph^3$. 
On a ainsi un isomorphisme entre $\Sph^3\smallsetminus$\nobreak\trefoil et la vari\'{e}t\'{e} modulaire~$\PSLR/\PSL$, qui nous permet de plonger les n{\oe}uds modulaires de mani\`ere naturelle dans $\Sph^3$\footnote{Consulter~\cite[p. 83]{Tali, TheseDehornoy, A:Paulin} pour des preuves plus visuelles d\'etaillant la topologie de $\USmod$.}.

\begin{definition}
Un n\oe ud est dit {\it modulaire} s'il peut \^etre r\'ealis\'e comme orbite p\'eriodique  du flot modulaire sur l'espace des r\'eseaux de covolume~1, plong\'e  dans~$\Sph^3$ via l'application~$(g_2, g_3)$.
\end{definition}


\subsection{N\oe uds modulaires et n\oe uds de Lorenz}
\label{S:ModulairesLorenz}

\`A ce point, nous avons introduits deux familles de n\oe uds, les n\oe uds de Lorenz et les n\oe uds modulaires, dont les \'el\'ements sont naturellement en bijection avec les mots de Lyndon. 
Le r\'esultat de Ghys affirme que ces familles co\"\i  ncident.

\begin{theorem}[\cite{A:Ghys3, A:Ghys4}]
	\label{T:Ghys}
	Plongeons $\USmod$ dans $\Sph^3$ \`a l'aide de l'application~$(g_2,g_3)\circ \resmod^{-1}$. 
	Alors, pour tout mot de Lyndon $w$, le n\oe ud modulaire et le n\oe ud de Lorenz associ\'es \`a~$w$ co\"\i ncident. 
\end{theorem}

\begin{likeproof}[Sch\'ema de la d\'emonstration]
L'id\'ee principale est de d\'eformer contin\^ument la m\'etrique sur $\surfmod$ et $\USmod$ pour amener tout le flot g\'eod\'esique au voisinage du segment reliant les images des points $i$ et $j$ de $\Hy$. 
En contractant la direction stable du flot, les orbites viennent alors s'accumuler sur un patron, qui se trouve \^etre le patron g\'eom\'etrique de Lorenz. 
\end{likeproof}

Pour plus de d\'etails, nous renvoyons \`a~\cite{A:Ghys3}, et \`a~\cite{A:Ghys4} o\`u l'on trouvera des repr\'esentations anim\'ees tr\`es spectaculaires. 

\begin{example}
Revenons au cas des id\'{e}aux de $\relatifs[2\sqrt6]$. 
On a vu qu'il y a trois  classes d'id\'{e}aux, dont trois repr\'{e}sentants sont  $\lid 1, 5+2\sqrt6 \rid, \lid 2, 4+2\sqrt6 \rid$ et $\lid 3, 4+2\sqrt6 \rid$. 
Les trois n{\oe}uds modulaires associ\'{e}s sont les n{\oe}uds correspondants aux mots $X^8Y, X^4Y^2$ et $X^2YXY$. 
Notre \'etude des n\oe uds de Lorenz et la correspondance de Ghys montrent que ceux-ci sont respectivement un n\oe ud trivial, un n\oe ud trivial, et un n\oe ud de tr\`efle.
\end{example}


\subsection{N\oe uds de Lorenz et classes d'id\'eaux}
\label{S:GroupeClasses}

On a maintenant une application naturelle associant \`{a} une classe d'id\'{e}aux dans un ordre d'un corps quadratique une unique classe de mots en $X$ et $Y$, d\'efinis \`a permutation circulaire des lettres et \`a \'echange des caract\`eres $X$ et $Y$ pr\`es. 
Or on ne change pas un n\oe ud de Lorenz en \'echangeant tous les  caract\`eres $\xx$ et $\yy$ dans le mot de Lyndon associ\'e\: cela correspond \`a une rotation autour de l'axe passant par le centre du patron et \'echangeant les deux points critiques, ou encore \`a une rotation de $180^\circ$ de la tresse de Lorenz autour d'un axe vertical. 
On a donc une bijection entre classes d'id\'eaux dans certains corps quadratiques et n\oe uds de Lorenz, modulo sym\'etrie du patron.
La question est alors de savoir si l'existence de cette bijection a des cons\'equences int\'eressantes. 

Celles-ci pourraient venir de  la structure de groupe qui existe naturellement sur les classes d'id\'eaux dans un corps quadratique.
Pr\'{e}cis\'{e}ment, soit $\rationnels[\alpha]$ un corps quadratique, $\grandordre$ son anneau des entiers, $\ordre$~un ordre d'indice $n$ dans $\grandordre$, et $\Ideal_1$ et $\Ideal_2$ deux id\'{e}aux de~$\ordre$. 
Alors la classe de l'id\'{e}al $\Ideal_1 \Ideal_2$ ne d\'{e}pend pas du choix de $\Ideal_1$ et $\Ideal_2$ dans leurs classes respectives. 
Par cons\'{e}quent, l'op\'{e}ration de multiplication est bien d\'{e}finie sur les classes d'id\'{e}aux. 
On v\'{e}rifie facilement que la classe des id\'{e}aux principaux est un \'{e}l\'{e}ment neutre pour cette op\'{e}ration, donc l'ensemble des classes d'id\'{e}aux forme un mono\"\i  de. 
De plus, toute classe contenant au moins un id\'{e}al $\Ideal$ dont la norme $N(\Ideal)$ est premi\`ere avec l'indice $n$ admet un inverse pour la multiplication, c'est-\`{a}-dire qu'il existe $\Ideal'$ tel que $\Ideal\Ideal'$ est un id\'{e}al principal, voir \cite[p. 122]{A:Cohn2}. 
L'ensemble de ces classes est stable par produit, par cons\'{e}quent l'ensemble des classes contenant au moins un id\'{e}al $\Ideal$ dont la norme est premi\`ere avec $n$ forme un groupe, appel\'{e} {\it groupe des classes} de l'ordre~$\ordre$. 

\begin{example}
\label{X:GroupeNoeuds}
Revenons une derni\`ere fois aux id\'{e}aux de $\relatifs[2\sqrt6]$. 
L'id\'{e}al $\lid 1, 5+2\sqrt6 \rid$ est principal, c'est donc un repr\'{e}sentant de la classe des id\'{e}aux principaux.
Pour ce qui est de $\lid 2, 4+2\sqrt6 \rid$, on v\'{e}rifie facilement que tout \'{e}l\'{e}ment de sa classe a une norme qui est multiple de 2. 
Or l'indice de~$\relatifs[2\sqrt6]$ dans $\relatifs[\!\sqrt6]$ est 2, donc la classe de $\lid 2, 4+2\sqrt6 \rid$ n'appartient pas au groupe des classes de~$\relatifs[2\sqrt6]$. 
De fait, on v\'{e}rifie l'\'egalit\'e $\lid 2, 4+2\sqrt6\rid^2=\lid 2\rid \lid 2, 4+2\sqrt6 \rid$, donc, si cette classe \'{e}tait dans le groupe des classes, elle en serait l'\'{e}l\'{e}ment neutre. 
Or cette position est d\'{e}j\`{a} occup\'{e}e par la classe des id\'{e}aux principaux.
Quant \`a $\lid 3, 4+2\sqrt6 \rid$, sa norme est 3, qui est premier \`a $2$, donc sa classe appartient au groupe des classes. 
Le groupe des classes de $\relatifs[2\sqrt6]$ est donc $\relatifs/2\relatifs$.
Notons que, dans la correspondance de Ghys, le n{\oe}ud associ\'{e} \`a l'\'el\'ement neutre est le n{\oe}ud trivial et le n\oe ud associ\'e \`a l'autre \'el\'ement est le n\oe ud de tr\`efle.

Pour un autre exemple, consid\'erons les matrices de trace 22 et de d\'{e}terminant 1. 
L'ordre associ\'{e} est~$\relatifs[2\sqrt30]$, qui compte six classes d'id\'{e}aux, correspondant aux matrices $X^{20}Y$, $X^{10}Y^2$, $X^5Y^4$, $X^6YXY$, $X^4YX^2Y$ et $X^2YX^2Y^2$. 
Les n{\oe}uds associ\'{e}s sont  le n{\oe}ud trivial pour les trois premi\`eres classes, le n{\oe}ud de tr\`efle pour les deux suivantes, et le n{\oe}ud torique $T(5,2)$ pour la derni\`ere.
Le groupe des classes de $\relatifs[2\sqrt30]$ a quatre \'el\'ements, car les classes correspondant aux mots $X^{10}Y^2$ et $ X^4YX^2Y$ ne sont pas inversibles. 
Ce groupe  est le groupe de Klein~$(\relatifs/2\relatifs)^2$. 
Vis-\`a-vis de la correspondance de Ghys, l'\'{e}l\'{e}ment neutre a pour n{\oe}ud associ\'{e} le n{\oe}ud trivial, et les trois autres \'{e}l\'{e}ments du groupe sont associ\'{e}s respectivement au n{\oe}ud trivial, au n{\oe}ud de tr\`efle et au n{\oe}ud~$T(5,2)$.
\end{example}


\subsection{Deux r\'esultats nouveaux}
\label{S:Nouveaux}

Comme on a associ\'{e} \`{a} chaque classe d'id\'eaux un n{\oe}ud de mani\`ere canonique, la question se pose naturellement de savoir si la multiplication du groupe  des classes d'id\'eaux a un pendant du c\^ot\'e des n{\oe}uds, c'est-\`a-dire si on peut d\'efinir directement une op\'eration de multiplication des n{\oe}uds et, par exemple, mieux comprendre les observations de l'exemple~\ref{X:GroupeNoeuds}. 
Les seules r\'eponses connues \`a ce jour, tr\`es partielles, sont les r\'esultats suivants, qui semblent nouveaux. 

Le premier r\'esultat d\'ecrit compl\`etement le sous-groupe du groupe des  classes d'id\'eaux form\'e par les classes associ\'ees au n{\oe}ud trivial. 
On note $\wedge$ l'op\'eration de pgcd de deux entiers et $\vee$ leur ppcm. 

\begin{theorem}
\label{T:groupetrivial}
Soit $t>2$ un entier, $\alpha$ une racine du polyn\^ome~$x^2-tx+1$, et $G$ le groupe des classes associ\'e \`a l'ordre~$\relatifs[\alpha]$. 
L'ensemble des \'el\'ements de~$G$ associ\'es au n\oe ud trivial forme un sous-groupe~$G_0$ de~$G$ qui est isomorphe \`a~$(\relatifs/2\relatifs)^{d-1}$ o\`u $d$ est le nombre de diviseurs premiers de~$t-2$. 
Les matrices correspondant aux \'el\'ements de~$G_0$ sont de la forme $X^mY^{m'}$ avec $mm'=t-2$ et~$m\wedge m'=1$. 
Exprim\'ee en termes de matrices, la multiplication de~$G_0$ est
donn\'ee par la formule
$$X^{m_1}Y^{(t-2)/m_1}\ \cdot\  X^{m_2}Y^{(t-2)/m_2}
= X^{m_3}Y^{(t-2)/m_3},\mbox{ \ avec } 
m_3 = \frac{m_1\vee m_2}{m_1\wedge m_2}.$$
\end{theorem}

\begin{proof}
Le n\oe ud trivial \'etant le seul n\oe ud obtenu comme cl\^oture d'une tresse \`a un brin, c'est le seul n\oe ud de Lorenz de \trip~1. 
Par cons\'equent, les orbites triviales du flot de Lorenz sont exactement celles qui sont associ\'ees aux mots de la forme $X^mY^{m'}$. 

Comme inverser dans le groupe des classes revient \`a transposer la matrice, les classes correspondant aux mots $X^mY^{m'}$ et $X^{m'}Y^{m}$  sont inverses l'une de l'autre. 
Or ces classes co\"\i  ncident, puisqu'au niveau des id\'eaux l'\'echange des caract\`eres $X$ et $Y$ et la permutation circulaire des lettres ne modifient pas la classe. 
Par cons\'equent, ces classes sont d'ordre 2 dans le groupe des classes~$G$.

L'ordre associ\'e aux matrices de trace~$t$ est~$\relatifs[\alpha]$,  o\`u $\alpha=\frac{-t+\sqrt{t^2-4}}{2}$, et l'id\'eal associ\'e \`a la matrice~$X^mY^{m'}$, c'est-\`a-dire la matrice~$\mama{1}{m}{m'}{1+mm'}$, est~$\lid m,\alpha-1 \rid$. 
D'apr\`es~\cite[p. 220]{A:Cohen}, la forme quadratique associ\'ee est $mx^2+mm'xy+m'y^2$, et, d'apr\`es la proposition 5.2.5 de~\cite{A:Cohen}, la condition pour qu'un id\'eal soit inversible est que les coefficients de la forme quadratique associ\'ee soient premiers entre eux, soit, dans le cas pr\'esent, 
$m\wedge m'=1$. 
Les n\oe uds triviaux apparaissant dans le groupe des classes correspondent donc aux mots de la forme $X^mY^{m'}$, avec $mm'=t-2$ et $m\wedge m'=1$. 

Il reste \`a v\'erifier la stabilit\'e par produit des \'el\'ements associ\'es \`a des n\oe uds triviaux, et la forme explicite du produit annonc\'ee. 
Calculant au niveau des id\'eaux, on~trouve
\begin{align*}
\lid m_1,\alpha-1 \rid
&\cdot \lid m_2,\alpha-1\rid \\
   &= \lid m_1m_2, m_1(\alpha-1),  m_2(\alpha-1), \alpha^2-2\alpha+1\rid
\\
   &= \lid (m_1\wedge m_2) (m_1\vee m_2), (m_1\wedge m_2)(\alpha-1),
(t-2)\alpha \rid \\
   &= (m_1\wedge m_2) \lid m_1\vee m_2, \alpha-1, \frac{t-2}{m_1\wedge
m_2}\alpha \rid.
 \end{align*}
Les diviseurs de $t-2$ se divisent en quatre cat\'egories\: 
ceux qui ne divisent ni $m_1$ ni $m_2$ dont nous noterons le 
pgcd~$p_{--}$,  ceux qui divisent $m_1$ et pas $m_2$
dont nous noterons le $\pgcd$~$p_{+-}$,  ceux qui divisent $m_2$ et pas
$m_1$ dont nous noterons le $\pgcd$~$p_{-+}$, et ceux qui divisent
$m_1$ et $m_2$ dont nous noterons le $\pgcd$ $p_{++}$. 
On a alors
$$m_1\vee m_2=p_{-+}p_{+-}p_{++} \mbox{\quad et \quad} 
	\frac{t-2}{m_1\wedge m_2} = p_{--}p_{-+}p_{+-},$$
d'o\`u 
$(m_1\vee m_2)\wedge {\displaystyle\frac{t-2}{m_1\wedge m_2}}
= p_{-+}p_{+-} = {\displaystyle\frac{m_1\wedge m_2}{m_1\vee m_2}}$. 
On en d\'eduit 
 $$(m_1\wedge m_2) \lid m_1\vee m_2, \alpha-1, \frac{t-2}{m_1\wedge m_2}\alpha \rid
	= (m_1\wedge m_2) \lid \frac{m_1\wedge m_2}{m_1\vee m_2}, \alpha-1\rid.$$
Ce dernier id\'eal est dans la m\^eme classe que $\lid {\displaystyle\frac{m_1\wedge m_2}{m_1\vee m_2}}, \alpha-1 \rid,$ dont la matrice id\'eale est 
$$X^{\frac{m_1\vee m_2}
	{m_1\wedge m_2}}Y^{(t-2)\big/\frac{m_1\vee m_2}{m_1\wedge m_2}}.$$
\end{proof}

Notre seconde observation est un r\'esultat de sym\'etrie tr\`es simple.
Si $w$ est un mot (de Lyndon), on appelle \emph{miroir} de~$w$ le mot obtenu en renversant l'ordre des lettres de~$w$.

\begin{proposition}
\label{T:inverse}
Pour tout mot de Lyndon~$w$, les n\oe uds  de Lorenz correspondants \`a $w$ et \`a son miroir sont les m\^emes. 
\end{proposition} 

\begin{proof}
Transpos\'e en termes de n\oe uds modulaires, l'\'enonc\'e devient\: les n\oe uds sur le fibr\'e unitaire tangent \`a la surface modulaire~$\USmod$ qui correspondent au parcours dans un sens et dans l'autre d'une g\'eod\'esique p\'eriodique de la surface modulaire~$\surfmod$ sont les m\^emes. 
Or, vis-\`a-vis de la m\'etrique riemannienne~$g$ introduite au d\'ebut de la section, l'involution de~$\USmod$ qui associe au vecteur tangent $v$ de $\USmodx$ le vecteur $-v$ de $\USmodx$ est une isom\'etrie de $\USmod$\footnote{En fait, c'est m\^eme une rotation d'angle~$\pi$ par rapport \`a l'axe constitu\'e des vecteurs tangents au point $i$.}, qui pr\'eserve l'orientation. 
Par cons\'equent, elle envoie tout lacet sur un lacet isotope, et les n{\oe}uds associ\'es co\"\i ncident.
\end{proof}

Ce r\'esultat montre l'int\'er\^et de la correspondance de Ghys et du lien entre n{\oe}uds de Lorenz et n{\oe}uds modulaires. 
Sur le patron de Lorenz, une seule involution naturelle appara\^it, \`a savoir la sym\'etrie par rapport \`a l'axe du patron, qui correspond \`a la sym\'etrie par rapport \`a la premi\`ere diagonale sur les diagrammes de Young associ\'es. 
La proposition~\ref{T:inverse} r\'ev\`ele une autre sym\'etrie sur les orbites du flot de Lorenz, naturelle dans le langage du fibr\'e tangent~$\USmod$ \`a la surface modulaire, mais cach\'ee sur le patron de Lorenz. 
Par exemple, elle montre que les orbites associ\'ees aux mots de Lyndon $\xx\xx\yy\yy\xx\yy\yy\yy$ et $\xx\xx\yy\yy\yy\xx\yy\yy$, correspondant aux tableaux compl\'et\'es
\begin{picture}(9,0)
	\put(0,-1){\includegraphics[scale=.35]{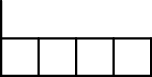}}
\end{picture}
et 
\begin{picture}(9.5,0)
	\put(0,-1){\includegraphics[scale=.35]{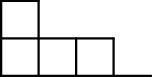}}
\end{picture}, 
sont isotopes.

\begin{example}
\label{X:groupetrivial}
Comme derni\`ere illustration, consid\'erons les n\oe uds associ\'es \`a des matrices de trace~$40$. 
L'ordre associ\'e est~$\relatifs[-20+3\sqrt{11}]$. 
Un logiciel de calcul alg\'ebrique, comme par exemple PARI/GP~\footnote{\tt{http:\!\!//pari.math.u-bordeaux.fr/}}, fournit la structure du groupe des classes associ\'e, ici $\relatifs/2\relatifs \times \relatifs/4\relatifs$, et des repr\'esentants des classes associ\'ees \`a chaque \'el\'ement. 
En termes de matrices, les huit \'el\'ements du groupe des classes d'id\'eaux de l'ordre~$\relatifs[-20+3\sqrt{11}]$ sont les suivants.

\begin{center}
\vspace{1mm}
\begin{tabular}{c|ccccl}
\vrule width0pt height12pt depth 8pt
$G = \relatifs/2\relatifs \times \relatifs/4\relatifs$
&\centering{$(\cdot, \bar0)$}
&\centering{$(\cdot, \bar1)$}
&\centering{$(\cdot, \bar2)$}
&\centering{$(\cdot, \bar3)$} 
&\\
\hline
\vrule width0pt height14pt depth 5pt
$(\bar0, \cdot)$
&\centering{$X^{38}Y$} 
&\centering{$X^3Y^3X^2Y$} 
&\centering{$XYXY^{12}$} 
&\centering{$X^3YX^2Y^3$}
&\\  
\vrule width0pt height12pt depth 5pt
$(\bar1, \cdot)$
&\centering{$X^{19}Y^2$}
&\centering{$X^2Y^7XY$} 
&\centering{$X^5YX^4Y$} 
&\centering{$X^2YXY^7$}                     
\end{tabular}
\vspace{2mm}
\end{center}
Traduit en termes de n\oe uds de Lorenz via la correspondance de Ghys, 
ce tableau devient\:

\begin{center}
\vspace{2mm}
\begin{tabular}{c|ccccl}
\vrule width0pt height12pt depth 9pt
$G = \relatifs/2\relatifs \times \relatifs/4\relatifs$
&\centering{$(\cdot, \bar0)$}
&\centering{$(\cdot, \bar1)$}
&\centering{$(\cdot, \bar2)$}
&\centering{$(\cdot, \bar3)$} 
&\\
\hline
\vbox to 11mm{\hsize=10mm$(\bar0, \cdot)$}
&\centering{\includegraphics[scale=0.7]{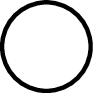}} 
&\centering{\includegraphics[scale=0.6]{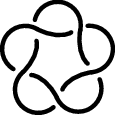}} 
&\centering{\includegraphics[scale=0.6]{K2.eps}} 
&\centering{\includegraphics[scale=0.6]{K1.eps}}
&\\  
\vbox to 11mm{\hsize=10mm$(\bar1, \cdot)$}
&\centering{\includegraphics[scale=0.7]{K0.eps}}
&\centering{\includegraphics[scale=0.6]{K2.eps}} 
&\centering{\includegraphics[scale=0.6]{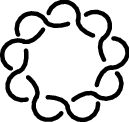}} 
&\centering{\includegraphics[scale=0.6]{K2.eps}}                     
\end{tabular}
\vspace{2mm}
\end{center}
On peut y voir illustr\'es les r\'esultats des th\'eor\`emes~\ref{T:groupetrivial}
et~\ref{T:inverse}. D'abord, comme 40 admet deux diviseurs premiers
distincts, les n\oe uds triviaux forment un sous-groupe isomorphe \`a
$\relatifs/2\relatifs$, qui correspond \`a la colonne~$(\cdot, \bar0)$ du tableau. Par
ailleurs, les colonnes~$(\cdot, \bar1)$ et~$(\cdot, \bar3)$
correspondent \`a des \'el\'ements oppos\'es dans le groupe, et on voit que, sur
chaque ligne, les n{\oe}uds correspondants sont en effet les m\^emes.
\end{example}

\section{Quelques questions ouvertes}

Nous terminerons ce tour d'horizon des n{\oe}uds de Lorenz par quelques-unes des multiples questions ouvertes les concernant.

\subsection{Codage par les mots de Lyndon et les diagrammes de Young}

On a vu dans la section~\ref{S:Combinatoire} comment tous les n{\oe}uds de Lorenz peuvent \^etre cod\'es par des mots de Lyndon, ou, de fa\c{c}on presque \'equivalente, par des diagrammes de Young. 
Ces codages ne sont pas injectifs\: m\^eme minimaux, deux mots de Lyndon distincts peuvent 
coder le m\^eme n{\oe}ud de Lorenz et, de m\^eme, deux diagrammes de Young (non compl\'et\'es) distincts peuvent coder le m\^eme n{\oe}ud de Lorenz. 
Reconna\^itre quand deux mots de Lyndon, ou deux diagrammes de Young, codent le m\^eme n\oe ud de Lorenz est une question pour le moment difficile. 
On a vu dans le corollaire~\ref{C:GenreTableau} que le genre du n\oe ud d\'etermine le nombre de cases du diagramme de Young, mais bien s\^ur pas le diagramme complet. 
Par exemple, on v\'erifie facilement que la tresse de Birman-Williams associ\'ee au deux diagrammes de Young $[3,3,2]$ et $[4,4]$ est la m\^eme, et donc les n{\oe}uds de Lorenz associ\'es co\"\i ncident.

\begin{question}
\label{Q:Codage}
Existe-t-il un algorithme permettant de d\'eterminer quand deux diagrammes de Young codent le m\^eme n{\oe}ud de Lorenz\?
\end{question}

La section~\ref{S:B-K} et l'article~\cite{A:B-K} donnent une reformulation plus g\'en\'erale.

\begin{question}
\label{Q:CodageToriqueItere}
Existe-t-il un algorithme permettant de d\'eterminer quand deux $T$-n\oe uds sont isotopes\?
\end{question}

Une premi\`ere \'etape pourrait consister \`a d\'emontrer la propri\'et\'e suivante, qui a \'et\'e v\'erifi\'ee sur un grand nombre d'exemples.

\begin{conjecture}[\cite{A:B-K}]
Soit $b$ et $b'$ deux tresses de Birman-Williams associ\'ees \`a deux orbites p\'eriodiques isotopes du flot de Lorenz. 
Alors $b$ et $b'$ sont conjugu\'ees dans $\braid_t$, o\`u $t$ est le \trip\ commun de ces orbites.
\end{conjecture}

\subsection{N{\oe}uds satellites}

Avec le th\'eor\`eme~\ref{T:satellite}, on a obtenu une description compl\`ete de ceux des n{\oe}uds de Lorenz qui sont satellites d'un n{\oe}ud de Lorenz. 
Par contre, on n'a pas de caract\'erisation g\'en\'erale de tous les n{\oe}uds de Lorenz qui sont des n{\oe}uds  satellites. 
Il semblerait en fait que seule la situation du th\'eor\`eme~\ref{T:satellite} soit possible.

\begin{conjecture} [H.\,Morton, communication personnelle]
\label{C:Satellite}
Tout n\oe ud de Lorenz qui est un n\oe ud satellite est un c\^ablage sur un n\oe ud de Lorenz, selon le sch\'ema de la
proposition~\ref{T:cablage}.
\end{conjecture}

Un important r\'esultat de W.\,Thurston~\cite{A:Thurston} montre que tout n{\oe}ud qui n'est ni torique, ni satellite est hyperbolique. 
Notons qu'une solution positive \`a la conjecture~\ref{C:Satellite} permettrait de d\'eterminer pr\'ecis\'ement quels n{\oe}uds de Lorenz sont des n{\oe}uds satellites, et, de l\`a, par compl\'ement, quels n{\oe}uds
de Lorenz sont hyperboliques.

Dans la lign\'ee de la question~\ref{Q:Codage}, il est \'egalement naturel de chercher \`a reconna\^itre le type d'un n{\oe}ud de Lorenz \`a partir d'un diagramme de Young le codant.

\begin{question}
\label{Q:Type}
Existe-t-il un algorithme qui, \`a partir d'un diagramme de Young codant un n\oe ud de Lorenz $K$, d\'etermine le type du n\oe ud torique associ\'e si $K$ est torique, et la suite des cabl\^ages successifs si $K$ est satellite\?
\end{question}

\subsection{Invariants polynomiaux et volume hyperbolique}

Il semble que les polyn\^omes d'Alexander ou de Jones des n\oe uds de Lorenz soient particuli\`erement simples.
Les simulations num\'eriques de~\cite{A:Dehornoy} montrent que les coefficients du polyn\^ome d'Alexander d'un n{\oe}ud de Lorenz donn\'e comme cl\^oture d'une tresse, tout comme ceux de leurs polyn\^omes de Jones, sont tr\`es petits compar\'es \`a ceux de la cl\^oture d'une tresse al\'eatoire de m\^eme taille.

\begin{question}
La petitesse observ\'ee des coefficients des polyn\^omes d'Alexander, de Jones, et HOMFLY associ\'es aux n{\oe}uds de Lorenz est-elle un ph\'enom\`ene g\'en\'eral, et peut-on l'expliquer\? Comment se r\'epartissent les racines des polyn\^omes d'Alexander et de Jones dans le plan complexe\?
\end{question}

Une réponse partielle pour la répartition des racines du polynôme d'Alexander est donnée dans~\cite{TheseDehornoy}, mais rien n'est dit sur les coefficients.
Dans~\cite{A:B-K}, J.\,Birman et I.\,Kofman ont calcul\'e les volumes hyperboliques des  compl\'ementaires d'un certain nombre de n\oe uds de Lorenz hyperboliques. 
Comme dans le cas des coefficients des invariants polynomiaux, il appara\^it que ce volume prend des valeurs \'etonnamment petites.

\begin{question}
Peut-on borner sup\'erieurement le volume hyperbolique des n\oe uds de Lorenz qui sont hyperboliques en fonction de la taille d'une tresse de Birman--Williams les repr\'esentant\?
\end{question}

Par ailleurs, il se trouve que, parmi les cent n\oe uds hyperboliques de plus petit volume, plus de la moiti\'e sont des n\oe uds de Lorenz. Une explication de cette observation \'etonnante serait bienvenue.

\subsection{Correspondance de Ghys et n{\oe}uds modulaires}

Le lien entre les n\oe uds de Lorenz et l'arithm\'etique par le biais du groupe des classes d'id\'eaux associ\'e \`a un ordre sur un corps quadratique est fascinant, mais encore tr\`es incompl\`etement compris.
Les propositions~\ref{T:groupetrivial} et~\ref{T:inverse} semblent n'\^etre que les parties \'emerg\'ees d'un analogue du groupe des classes qui serait d\'efini sur les n\oe uds modulaires, donc, de fa\c{c}on \'equivalente, sur les n{\oe}uds de Lorenz.

\begin{question}
Peut-on d\'efinir directement au niveau des n{\oe}uds modulaires, ou des n{\oe}uds de Lorenz, une multiplication (d\'ependant de $\alpha$) qui induise la structure de groupe associ\'ee sur les n{\oe}uds correspondant au groupe de classes d'id\'eaux de~$\relatifs[\alpha]$\?
\end{question}


\bibliographystyle{amsplain}

\providecommand{\bysame}{\leavevmode\hbox to3em{\hrulefill}\thinspace}
\providecommand{\MR}{\relax\ifhmode\unskip\space\fi MR }
\providecommand{\MRhref}[2]{%
  \href{http://www.ams.org/mathscinet-getitem?mr=#1}{#2}
}
\providecommand{\href}[2]{#2}
\begin{thebibliography}{}

\end{thebibliography}


\begin{thebibliography}{30}

\bibitem{A:Alexander}{\sc J.\,W.\,Alexander}, {\it A lemma on systems of knotted curves}, {Proc. Amer. Sci. U.\,S.\,A.} {\bf 9} (1923), 93--95.

\bibitem{A:Bediant}{\sc R.\,Bedient}, {\it Classifying 3-trip Lorenz knots}, {Topology Appl.} {\bf 20} (1985), 89--96.

\bibitem{A:B-P}{\sc J.\,Berstel, D.\,Perrin}, {\it The origins of combinatorics on words}, European J. Combinatorics
{\bf 28} (2007), 996--1022.

\bibitem{A:B-K}{\sc J.\,S.\,Birman, I.\,Kofman}, {\it A new twist on Lorenz knots}, {\tt arXiv:0707.4331v3 [math.GT]}.

\bibitem{A:B-W}{\sc J.\,S.\,Birman, R.\,F.\,Williams}, {\it Knotted periodic orbits in dynamical systems--I: Lorenz's Equations} {Topology} {\bf 22} (1) (1983), 47--82. (erratum at {\tt http:\!\!//www.math.columbia.edu/~jb/bw-KPO-I-erratum.pdf})

\bibitem{A:B-Z}{\sc G.\,Burde, H.\,Zieschang}, {\it Knots}, de Gruyter Studies in Mathematics 5, Walter de Gruyter \& Co, 1985.

\bibitem{A:C-F-L}{\sc K.\,Chen, R.\,Fox, R.\,Lyndon}, {\it Free differential calculus IV: The quotient groups of the
lower central series}, {Annals of Math.} {\bf 58} (1958), 81--95.

\bibitem{A:Cohen}{\sc H.\,Cohen}, {\it A Course in computational algebraic number theory}, GTM 138, Springer Verlag, 1993.

\bibitem{A:Cohn}{\sc H.\,Cohn}, {\it A classical invitation to algebraic numbers and class fields}, Universitext, Springer Verlag, 1978.

\bibitem{A:Cohn2}{\sc H.\,Cohn}, {\it Advanced number theory}, Dover Publications, 1980.

\bibitem{A:Paulin}{\sc F.\,Dal`Bo, F.\,Paulin, G.\,Courtois}, {\it Sur la dynamique des groupes de matrices et applications arithm\'etiques}, Journ\'ees math\'ematiques X-UPS 2007.

\bibitem{A:atlas} {\sc Pi.\,Dehornoy}, {\it Atlas of Lorenz knots},
{\tt http:\!\!//www.eleves.ens.fr/home/dehornoy/atlaslorenz.txt}

\bibitem{A:Dehornoy} {\sc Pi.\,Dehornoy}, {\it N\oe uds de Lorenz},
M\'emoire de M2R, 

{\tt http:\!\!//www.eleves.ens.fr/home/dehornoy/maths/Lorenz4.pdf}

\bibitem{TheseDehornoy} {\sc Pi.\,Dehornoy}, {\it Invariants topologiques des orbites p\'eriodiques d'un champ de vecteurs}, th\`ese en pr\'eparation.

\bibitem{A:E-N}{\sc D.\,Eisenbud, W.\,Neuman}, {\it Three-dimensional link theory and invariants of plane curve singularities}, Ann. Math. Studies 110, Princeton University Press, 1985.

\bibitem{A:ElRifai}{\sc E.\,A.\,El\,Rifai}, {\it Positive braids and Lorenz links}, PhD. Thesis, Univ. Liverpool, 1988.

\bibitem{A:ElRifai2}{\sc E.\,A.\,El\,Rifai}, {\it Necessary and sufficient condition for Lorenz knots to be closed under satellite construction}, {Chaos, Solitons and Fractals} {\bf 10} (1) (1999), 137--146.

\bibitem{A:F-W}{\sc J.\,Franks, R.\,F.\,Williams}, {\it Braids and the Jones polynomial}, {Trans. Amer. Math. Soc.} {\bf 303} (1) (1987), 97--108.

\bibitem{A:Gabai1}{\sc D.\,Gabai}, {\it The Murasugi sum is a natural geometric operation}, {Cont. Math.} {\bf 20} (1983), 131--143.

\bibitem{A:Gabai2}{\sc D.\,Gabai}, {\it Detecting fibred links in $\Sph^3$}, {Comment. Math. Helv.} {\bf 61} (1986), 519--555.

\bibitem{A:Garside}{\sc A.\,Garside}, {\it The braid group and other groups}, {Quart. J. Math. Oxford} {\bf 20} (1969), 235--254.

\bibitem{A:Ghrist}{\sc R.\,W.\,Ghrist}, {\it Branched two-manifolds supporting all links}, {Topology} {\bf 36}(2) (1997), 423--488.

\bibitem{A:G-H-S}{\sc R.\,W.\,Ghrist, Ph.\,J.\,Holmes, M.\,C.\,Sullivan}, {\it Knots and links in three-dimensional flows}, Lect. Notes Math. 1654, Springer Verlag, 1997.

\bibitem{A:Ghys3}{\sc \'E.\,Ghys}, {\it Knots and dynamics}, {Proc. of the International Congress of Mathematicians} {\bf I}, Madrid, 2006.

\bibitem{A:Ghys4}{\sc \'E.\,Ghys, J.\,Leys}, {\it Lorenz and modular flow: a visual introduction}, {Feature Column Amer. Math. Soc.} (2006).

\bibitem{A:G-P}{\sc N.\,D.\,Gilbert, T.\,Porter}, {\it Knots and Surfaces}, Oxford U.P. (1994).

\bibitem{A:Guckenheimer}{\sc J.\,Guckenheimer}, {\it A strange, strange attractor}, {The Hopf Bifurcation Theorem and its Applications}, Springer-Verlag (1976), 368--381. 

\bibitem{A:G-W}{\sc R.\,Guckenheimer, R.\,F.\,Williams}, {\it Structural stability of Lorenz attractors}, {Publ. Math. IHES} {\bf 50} (1979), 59--72.

\bibitem{A:Hadamard}{\sc J. Hadamard}, {\it Les surfaces \`a courbures oppos\'ees et leurs lignes g\'eod\'esiques}, {J.~Math. Pures Appl.} {\bf 4} (1898), 27--74.

\bibitem{A:Katok}{\sc S.\,Katok}, {\it Fuchsian Groups}, University of Chicago Press, 1992.

\bibitem{A:Lickorish}{\sc W.\,B.\,R.\,Lickorish}, {\it An introduction to knot theory}, GTM 175, Springer Verlag, 1997.

\bibitem{A:Lorenz}{\sc E.\,N.\,Lorenz}, {\it Deterministic nonperiodic flow}, {J.~Atmospheric Sci.} {\bf 20} (1963), 130--141.

\bibitem{A:Lorenz2}{\sc E.\,N.\,Lorenz}, {\it Predictability : does the flap of a butterfly's wings in Brazil set off a tornado in Texas\?},  139th Meeting of the AAAS (1972).

\bibitem{A:Lyndon}{\sc R.\,Lyndon}, {\it On Burnside problem I}, {Trans. AMS} {\bf 77} (1954), 202--215.

\bibitem{Meleshuk}{\sc V.\,Meleshuk}, {\it Embedding templates in flows}, PhD. Thesis, Northwestern Univ., 2002.

\bibitem{A:Milnor}{\sc J.\,Milnor}, {\it Singular points of complex hypersurfaces}, Ann. Math. Studies 61, Princeton University Press (1965).

\bibitem{A:Morton}{\sc H.\,R.\,Morton}, {\it Seifert circles and knot polynomials}, 
Math. Proc. Camb. Phil. Soc. {\bf 99} (1986), 107--109. 

\bibitem{A:Murasugi}{\sc K.\,Murasugi}, {\it On a certain subgroup of the group of an alternating link}, {Amer.~J. Math.} {\bf 85} (1963), 544--550.

\bibitem{Tali}{\sc T.\,Pinsky}, {\it Templates for geodesic flows}, phD. Thesis (2011).

\bibitem{A:PoincarŽ}{\sc H.\,Poincar\'e}, {\it Sur le probl\`eme des trois corps et les \'equations de la dynamique}, {Acta. Math.} {\bf 13}  (1890), 1--270.
 
\bibitem{A:Reut}{\sc C.\,Reutenauer}, {\it Free Lie algebras}, Clarendon Press, Oxford, 1993.

\bibitem{A:Rolfsen}{\sc D.\,Rolfsen}, {\it Knots and links}, Publish or Perish, 1976.

\bibitem{A:Rudolph}{\sc L.\,Rudolph}, {\it Non-trivial positive braids have positive signature}, Topology {\bf 21} (1983), 325--327.

\bibitem{A:Serre}{\sc J.-P.\,Serre}, {\it Cours d'arithm\'etique}, Presses universitaires de France, 1970.

\bibitem{A:Sparrow}{\sc C.\,Sparrow}, {\it The Lorenz equations\: bifurcation, chaos and strange attractors}, Applied Mathematical Science 41, Springer Verlag, 1982.

\bibitem{A:Stallings}{\sc J.\,Stallings}, {\it Constructions of fibered knots and links}, {Symp. in Pure math. Am. Math. Soc.} Part 2, (1978), 55--59.

\bibitem{A:Sullivan1}{\sc M.\,C.\,Sullivan}, {\it The prime decomposition of knotted periodic orbits in Lorenz-like templates}, J. Knot Thy. and Ram. {\bf 3}(1) (1994), 83--120.

\bibitem{A:Thurston}{\sc W.\,Thurston}, {\it The geometry and topology of 3-manifolds}, Princeton lecture notes (1978-1981).

\bibitem{A:Tucker1}{\sc W.\,Tucker}, {\it The Lorenz attractor exists}, PhD. Thesis, Univ. Uppsala, 1998.

\bibitem{A:Tucker2}{\sc W.\,Tucker}, {\it A rigorous ODE solver and Smale's $14$th problem}, {Found. Comput. Math.} {\bf 2} (2002), 53--117.

\bibitem{A:Williams2}{\sc R.\,F.\,Williams}, {\it The structure of Lorenz attractors}, {Publ. Math. IHES} {\bf 50} (1979), 73--99.

\bibitem{A:Williams}{\sc R.\,F.\,Williams}, {\it Lorenz knots are prime}, {Ergod. Th. Dynam. Sys.} {\bf 4} (1983), 147--163.

\end{thebibliography}

\vspace{2cm}

\hspace{10cm}Lyon, le 27 novembre 2009

\end{document}